\documentclass[preprint,12pt]{elsarticle}

\usepackage{amssymb}
\usepackage{indentfirst}
\usepackage{graphicx}
\usepackage{epstopdf}
\usepackage{epsfig}
\usepackage{overpic}
\usepackage{array}
\usepackage{color}
\usepackage{multirow}
\usepackage{hyperref}
\usepackage{float}
\usepackage{subfigure}
\usepackage{amssymb}
\usepackage{amsmath,amssymb, amsthm}
\usepackage[mathscr]{euscript}
\usepackage{latexsym,bm}
\usepackage{booktabs}
\usepackage{bbm}
\usepackage{diagbox}
 \usepackage[all]{xy}
\usepackage{diagbox}
\usepackage{algorithm}  
\usepackage{algorithmicx} 
\usepackage{algpseudocode} 
\usepackage{mathtools}

\newcommand\cf{\color{red}}

\def\mi{\mathbbm{i}}
\def\me{\mathbbm{e}}

\def\bx{\bm{x}}

\def\bv{\bm{v}}

\def\D{\textup{d}}

\makeatletter
\@addtoreset{equation}{section}
\makeatother

\begin{document}

\newsavebox{\tablebox}

\begin{frontmatter}



\title{Distributed local spline simulator for wave propagation}

  \author{Xu Guo\corref{cofirst}\fnref{label1}}
 \author{Yaomeng Li\corref{cofirst}\fnref{label1}}
\author[label2]{Yunfeng Xiong\corref{cor}}
\ead{yfxiong@bnu.edu.cn}

\affiliation[label1]{organization={Geotechnical and Structural Engineering Center, Shandong University},
            city={Jinan},
            postcode={250061}, 
            state={Shandong},
            country={China}}

\affiliation[label2]{organization={School of Mathematical Sciences, Beijing Normal University},
            postcode={100871}, 
            state={Beijing},
            country={China}}
 
\cortext[cofirst]{These authors contribute equally to this paper}           
\cortext[cor]{To whom correspondence}

\begin{abstract}
Numerical simulation of wave propagation in elastic media faces the challenges arising from increasing demand of high resolution in modern 3-D imaging applications, which requires a balance between efficiency and accuracy in addition to being friendly to the distributed high-performance computing environment. In this paper, we propose a distributed local spline simulator (LOSS) for solving the wave equation. LOSS uses patched cubic B-splines to represent the wavefields and attains an accurate evaluation of spatial derivatives with linear complexity.  In order to link the adjacent patches, a perfectly matched boundary condition is introduced to give a closure of local spline coefficients. Owing to the rapid decay property of the local wavelets in dual space, it can recover the global spline as accurately as possible only at the cost of local communications among adjacent neighbors. Several typical numerical examples, including 2-D acoustic wave equation and $P$- and $S$- wave propagation in 3-D homogenous or heterogenous media, are provided to validate its convergence, accuracy and parallel scalability. 
\end{abstract}


%
%
%
%
%
%
%


\end{frontmatter}


\section{Introduction}

The wave propagation in elastic media plays an essential role in the field of geological imaging techniques \cite{VirieuxOperto2009,Tromp2020,MirzanejadTranWang2022} and new-trend neuroimaging \cite{GuaschAgudoTangNachevWarner2020}. Nowadays, a large collection of efficient numerical techniques is available for both forward and inverse wavefield modelings, including the widely used staggered-grid finite difference method \cite{Graves1996,KristekMoczo2003}, the optimal difference method \cite{Liu2014}, the pseudo-spectral method \cite{Carcione1994,bk:Fornberg1998,Carcione2009,XiongGuo2022}. In spite of their great success and extensive applications, these standard techniques face new challenges by tremendous memory demand and significant computational cost especially for 3-D problems, arising from the increasing need for high resolution in modern imaging applications, e.g., the teleseismic datasets for waveform inversion and deep lithospheric structures \cite{Tromp2020,LiuYangDongLiuZheng2017} or the sub-millimetre-resolution in brain and surrounding tissue \cite{GuaschAgudoTangNachevWarner2020}. Thus, it urgently calls for discretized methods to account for sharp variations of solutions induced by material discontinuities accurately \cite{KomatitschTromp1999} and to be friendly to large-scale high-performance computing environment \cite{BaoBielakGhattasKallivokasOHallaronShewchukXu1998}.

In recent years, the spectral element method (SEM) and discontinuous Galerkin method (DG) \cite{ChungEngquist2006,StanglmeierNguyenPeraireCockburn2016} have gained an increasing attention \cite{LiuYangDongLiuZheng2017,KomatitschTromp1999,KomatitschTsuboTromp2002,KomatitschTsuboTromp2005,TrinhBrossierMetivierTavardVirieux2019} as they take advantage of both flexibility of the finite element method (FEM) in resolving multiscale phenomena \cite{BaoBielakGhattasKallivokasOHallaronShewchukXu1998} and high accuracy of the spectral method. A significant merit of SEM and DG over FEM is that the mass matrix is exactly diagonal by construction, which drastically simplifies the implementation and the temporal integration \cite{KomatitschTromp1999}. Moreover, the assembly of the stiffness matrix can be performed in an element-by-element manner, thereby greatly facilitating the parallelization in a distributed computing environment \cite{LiuYangDongLiuZheng2017}. But the evaluation of the stiffness matrix at the elemental level has a relatively high computational cost due to the matrix multiplications involved, e.g., the complexity scales as $O(n_l^4)$ for SEM in three dimensions with $n_l$ the polynomial degree used to present the functions in each direction \cite{KomatitschTromp1999}. It may somehow pose a limitation on $n_l$ to achieve a trade-off in accuracy and efficiency.

As forward wavefield modelings have to be simulated thousands of times in waveform-fitting imaging, reducing the computational complexity of numerical solvers is always a central issue especially for 3-D problems \cite{VirieuxOperto2009}. Thus it is natural to seek a local polynomial basis that can calculate the spatial derivatives of wavefields accurately with relatively low computational cost. To achieve this, we propose a distributed local spline simulator (LOSS) for solving the wave propagation, using the local cubic B-spline wavelet as the basis  \cite{bk:Chui1992}. The name LOSS  comes from  the semi-Lagrangian methods in computational fluid dynamics \cite{StaniforthCote1991,MalevskyThomas1997} and kinetic theory \cite{SonnendruckerRocheBertrand1999,CrouseillesLatuSonnendrucker2009,CrouseillesMehrenbergerSonnendrucker2010}, where the cubic spline has been ubiquitously applied for interpolating the advection and is believed to strike the best balance between accuracy and cost \cite{CrouseillesLatuSonnendrucker2009}. The cubic spline achieves spatial fourth-order convergence \cite{Bermejo1990} and its construction can be realized by the standard sweeping method, where the complexity scales linearly with respect to mesh size \cite{bk:Boor2001}. For these reasons, the most recent semi-Lagrangian methods are even capable to resolve the kinetic-type equation in full 6-D phase space, e.g., the Vlasov equation  \cite{KormannReuterRampp2019} and quantum Wigner equation \cite{XiongZhangShao2022_arXiv}.  
 
In structural mechanics, the cubic spline has been also used to solve the wave propagation in cracked rod \cite{ChenYangZhangHe2012}
 and the wave-structure interaction \cite{SriramSannasirajSundar2006} within the framework of FEM, while the spline coefficients are globally dependent in principle. A major advantage of LOSS lies in its distributed construction with only local communication cost. This is based on a key observation that the wavelet basis decay exponentially in the dual space \cite{bk:Chui1992,MalevskyThomas1997}, so that  the vanished off-diagonal elements in the inverse coefficient matrix can be truncated. As a consequence, a perfectly matched boundary condition (PMBC) can be introduced to give a closure of patched spline coefficients and allows local splines to recover the global one as accurately as possible \cite{XiongZhangShao2022_arXiv}. Since only local communications in adjacent neighbors are needed, LOSS is expected to be suitable for the computational clusters with high-latency network, known as the Beowulf machines \cite{KomatitschTsuboTromp2002}. We will also show that the natural boundary conditions on two ends of splines are fully compatible with the absorbing perfectly matched layers (PML)  \cite{CollinoTsogka2001,KomatitschMartin2007,MartinKomatitsch2009,MartinKomatitschGedneyBruthiaux2010} in outer domain. At present, LOSS is readily implemented in the standard architecture and may potentially alleviate both the memory limitation and the computational burden for high-resolution wave propagation.

The rest of this paper is organized as follows. In Section \ref{sec.background}, we briefly review the background of the elastic wave propagation. Section \ref{sec.LOSS} illustrates the formulation of LOSS in both serial and parallel settings, as well as the exponential integrator for the auxiliary differential equation form of the perfectly matched layer (ADE-PML).  In Section \ref{sec.num}, we provide a series of benchmark by simulating 2-D acoustic wave equation to test the convergence and accuracy of LOSS, as well as its compatibility with ADE-PML. The  $P$- and $S$- wave propagation in 3-D homogenous or heterogenous media will also be investigated to validate the performance and parallel scalability of LOSS. Finally, conclusions and discussions are drawn in Section \ref{sec.con}.

\section{Background}
\label{sec.background}

The dynamics of wave propagation is governed by three sets of equations with $\bx \in \mathbb{R}^d$, $d\le 3$ \cite{Carcione2009}.  The first set is the conservation of linear momentum:
\begin{equation}\label{conservation_momentum}
\rho(\bx) \frac{\partial^2}{\partial t^2} u_i(\bx, t) = \frac{\partial}{\partial x_j} \sigma_{ij}(\bx, t) + f_i(\bx, t), \quad i, j = 1, \dots, 3,
\end{equation}
where $\sigma_{ij}$ are the components of the stress tensor, $u_i$ are the components of the displacement vector, $\rho$ is the mass density and $f_i$ are components of the body forces per unit (source term). The summation over repeated indices $j$ is assumed in Eq.~\eqref{conservation_momentum}. The second set is the definition of strain tensor $\varepsilon_{ij}$, which can be obtained in terms of the displacement components as
\begin{equation}\label{defintion_strain}
\varepsilon_{ij}(\bx, t) = \frac{1}{2} \left(  \frac{\partial}{\partial x_i} u_j(\bx, t) +  \frac{\partial}{\partial x_j} u_i(\bx, t) \right), \quad i, j = 1, \dots, 3.
\end{equation}

The constitutive equation reads that
\begin{equation}\label{stress_strain_relation}
\sigma_{ij}(\bx, t) = M_P(\bx) \varepsilon_{kk}(\bx, t) \delta_{ij} + 2 M_S(\bx) ( \varepsilon_{ij}(\bx, t) - \varepsilon_{kk}(\bx, t) \delta_{ij} ),
\end{equation}
where the summation over the repeated indices $k$ is assumed, $\delta_{ij}$ is the Kronecker symbol,  $M_P(\bx) =  c_P^2(\bx) \rho(\bx)$ and $M_S(\bx) =  c_S^2(\bx) \rho(\bx)$ are moduli with $c_P(\bx)$ and $c_S(\bx)$ the $P$- and $S$-velocities, respectively.

For instance, the velocity-stress form of the full elastic wave equation in 3-D media involves 15 wavefields. The conservation of momentum is 
\begin{equation*}
\begin{split}
\rho
\frac{\partial }{\partial t} 
\begin{pmatrix}
v_1 \\ v_2 \\ v_3
\end{pmatrix}
= &
\begin{pmatrix}
\frac{\partial }{\partial x_1}  & 0 &0 \\
0 & \frac{\partial }{\partial x_2}  &0\\
0 & 0  & \frac{\partial }{\partial x_3} \\
\end{pmatrix}
\begin{pmatrix}
\sigma_{11} \\
\sigma_{22} \\
\sigma_{33} \\
\end{pmatrix}
+
\begin{pmatrix}
\frac{\partial }{\partial x_2}  & \frac{\partial }{\partial x_3} &0 \\
\frac{\partial }{\partial x_1} & 0  &\frac{\partial }{\partial x_3}\\
0 & \frac{\partial }{\partial x_1}  & \frac{\partial }{\partial x_2} \\
\end{pmatrix}
\begin{pmatrix}
\sigma_{12} \\
\sigma_{13} \\
\sigma_{23} \\
\end{pmatrix}
+ 
\begin{pmatrix}
f_1 \\ f_2 \\ f_3
\end{pmatrix}.
\end{split}
\end{equation*}
The definition of strain tensor is given by
\begin{equation*}
\frac{\partial }{\partial t} 
\begin{pmatrix} 
\varepsilon_{11} \\ \varepsilon_{22} \\ \varepsilon_{33}
\end{pmatrix}
=
\begin{pmatrix}
\frac{\partial }{\partial x_1}  &0 &0 \\
0 & \frac{\partial }{\partial x_2}  & 0 \\
0 & 0  & \frac{\partial }{\partial x_3} \\
\end{pmatrix}
\begin{pmatrix} 
v_1 \\ v_2 \\v_3
\end{pmatrix},
~~
\frac{\partial }{\partial t} 
\begin{pmatrix} 
\varepsilon_{12} \\ \varepsilon_{13} \\ \varepsilon_{23}
\end{pmatrix}
=
\frac{1}{2}
\begin{pmatrix}
\frac{\partial }{\partial x_2}  & \frac{\partial }{\partial x_1} &0 \\
 \frac{\partial }{\partial x_3}  &0  & \frac{\partial }{\partial x_1}  \\
0 &  \frac{\partial }{\partial x_3}  & \frac{\partial }{\partial x_2} \\
\end{pmatrix}
\begin{pmatrix} 
v_1 \\ v_2 \\v_3
\end{pmatrix}.
\end{equation*}
And the stress-strain relation reads
\begin{equation*}
\begin{pmatrix} 
\sigma_{11} \\ \sigma_{22} \\ \sigma_{33} \\  \sigma_{12} \\ \sigma_{13} \\ \sigma_{23}
\end{pmatrix}
=
\rho
\begin{pmatrix} 
c^2_P 		& c^2_P - 2c^2_S 	&  c^2_P - 2c^2_S & 0 & 0 & 0 \\ 
c^2_P - 2c^2_S	& c^2_P	  		&  c^2_P - 2 c^2_S & 0 & 0 & 0 \\ 
c^2_P -2 c^2_S	& c^2_P - 2c^2_S	&  c^2_P		 & 0 & 0 & 0 \\ 
0 & 0 & 0 & 2 c^2_S & 0 & 0 & \\
0 & 0 & 0 & 0 & 2 c^2_S & 0 & \\
0 & 0 & 0 & 0 & 0 & 2 c^2_S & \\
\end{pmatrix}
\begin{pmatrix} 
\varepsilon_{11} \\ \varepsilon_{22} \\ \varepsilon_{33} \\  \varepsilon_{12} \\ \varepsilon_{13} \\ \varepsilon_{23}
\end{pmatrix}.
\end{equation*}

In some situations,  the $P$- wave propagation can be approximated by the acoustic wave equation based on the acoustic media assumption  \cite{Alkhalifah2000}, so that the wavefield is described by a scalar function instead of a vector, 
\begin{equation}
\rho(\bx) \frac{\partial^2}{\partial t^2} u(\bx, t) -  \nabla \cdot \left( \rho(\bx) c_P^2(\bx) \nabla u(\bx, t)\right) = f(\bx, t).
\end{equation}
Taking its two-dimensional case as an example. By introducing the velocities $v_1 = \frac{\partial }{\partial x}u$, $v_3 = \frac{\partial }{\partial z} u$ and the scalar pressure field $\sigma(x, z, t)$, it can be cast into velocity-stress form,
\begin{equation}\label{acoustic_wave}
    \begin{split}
        \frac{\partial v_1(x, z, t)}{\partial t}&=-\frac{1}{\rho(x, z)}\frac{\partial \sigma(x, z, t)}{\partial x},\\
        \frac{\partial v_3(x, z, t)}{\partial t}&=-\frac{1}{\rho(x, z)}\frac{\partial \sigma(x, z, t)}{\partial z},\\
        \frac{\partial \sigma(x, z, t)}{\partial t}&=-  \rho(x, z) c_P^2(x, z)\left(\frac{\partial v_1(x, z, t)}{\partial x}+\frac{\partial v_3(x, z, t)}{\partial z}\right).
    \end{split}
\end{equation}
It is seen that the complexity for solving the wave equation  lies in calculations of  first-order spatial derivatives of wavefields.

\section{Distributed Local spline simulator}
\label{sec.LOSS}

As a powerful tool for curve fitting, the cubic spline has been applied for solving PDEs under the framework of FEM \cite{Bermejo1990,ChenYangZhangHe2012,SriramSannasirajSundar2006}. The spline expansion is essentially global as  it requires solving global algebraic equations with tridiagonal coefficient matrice. Nonetheless, we will show that the cubic spline can be reconstructed by imposing effective inner boundary conditions on the junctions of local patches \cite{CrouseillesMehrenbergerSonnendrucker2010,XiongZhangShao2022_arXiv}.

The spline collocation method will be derived for solving the strong form of the wave equation in unidimensional space in Section \ref{sec_spline}, while  multidimensional wavefields can be constructed by the tensor product of unidimensional splines successively. For brevity, a uniform grid mesh will be adopted hereafter, but the idea is straightforward to be generalized to the non-uniform grid  due to the scaling property of wavelets \cite{ChenYangZhangHe2012}. It follows by the parallel setting of LOSS in Section  \ref{sec_distributed}, where the global spline is distributed into several local patches. The junctions are shared by adjacent nodes and the patched splines are linked by PMBCs. In the meantime, the  natural boundary conditions imposed on both ends are fully compatible with ADE-PML, where the stiffness  can be largely alleviated by the usage of exponential integrators as discussed in Section \ref{sec_exponential_integrator}.

\subsection{Spline collocation method}
\label{sec_spline}

Without loss of generality, we adopt a uniform grid mesh for the domain $[x_{\min}, x_{\max}]$ with $N+1$ points $x_{\min} = x_0 \le x_1 \le \dots \le x_N = x_{\max}$. Denote by $B_i(x)$ the cubic B-spline with compact support over four grid points \cite{CrouseillesLatuSonnendrucker2009},
\begin{equation}
B_{i}(x) = 
\left\{
\begin{split}
&\frac{(x - x_{i-2})^3}{6h^3}, \quad  x \in [x_{i-2}, x_{i-1}],\\
&-\frac{(x - x_{i-1})^3}{2h^3} + \frac{(x - x_{i-1})^2}{2h^2} + \frac{(x - x_{i-1})}{2h} + \frac{1}{6}, \quad x \in [x_{i-1}, x_{i}],\\
&-\frac{(x_{i+1} - x)^3}{2h^3} +\frac{(x_{i+1} - x)^2}{2h^2} + \frac{(x_{i+1} - x)}{2h} + \frac{1}{6}, \quad x \in [x_{i}, x_{i+1}],\\
&\frac{(x_{i+2} - x)^3}{6h^3}, \quad x \in [x_{i+1}, x_{i+2}],\\
&0, \quad \textup{otherwise},
\end{split}
\right.
\end{equation}
implying $B_{i - 1}, B_{i}, B_{i+1}, B_{i+2}$ overlap a grid interval $(x_{i}, x_{i+1})$ \cite{MalevskyThomas1997}, and
\begin{equation}\label{spline_coefficient_relation}
B_{i-1}(x_i) =  \frac{1}{6}, \quad B_{i}(x_i) =  \frac{2}{3}, \quad B_{i+1}(x_i) = \frac{1}{6}.
\end{equation}

Now the velocity wavefield can be expanded by $N+3$ splines with $N+3$ coefficients $\bm{\widetilde v}(t) = (\widetilde v_{-1}(t), \dots, \widetilde v_{N+1}(t))^T$
 \begin{equation}
\begin{split}
v_i(t) \approx \sum_{i = -1}^{N+1} \widetilde v_{i}(t) B_{i}(x), 
\end{split}
\end{equation}
where $v_i(t)$ is short for $v(x_i, t)$. In order to determine the coefficients, it suggests imposing the natural boundary conditions on two ends  to minimize the effect of boundary constraints \cite{SriramSannasirajSundar2006}, namely, $\frac{\partial^2}{\partial x^2} v(x_0, t) = 0$ and $\frac{\partial^2}{\partial x^2} v(x_N, t) = 0$. By omitting the time variable for brevity, it has that
\begin{equation}\label{natural_boundary}
\begin{split}
&\frac{1}{h^2}  \widetilde v_{-1} - \frac{2}{h^2}  \widetilde v_{0} +  \frac{1}{h^2}  \widetilde v_{1} = 0, \quad  \frac{1}{h^2}  \widetilde v_{N-1} - \frac{2}{h^2}  \widetilde v_{N} +  \frac{1}{h^2}  \widetilde v_{N+1} = 0.
 \end{split}
\end{equation}
Combining with the relation \eqref{spline_coefficient_relation}, it remains to solve the algebraic equation by the sweeping method with complexity $\mathcal{O}(N)$ \cite{bk:Boor2001},
\begin{equation}\label{spline_cofficient_eq}
A
\begin{pmatrix}
\widetilde v_{-1} \\
\widetilde v_{0} \\
\widetilde v_{1} \\
\vdots \\
\widetilde v_{N} \\
\widetilde v_{N+1} \\
\end{pmatrix}
=
\frac{1}{6}
\begin{pmatrix}
\frac{6}{h^2} & -\frac{12}{h^2} & \frac{6}{h^2} & 0 & \cdots & 0 \\
1     & 4 & 1     & 0 &           & \vdots \\
0     & 1 & 4     & 1 &  & \vdots \\
\vdots & \vdots & \vdots & \vdots & \ddots & \vdots \\ 
\vdots &           & 0 & 1 & 4 & 1 \\
0         &  0       & 0 & \frac{6}{h^2} & -\frac{12}{h^2} & \frac{6}{h^2} \\
\end{pmatrix}
\begin{pmatrix}
\widetilde v_{-1} \\
\widetilde v_{0} \\
\widetilde v_{1} \\
\vdots \\
\widetilde v_{N} \\
\widetilde v_{N+1} \\
\end{pmatrix}
=
\begin{pmatrix}
0 \\
v_0 \\
 v_1 \\
\vdots \\
v_N \\
0 \\
\end{pmatrix}.
\end{equation}

Once the spline coefficients are obtained, the spatial first-order derivatives can be directly approximated by 
\begin{equation}\label{derivative_approximation}
\frac{\partial }{\partial x}v(x_i, t) \approx -\frac{1}{2h}\widetilde v_{i-1}(t) + \frac{1}{2h}\widetilde v_{i+1}(t)
\end{equation}
as
\begin{equation}
\frac{\partial}{\partial x} B_{i}(x) = 
\left\{
\begin{split}
&\frac{(x - x_{i-2})^2}{2h^3}, \quad  x \in [x_{i-2}, x_{i-1}],\\
&-\frac{3(x - x_{i-1})^2}{2h^3} + \frac{(x - x_{i-1})}{h^2} + \frac{1}{2h}, \quad x \in [x_{i-1}, x_{i}],\\
&\frac{3(x_{i+1} - x)^2}{2h^3} - \frac{(x_{i+1} - x)}{h^2} - \frac{1}{2h}, \quad x \in [x_{i}, x_{i+1}],\\
&-\frac{(x_{i+2} - x)^2}{2h^3}, \quad x \in [x_{i+1}, x_{i+2}],\\
&0, \quad \textup{otherwise}.
\end{split}
\right.
\end{equation}

\subsection{Distributed local spline}
\label{sec_distributed}

For distributed parallelization,  the spline needs to be decomposed into several pieces and stored in multiple processors. Here we simply divide $N+1$ grid points on a line into $p$ uniform parts, with $M = N/p$,
\begin{align*}
\underbracket{v_0 < v_1 < \cdots < v_{M-1}}_{\textup{the first processor}} < \underbracket{v_{M}}_{\textup{shared}} < \cdots < \underbracket{ v_{(p-1)M}}_{\textup{shared}} < \underbracket{v_{(p-1)M+1} < \cdots < v_{pM}}_{\textup{$p$-th processor}},
\end{align*}
The grid points $v_{M}, v_{2M}, \dots, v_{(p-1)M}$ are shared by the adjacent patches. Recall that our target is to recover the global B-spline by the local spline coefficients $\bm{\widetilde v}^{(l)}$ for $l$-th piece without global communications, that is, 
\begin{equation}\label{relation}
\bm{\widetilde v}^{(l)}= (\widetilde v_{-1}^{(l)}, \dots, \widetilde v_{M+1}^{(l)}) = (\widetilde v_{-1 +(l-1)M}, \dots, \widetilde v_{(l-1)M+M+1}), \quad l = 1, \dots, p.
\end{equation}

\begin{figure}[!h]
\centering
\subfigure[The element $|b_{ij}|$ of $A^{-1}$ in log10 scale.\label{element_decay}]{
\includegraphics[width=0.49\textwidth,height=0.24\textwidth]{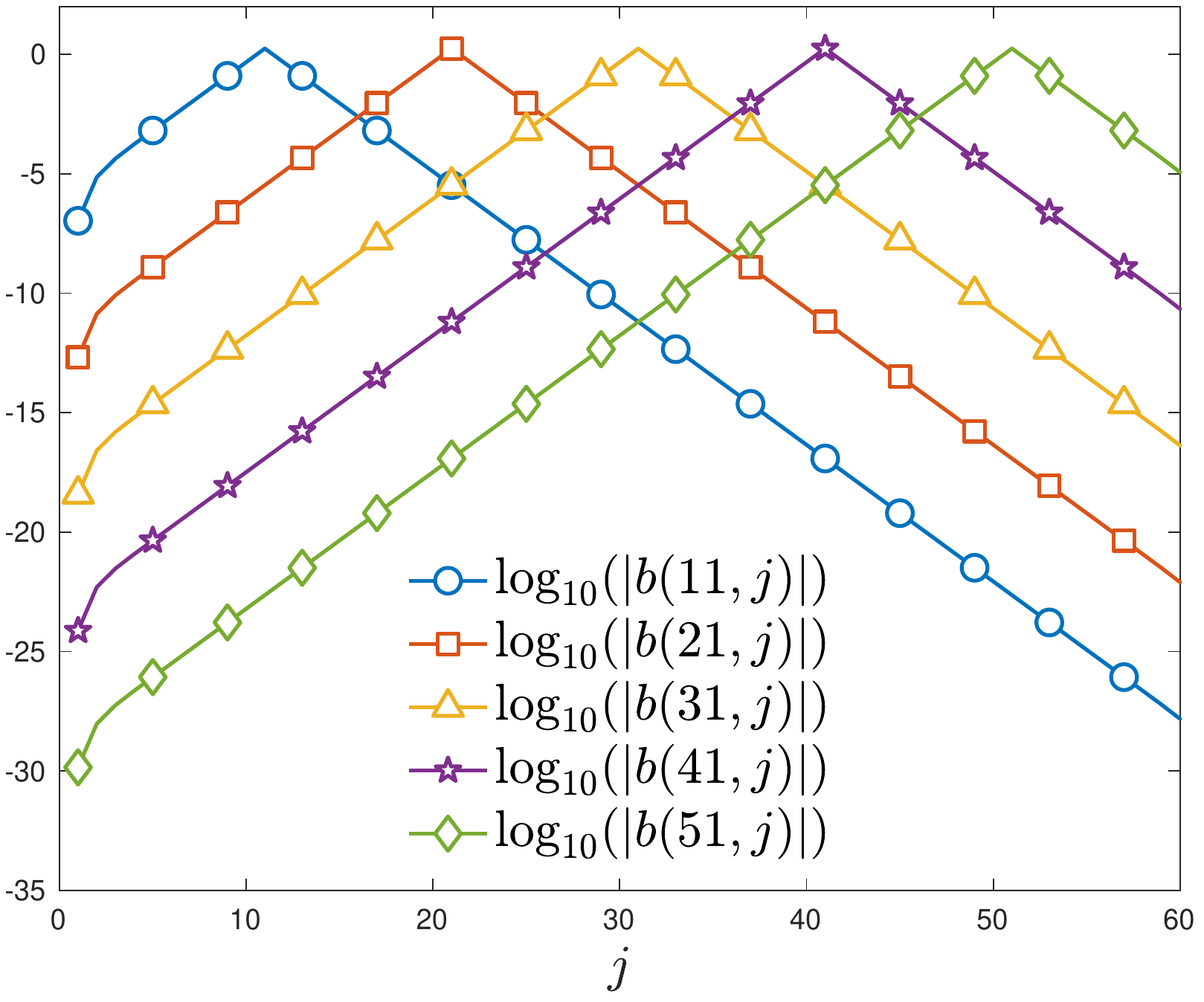}}
\subfigure[ Approximation of $A^{-1}$ by $(A_M^{(l)})^{-1}$.\label{domain_decomposition}]
{\includegraphics[width=0.49\textwidth,height=0.24\textwidth]{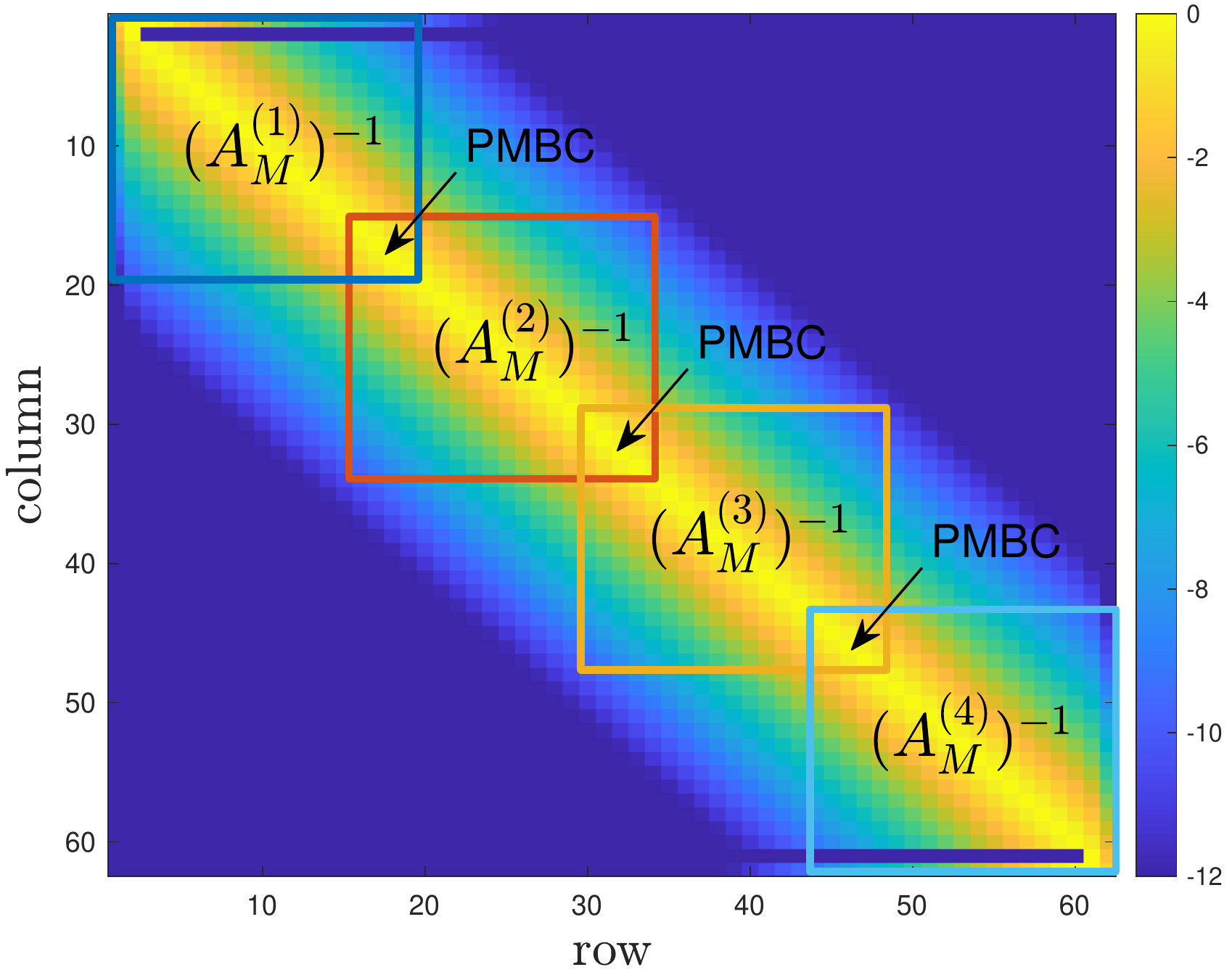}}
\caption{Since the non-diagonal elements $b_{ij}$ in $A^{-1}$ decays exponentially away from the main diagonal $b_{ii}$, the coefficients $\widetilde{\bv} = A^{-1}(v_0, \dots, v_N)^{\tau}$ can be approximated by $\widetilde{\bv}^{(l)} =  (A_M^{(l)})^{-1}(v_0^{(l)}, \dots, v_M^{(l)})^\tau$ when PMBCs are adopted.}
\end{figure}

This can be realized by imposing effective Hermite boundary conditions on two ends of local splines (see Figure \ref{domain_decomposition}) \cite{CrouseillesLatuSonnendrucker2009,XiongZhangShao2022_arXiv},
\begin{equation}
\frac{\partial v}{\partial x} \Big |_{x = x_{(l-1)M}} = \phi_L^{(l)}, \quad \frac{\partial v}{\partial x} \Big |_{x = x_{lM}} = \phi_R^{(l)}, 
\end{equation}
which is equivalent to
\begin{equation}
-\frac{1}{2h} \widetilde v_{-1}^{(l)} + \frac{1}{2h} \widetilde v_{1}^{(l)}= \phi_L^{(l)}, \quad -\frac{1}{2h} \widetilde v_{M-1}^{(l)} + \frac{1}{2h} \widetilde v_{M+1}^{(l)}= \phi_R^{(l)}.
\end{equation}
Thus all the coefficients $ \bm{\widetilde v}^{(l)}= (\widetilde v_{-1}^{(l)}, \dots, \widetilde v_{M+1}^{(l)})$ can be obtained straightforwardly by solving the algebraic equation 
\begin{equation}\label{spline_cofficient_eq_piece}
A_{M}^{(l)}
\begin{pmatrix}
\widetilde v^{(l)}_{-1} \\
\widetilde v^{(l)}_{0} \\
\widetilde v^{(l)}_{1} \\
\vdots \\
\widetilde v^{(l)}_{M} \\
\widetilde v^{(l)}_{M+1} \\
\end{pmatrix}
=
\frac{1}{6}
\begin{pmatrix}
-\frac{3}{h} & 0 & \frac{3}{h} & 0 & \cdots & 0 \\
1     & 4 & 1     & 0 &           & \vdots \\
0     & 1 & 4     & 1 &  & \vdots \\
\vdots & \vdots & \vdots & \vdots & \ddots & \vdots \\ 
\vdots &           & 0 & 1 & 4 & 1 \\
0         &  0       & 0 & -\frac{3}{h} & 0 & \frac{3}{h} \\
\end{pmatrix}
\begin{pmatrix}
\widetilde v^{(l)}_{-1} \\
\widetilde v^{(l)}_{0} \\
\widetilde v^{(l)}_{1} \\
\vdots \\
\widetilde v^{(l)}_{M} \\
\widetilde v^{(l)}_{M+1} \\
\end{pmatrix}
=
\begin{pmatrix}
\phi^{(l)}_L \\
v^{(l)}_0 \\
 v^{(l)}_1 \\
\vdots \\
v^{(l)}_M \\
\phi^{(l)}_R \\
\end{pmatrix}.
\end{equation}
where $(M+3)\times (M+3)$ coefficient matrix $A^{(l)}_{M}$ has an explicit LU decomposition, namely, $A^{(l)}_{M} = L U$,
\begin{equation}\label{explicit_L}
L = 
\begin{pmatrix}
1 & 0 & 0 &\cdots & \cdots & 0 \\
-\frac{h}{3} & 1 & 0  &\ddots & & \vdots \\
0 & l_1 & 1  &\ddots & & \vdots \\
0 & 0    & l_2  & \ddots & & \vdots \\
\vdots & \vdots   & & l_{M} & 1 & 0 \\
0 & 0 & \cdots  & -\frac{3l_{M}}{h} & \frac{3l_{M+1}}{h} & 1 \\
\end{pmatrix}
\end{equation}
and
\begin{equation}\label{explicit_U}
U = \frac{1}{6}
\begin{pmatrix}
-\frac{3}{h} & 0 & \frac{3}{h} & 0 &\cdots & \cdots & 0 \\
0 & d_1 & 2 & 0 &\ddots & & \vdots \\
0 & 0 & d_2 & 1 &\ddots & & \vdots \\
0 & 0    & 0  & d_3 & \ddots & & \vdots \\
\vdots & \vdots  & & & 0 & d_{M+1} & 0 \\
0 & 0 & \cdots  & & 0 & 0 & \frac{3 d_{M+2}}{h} \\
\end{pmatrix},
\end{equation}
with
\begin{equation}
\begin{split}
& d_1 = 4, \quad l_1 = 1/4, \quad d_2 = 4 - 2l_1 = 7/2, \\
& l_i = 1/d_i, \quad d_{i+1} = 4 - l_i, \quad i = 2, \dots, M+1, \\
&  l_{M+1} = 1/(d_{M} d_{M+1}), \quad d_{M+2} = 1 - l_{M+1}.
\end{split}
\end{equation}

Now the solution of Eq.~\eqref{spline_cofficient_eq_piece}  should be equivalent to that of Eq.~\eqref{spline_cofficient_eq} by choosing appropriate $\phi^{(l)}_L$ and $\phi^{(l)}_R$. Denote by $(b_{ij}) = {A}^{-1}, -1\le i, j \le pM+1$, thus the solution $\widetilde v_i$ of Eq.~\eqref{spline_cofficient_eq} can be represented by
\begin{equation}\label{exact_solution}
\widetilde v_i = b_{ii} v_i + \sum_{j=-1}^{i-1} b_{ij} v_j + \sum_{j = i+1}^{pM+1} b_{i j} v_j, \quad i = -1, \dots, pM+1.
\end{equation}
Using the constraints  \eqref{relation}, it directly solves $\phi_L^{(l)}$ and $\phi_R^{(l)}$ by
 \begin{equation}\label{boundary_condition}
 \begin{split}
& \phi^{(l)}_L = -\frac{1}{2h} \widetilde v_{-1}^{(l)} + \frac{1}{2h} \widetilde v_{1}^{(l)}=  -\frac{1}{2h} \widetilde v_{-1 + (l-1)M} + \frac{1}{2h} \widetilde v_{1+ (l-1)M}, \\
& \phi^{(l)}_R = -\frac{1}{2h} \widetilde v_{M-1}^{(l)} + \frac{1}{2h} \widetilde v_{M+1}^{(l)}=  -\frac{1}{2h} \widetilde v_{M-1 + (l-1)M} + \frac{1}{2h} \widetilde v_{M+1+ (l-1)M},
 \end{split}
 \end{equation}
 where $\widetilde v_{i}$ are given by Eq.~\eqref{exact_solution}. 
 
 At first glance,  it still requires the information from all other pieces. Fortunately, there is a key observation that the non-diagonal elements $b_{ij}$ decays exponentially away from the main diagonal $b_{ii}$ due to the rapid decay of the wavelet basis in its dual space \cite{bk:Chui1992,MalevskyThomas1997}, which is clearly visualized in Figure \ref{element_decay}. Therefore, it allows us to truncate $b_{ij}$ when $|i - j| \ge n_{nb}$ for sufficiently large $n_{nb}$,
 \begin{equation}\label{truncated_solution}
\widetilde v_i \approx b_{ii} v_i + \sum_{j=i - (n_{nb}-1)}^{i-1} b_{ij} v_j + \sum_{j = i+1}^{i + n_{nb} - 1} b_{i j} v_j,  \quad i = -1, \dots, pM+1.
\end{equation}
Now using the truncated stencils \eqref{truncated_solution}, it has that
\begin{equation}\label{truncation}
\begin{split}
& \widetilde v_{lM-1} \approx \sum_{j=(lM-1)-n_{nb}+1}^{(lM-1)+n_{nb}-1} b_{lM-1, j} v_j =  \sum_{j=-n_{nb}}^{n_{nb}-2} b_{lM-1, lM+j} v_{lM+j}, \\
& \widetilde v_{lM+1} \approx \sum_{j=(lM+1)-n_{nb}+1}^{(lM+1)+n_{nb}-1} b_{lM+1, j} v_j = \sum_{j=-n_{nb}+2}^{n_{nb}} b_{lM+1, lM+ j} v_{lM+j}. \\ 
\end{split}
\end{equation}
By further adding four more terms in Eq.~\eqref{truncation} to complete the summations from $-n_{nb}$ to $n_{nb}$, it yields that 
\begin{equation*}
\begin{split}
-\frac{1}{2h} \widetilde v_{lM-1} + \frac{1}{2h} \widetilde v_{lM+1} \approx &\underbracket{\left(-\frac{1}{2h} b_{lM-1, lM} + \frac{1}{2h} b_{lM+1, lM}\right) v_{lM}}_{\textup{shared by adjacent two processors}} \\
= &\underbracket{\sum_{j=-n_{nb}}^{-1} \left(-\frac{1}{2h} b_{lM-1, lM+j} + \frac{1}{2h} b_{lM+1, lM+j}\right) v_{lM+j}}_{\textup{stored in left processor}} \\
&~+\underbracket{\sum_{j=1}^{n_{nb}} \left(-\frac{1}{2h} b_{lM-1, lM+j} + \frac{1}{2h} b_{lM+1, lM+j}\right) v_{lM+j}}_{\textup{stored in right processor}}. 
\end{split}
\end{equation*}
Thus it arrives at the formulation of PMBC for $1 \le l \le p-1$,
\begin{equation}\label{PMBC}
\begin{split}
\phi_{R}^{(l)} =  \phi_{L}^{(l+1)} \approx & \underbracket{\frac{1}{2}c_{0,l} v_{lM} +  \sum_{j = 1}^{n_{nb}} c_{j, l}^{-} v_{lM-j}}_{\textup{stored in left processor}} + \underbracket{\frac{1}{2}c_{0,l} v_{lM} + \sum_{j = 1}^{n_{nb}} c_{j, l}^{+} v_{lM+j}}_{\textup{stored in right processor}},
\end{split}
\end{equation}
where $c_{0, l} = -\frac{b_{lM-1, lM}}{2h} + \frac{b_{lM+1, lM}}{2h}$ and 
\begin{equation}\label{PMBC_coeffcients}
\begin{split}
&c_{j, l}^+ =  -\frac{b_{lM-1, lM+j}}{2h} + \frac{b_{lM+1, lM+j}}{2h}, \quad c_{j, l}^- =  -\frac{b_{lM-1, lM-j}}{2h} + \frac{b_{lM+1, lM-j}}{2h}.
\end{split}
\end{equation}

Finally, the effective boundary conditions $\phi_{L}^{(1)}$ and $\phi_{R}^{(p)}$ on two ends should match the natural boundary conditions of the global cubic spline, yielding
\begin{align}
\phi_{L}^{(1)} &=\frac{\widetilde v_1-\widetilde v_{-1}}{2h}\approx \sum_{j=0}^{n_{nb}} c_{j, 0}^{-} v_j, \quad c_{j, 0}^{-} =  \frac{-{b}_{-1, j}  + {b}_{1, j}}{2h}, \\
\phi_{R}^{(p)} &= \frac{\widetilde v_{N+1}-\widetilde v_{N-1}}{2h} \approx \sum_{j=0}^{n_{nb}} c_{j, p}^{+}v_{N-j}, \quad c_{j, p}^{+} = \frac{-{b}_{N-1, N-j}  + {b}_{N+1, N-j}}{2h}.
\end{align}

\subsection{Exponential integrator for the wave equation and ADE-PML}
\label{sec_exponential_integrator}

PML is usually adopted in finite computational domain to model wave propagation in unbounded media \cite{CollinoTsogka2001,KomatitschMartin2007,MartinKomatitsch2009,MartinKomatitschGedneyBruthiaux2010}, which suggests adding a layer of the thickness $L$ outside to attenuate the wave and avoid its artificial reflection. It will be shown that the natural boundary condition of cubic spline can be imposed on the differential equation form of PML (known as ADE-PML \cite{MartinKomatitschGedneyBruthiaux2010}), and the resulting rigid dynamics can be solved efficiently by exponential integrators \cite{HochbruckOstermann2005,XiongGuo2022}.

For convenience, we use the 2-D acoustic wave equation to illustrate the basic idea of ADE-PML. It starts by taking the Fourier transform of the acoustic wave equation in the tensile coordinates ($t \to \omega$) and adding two stretching terms $1/s_x$ and $1/s_z$ \cite{MartinKomatitschGedneyBruthiaux2010},
\begin{align}
        \mi\omega\widehat{v}_1(x, z, \omega)&=-\frac{1}{\rho} \frac{1}{s_x}  \widehat{\frac{\partial \sigma}{\partial x}}(x, z, \omega),\label{acoustic_eq1}\\
        \mi\omega\widehat{v}_3(x, z, \omega)&=-\frac{1}{\rho} \frac{1}{s_z}  \widehat{\frac{\partial \sigma}{\partial z}}(x, z, \omega),\label{acoustic_eq2}\\
        \mi\omega\widehat{\sigma}(x, z, \omega)&=-\rho c_P^2 \left( \frac{1}{s_x}  \widehat{\frac{\partial v_1}{\partial x}}(x, z, \omega)
                                + \frac{1}{s_z} \widehat{\frac{\partial v_3}{\partial z}}(x, z, \omega)\right),\label{acoustic_eq3}
\end{align}
where $\widehat{\cdot}$ denotes the Fourier transform, and $s_i$ ($i=x, z$) are introduced to present a complex stretching of the coordinate system,
\begin{equation}
    \begin{split}
  s_i=k_i+\frac{d_i}{\alpha_i+ \mi \omega}, \quad s_i^{-1}=\frac{1}{k_i}-\frac{d_i}{k_i^2}\frac{1}{\frac{d_i}{k_i}+\alpha_i+\mi \omega}, ~~ \quad i = x, z.
    \end{split}
\end{equation}
Here $\alpha_i$, $d_i$ and $k_i$ are flexible parameters that adjust the absorbing effect,
\begin{equation}
        d_i=d_0\left(\frac{i}{L}\right)^2, \quad k_i=  1 + (k_{\max} - 1)^m, \quad  \alpha_i=\alpha_{\max}\left(1-\frac{i}{L}\right)^p
\end{equation}
and $d_0=-3c_{P, {\max}}\log(R)/2L$, with $c_{P, {\max}}$ equal to the maximal velocity of the pressure wave.
$R$ is the theoretical reflection coefficient of the target, 
$\alpha_{\max} = \pi f_0$ with $f_0$ as the main frequency 
of the hypocenter, and $k_{\max}$ is the cut-off wavenumber \cite{MartinKomatitschGedneyBruthiaux2010}. Here it suffices to take $m=1$ and $p=1$ for simplicity.

For the first equation \ref{acoustic_eq1} of $\widehat{\frac{\partial \sigma}{\partial x}}(x, z, \omega)$, it is equivalent to
\begin{equation}
    \frac{1}{s_x} \widehat{\frac{\partial \sigma}{\partial x}}(x, z, \omega)
    =\frac{1}{k_x}\widehat{\frac{\partial \sigma}{\partial x}}(x, z, \omega)
    -\frac{d_x}{k_x^2}\frac{1}{\frac{d_x}{k_x}+\alpha_x+\mi \omega}
    \widehat{\frac{\partial \sigma}{\partial x}}(x, z, \omega).
\end{equation}
By introducing a memory variable  $\widehat{\Psi}_x(x, z, \omega)$,
\begin{equation}
    \widehat{\Psi}_x(x, z, \omega)=-\frac{d_x}{k_x^2}\frac{1}{\frac{d_x}{k_x}+\alpha_x+ \mi \omega}
    \widehat{\frac{\partial \sigma}{\partial x}}(x, z, \omega),
\end{equation}
it arrives at an auxiliary differential equation by taking inverse Fourier transform $\widehat{\Psi}_x(x, z, \omega) \to \Psi_x(x, z, t)$, 
\begin{equation}\label{ADE_PML_1}
    \frac{\partial \Psi_x}{\partial t}(x, z, t) + \left(\frac{d_x}{k_x}+\alpha_x\right)\Psi_x(x, z, t)
    =-\frac{d_x}{k_x^2}\frac{\partial \sigma(x, z, t)}{\partial x}.
\end{equation}

Similarly, we can define memory variables ${\Psi}_z(x, z, t)$ for $ \widehat{\frac{\partial \sigma}{\partial z}}(x, z,\omega)$ in Eq.~\ref{acoustic_eq2} and $\Phi_x(x, z, t)$, $\Phi_z(x, z, t)$ for $ \widehat{\frac{\partial v_1}{\partial x}}(x, z,\omega)$, $ \widehat{\frac{\partial v_3}{\partial z}}(x, z,\omega)$ in Eq.~\ref{acoustic_eq3}, respectively, yielding the ADE-PML for velocity components outside the computational domain,
\begin{equation}\label{ADE_PML_A}
\textup{(A)}\left\{
\begin{split}
&\frac{\partial v_1}{\partial t} = -\frac{1}{\rho} \left(\frac{1}{k_x} \frac{\partial \sigma}{\partial x} + \Psi_x\right), \\
&\frac{\partial v_3}{\partial t} = -\frac{1}{\rho} \left(\frac{1}{k_z} \frac{\partial \sigma}{\partial z} + \Psi_z\right), \\
&\frac{\partial \Psi_x}{\partial t} = -\left(\frac{d_x}{k_x} + \alpha_x\right) \Psi_x - \frac{d_x}{k_x^2} \frac{\partial \sigma}{\partial x}, \\
&\frac{\partial \Psi_z}{\partial t}  = -\left(\frac{d_z}{k_z} + \alpha_z\right) \Psi_z - \frac{d_z}{k_z^2} \frac{\partial \sigma}{\partial z}, \\
\end{split}
\right.
\end{equation}
and those for the stress component,
\begin{equation}\label{ADE_PML_B}
\textup{(B)}\left\{
\begin{split}
&\frac{\partial \sigma}{\partial t} = -\rho c_P^2 \left(\frac{1}{k_x} \frac{\partial v_1}{\partial x} + \Phi_x + \frac{1}{k_z} \frac{\partial v_3}{\partial z} + \Phi_z\right),  \\
&\frac{\partial \Phi_x}{\partial t}  = -\left(\frac{d_x}{k_x} + \alpha_x\right) \Phi_x - \frac{d_x}{k_x^2} \frac{\partial v_1}{\partial x}, \\
&\frac{\partial \Phi_z}{\partial t}  = -\left(\frac{d_z}{k_z} + \alpha_z\right) \Phi_z - \frac{d_z}{k_z^2} \frac{\partial v_3}{\partial z}. \\
\end{split}
\right.
\end{equation}


The stiff terms in the auxiliary differential equations might pose severe limitation on the time step when using explicit numerical integrators. Fortunately, this can be alleviated by the exponential integrator \cite{HochbruckOstermann2005,XiongGuo2022}.

For the time interval $[t_n, t_{n+1}]$ with $\Delta t = t_{n+1} - t_n$, it starts from the variation-of-constant formula of Eq.~\eqref{ADE_PML_1},
\begin{equation*}
    \Psi_x(\bx, t_{n+1})=\me^{- (\frac{d_x}{k_x}+\alpha_x)\Delta t}\Psi_x(\bx, t)-\frac{d_x}{k_x^2} \int^{\Delta t}_0 \me^{- (\frac{d_x}{k_x}+\alpha_x)(\Delta t-\tau)}\frac{\partial \sigma}{\partial x}(\bx, t_n+\tau) \mathrm{d} \tau.
\end{equation*}
When the Euler approximation $\frac{\partial \sigma}{\partial x}(\bx, t+\tau) \approx \frac{\partial \sigma}{\partial x}(\bx, t)$ is used, it yields \cite{HochbruckOstermann2005},
\begin{equation}\label{Euler_approximation}
    \Psi_x(\bx, t_{n+1})\approx \me^{-  (\frac{d_x}{k_x}+\alpha_x) \Delta t}\Psi_x(\bx, t_n)- \frac{d_x}{k_x^2}
    \frac{1-\me^{-  (\frac{d_x}{k_x}+\alpha_x)\Delta t }}{ (\frac{d_x}{k_x}+\alpha_x) }
    \frac{\partial \sigma}{\partial x}(\bx, t_n).
\end{equation}
The other three auxiliary equations for $\Phi_z$, $\Psi_x$ and $\Psi_z$ can be tackled in a similar way. Combining with the temporal finite difference scheme for $v_1$, $v_3$ and $\sigma$, one can obtain the non-splitting exponential Euler scheme for ADE-PML. Because the exact flow of the stiff term in the auxiliary dynamics \eqref{ADE_PML_1} is exploited, it can largely alleviate the restriction on the time step.

Alternatively, one can utilize the exponential operator splitting scheme, which also utilizes the exact stiff flow of the auxiliary equations. The basic idea is alternating update of velocities and stress based on the splitting of two subproblems \eqref{ADE_PML_A} and \eqref{ADE_PML_B}, which is similar to our short-memory operator splitting for time-fractional constant-Q wave equation  \cite{XiongGuo2022}. 

One can first solve (A) exactly by assuming that $\sigma$, $\Phi_x$ and $\Phi_z$ are invariant in small time step. For instance, when $\frac{\partial \sigma}{\partial x}$ is invariant in $[t_n, t_{n+1}]$, it yields the exact solution of $\Psi_x$,
\begin{equation}
\begin{split}
\Psi_x(\bx, t_{n+1}) = \me^{-(\frac{d_x}{k_x} + \alpha_x) \Delta t} \Psi_x(\bx, t_n) - \frac{d_x}{k_x^2} \frac{1 - \me^{-(\frac{d_x}{k_x} + \alpha_x) \Delta t}}{(\frac{d_x}{k_x} + \alpha_x) } \frac{\partial \sigma}{\partial x}(\bx, t_n). \\
\end{split}
\end{equation}
In addition, since
\begin{equation*}
\begin{split}
\int_{t_n}^{t_{n+1}}\frac{\partial \Psi_x}{\partial t} \D t &= -\left(\frac{d_x}{k_x} + \alpha_x\right) \int_{t_n}^{t_{n+1}} \Psi_x(\bx, t) \D t - \frac{d_x}{k_x^2}  \int_{t_n}^{t_{n+1}} \frac{\partial \sigma}{\partial x}(\bx, t) \D t \\
& =  -\left(\frac{d_x}{k_x} + \alpha_x\right) \int_{t_n}^{t_{n+1}} \Psi_x(\bx, t) \D t  - \Delta t \frac{d_x}{k_x^2}  \frac{\partial \sigma}{\partial x}(\bx, t_n),
\end{split}
\end{equation*}
it further yields  
\begin{equation*}
\int_{t_n}^{t_{n+1}} \Psi_x(\bx, t) \D t = - \frac{1}{(\frac{d_x}{k_x} + \alpha_x)} \left( \Psi(\bx, t_{n+1}) - \Psi(\bx, t_n) +  \Delta t \frac{d_x}{k_x^2}  \frac{\partial \sigma}{\partial x}(\bx, t_n) \right).
\end{equation*}
Substituting it into Eq.~\eqref{ADE_PML_A}, then the velocity $v_1$ can be solved by
\begin{equation}
\begin{split}
v_1(\bx, t_{n+1}) = &v_1(\bx, t_{n}) - \frac{1}{\rho(\bx)} \left(  \frac{\Delta t}{k_x} \frac{\partial \sigma}{\partial x}(\bx, t_n)  + \int_{t_n}^{t_{n+1}} \Psi_x(\bx, t) \D t \right)\\
 = &v_1(\bx, t_{n}) + \frac{1}{\rho(\bx)}\frac{k_x}{ (d_x + \alpha_xk_x)} ( \Psi_x(\bx, t_{n+1}) - \Psi_x(\bx, t_n)  ) \\
 &-\frac{\Delta t}{\rho(\bx)} \frac{\alpha_x }{(d_x + \alpha_x k_x) }  \frac{\partial \sigma}{\partial x}(\bx, t_n).
\end{split}
\end{equation}
The solution of $v_3$ and $\Phi_z$ can be obtained in the same way. Besides, one can also solve the subsystem (B) exactly when  $v_1$, $v_3$, $\Psi_x$ and $\Psi_z$ are assumed to be invariant in small time step.  Specifically, when the Strang splitting is used, say, half step evolution of (A) + full step evolution of (B) +  half step evolution of (A), it is expected to achieve  global second-order convergence as two exact flows are exploited.

\section{Numerical experiments}
\label{sec.num}

From this section, several benchmarks have been performed to evaluate the performance of LOSS. In the first example, we made a series of benchmarks on 2-D acoustic wave equation in homogenous media to test the convergence of the spline collocation method, where  the ADE-PML associated with the natural boundary condition was adopted. In particular, the influence of several key parameters, including the layer thickness $L$, the reflection coefficient $R$ and the cut-off wavenumber $k_{\max}$, were carefully studied. After that, we used LOSS to simulate 3-D wave propagation in either a homogenous media or a  double-layer media activated by the impulse of a Ricker-type wavelet history.  These typical examples may validate the performance of LOSS when the coefficients are either smooth or of a large variation, as well as its parallel scalability.

To evaluate the errors of LOSS, we adopted two metrics: the relative $l^2$-error $\varepsilon_2(t)$ and the relative maximal error $\varepsilon_{\infty}(t)$, where $v_3^{\textup{num}}$ and $v_3^{\textup{ref}}$ denoted the numerical and reference velocity, respectively. 
\begin{equation*}
\begin{split}
&\varepsilon_2(t) = \frac{(\sum_{i} |v_3^{\textup{num}}(\bx_i, t) - v_3^{\textup{ref}}(\bx_i, t)|^2)^{1/2}}{\max_{\bx}|v_3^{\textup{ref}}(\bx, t)|}, \quad \varepsilon_\infty(t) = \frac{\max_{\bx}|v_3^{\textup{num}}(\bx, t) - v_3^{\textup{ref}}|}{\max_{\bx}|v_3^{\textup{ref}}(\bx, t)|},
\end{split}
\end{equation*}

All the 2-D simulations were realized by MATLAB, while all the 3-D simulations performed via our Fortran implementations  ran on the platform: AMD Ryzen 7950X (4.50GHz, 64MB Cache, 16 Cores, 32 Threads) with 64GB Memory (4800Mhz). The parallelization was realized by the Message Passing Interface (MPI).

\subsection{2-D acoustic wave equation}

First we need to validate the convergence of the spline collocation method. The model parameters were set as: the wave speed $c_P = 50$m/s and the density $\rho = 1$kg$/\mathrm{m}^3$. The computational domain was $[-5, 5]^2$ ($10$m$\times$$10$m). The Strang operator splitting was adopted with time step $\Delta t=0.0001$s and the final time was $T = 0.2$s. The parameters of PML were given by: the thickness $L = 50$, the reflection coefficient $R = 10^{-6}$, the cut-off wavenumber $k_{\max}=1$ and $f_P = c_P/L$, $a_{\max} = \pi f_P$. 

Five groups of spline simulations under the grid size $N = N_x \times N_z = 32^2, 64^2, 128^2, 256^2, 512^2$ were performed, with zero initial velocity and the initial pressure
\begin{equation}
\sigma(x, z, 0) = \me^{-5(x^2+ z^2)}.
\end{equation}
The reference solutions were produced by the Fourier spectral method (FSM) with a $N_x\times N_z = 256^2$ grid, where the domain was extended to $27$m$\times27$m with $712^2$ grid points to avoid the reflection of waves. 

\begin{figure}[!h]
    \centering
    \subfigure[$t=0.02$s.]{
    \begin{minipage}[t]{0.33\linewidth}
    \centering
    \includegraphics[scale=0.35]{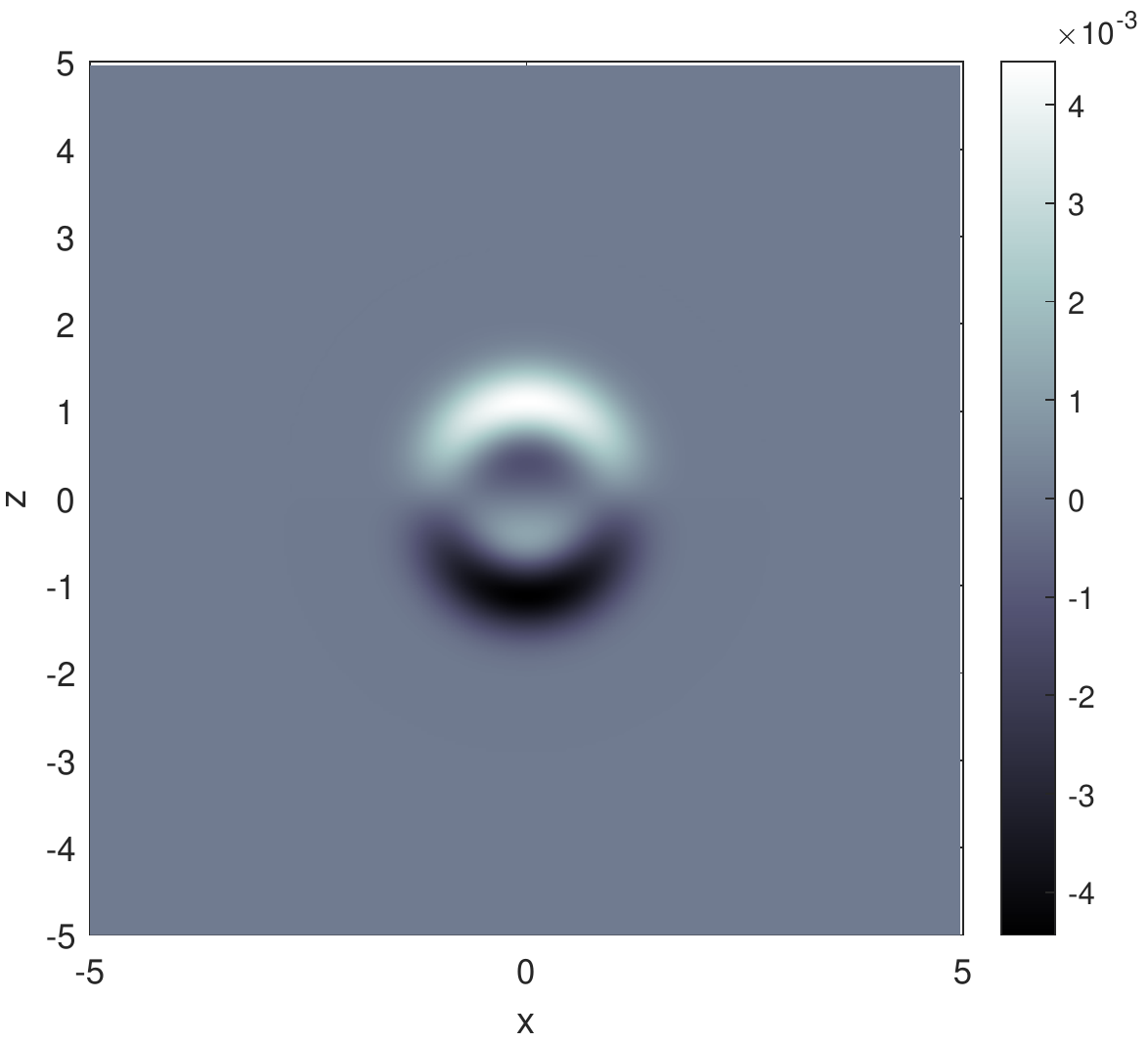}
    \end{minipage}%
    }%
    \subfigure[$t=0.06$s.]{
    \begin{minipage}[t]{0.33\linewidth}
    \centering
    \includegraphics[scale=0.35]{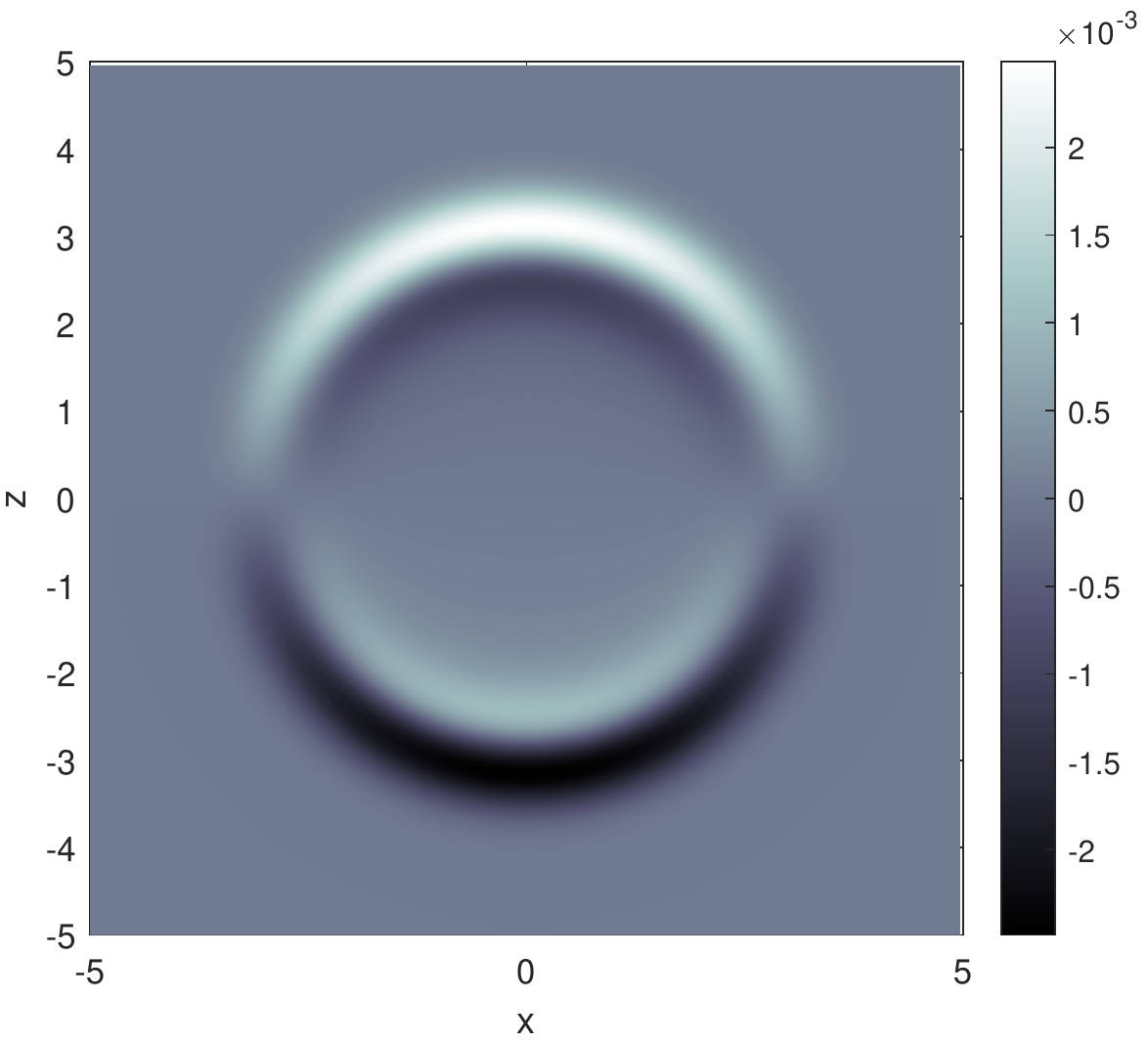}
    \end{minipage}%
    }%
    \subfigure[$t=0.10$s.]{
    \begin{minipage}[t]{0.33\linewidth}
    \centering
    \includegraphics[scale=0.35]{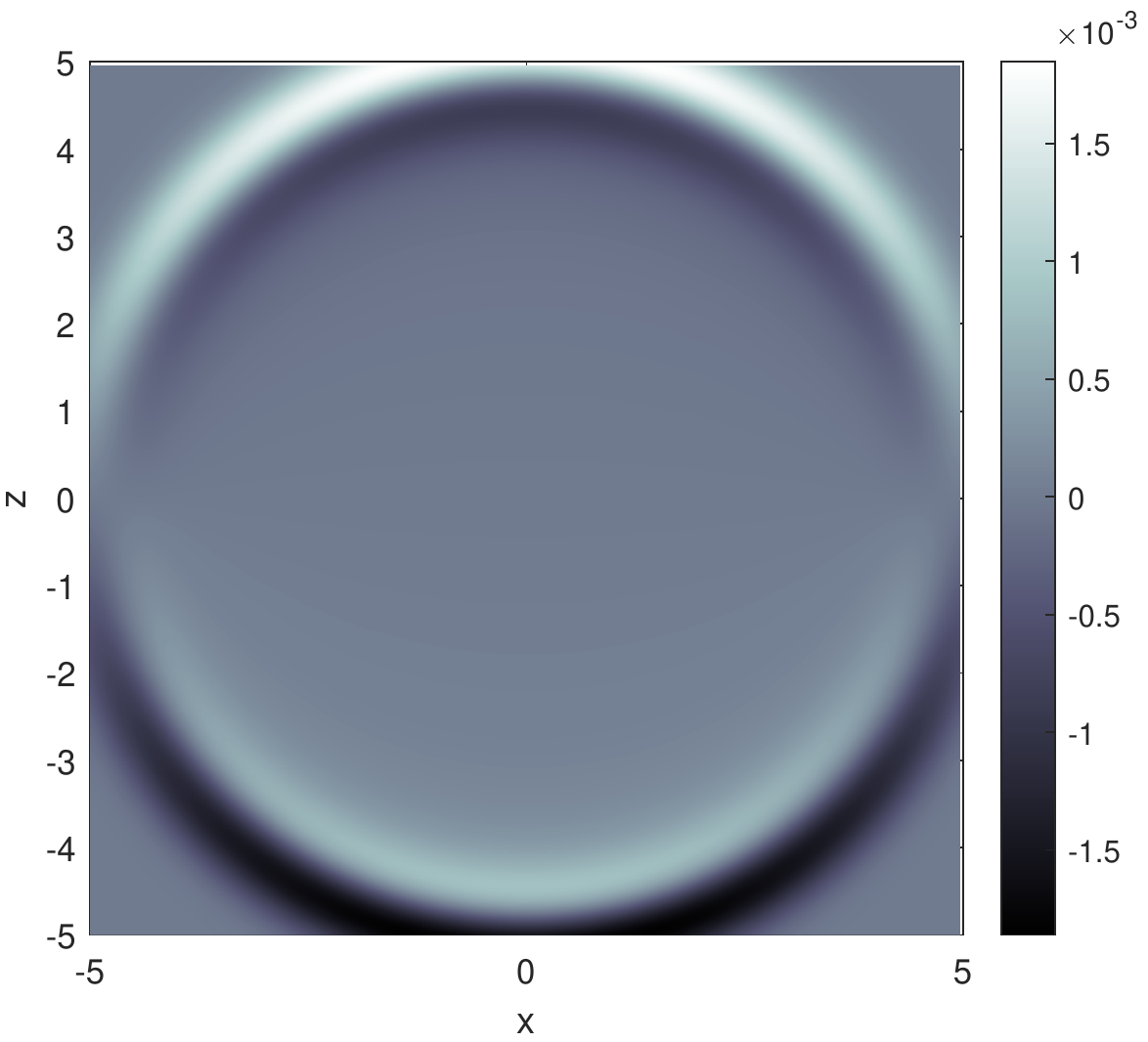}
    \end{minipage}%
    }%

    \subfigure[$t=0.12$s.]{
    \begin{minipage}[t]{0.33\linewidth}
    \centering
    \includegraphics[scale=0.35]{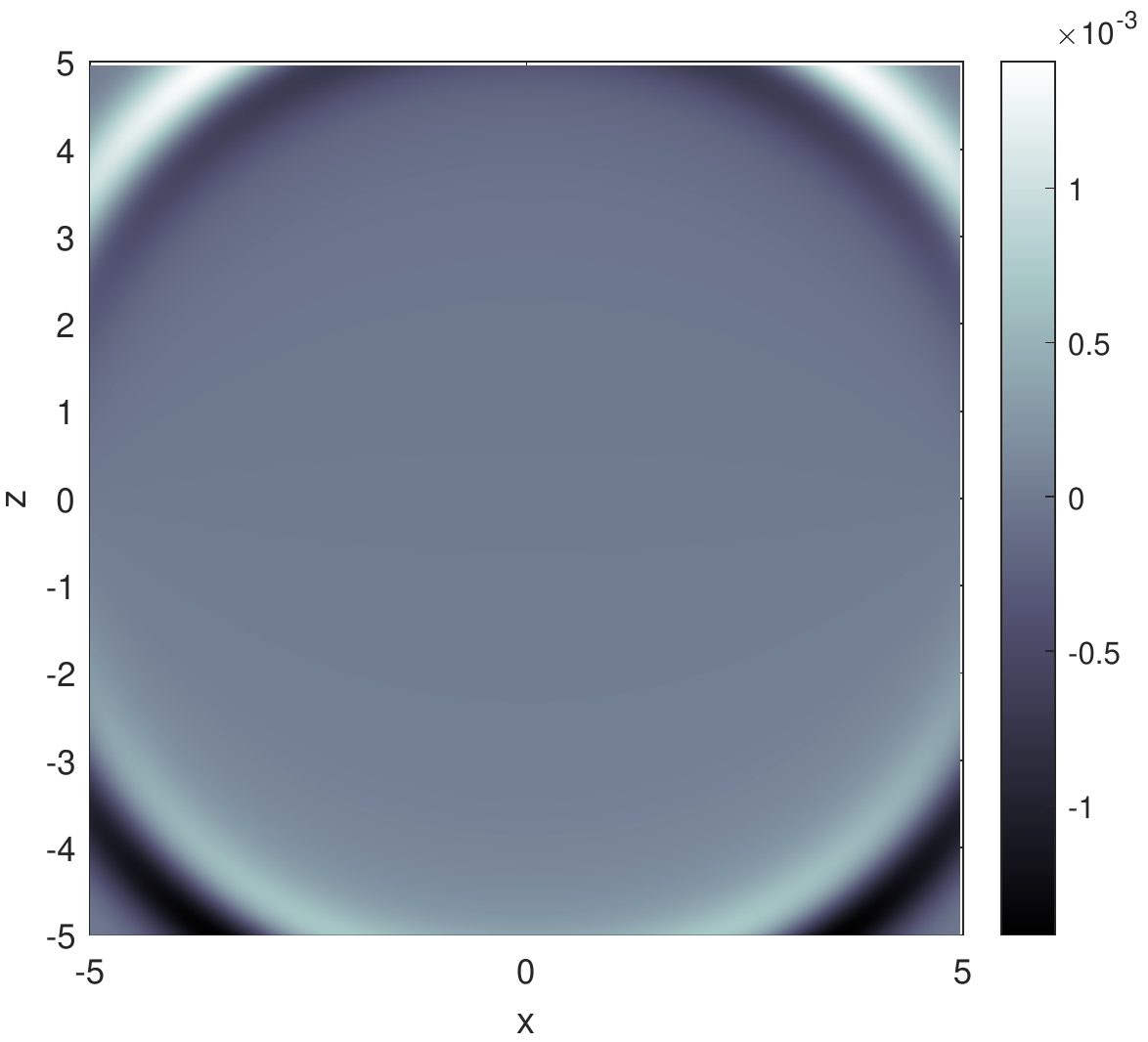}
    \end{minipage}%
    }%
    \subfigure[$t=0.16$s.]{
    \begin{minipage}[t]{0.33\linewidth}
    \centering
    \includegraphics[scale=0.35]{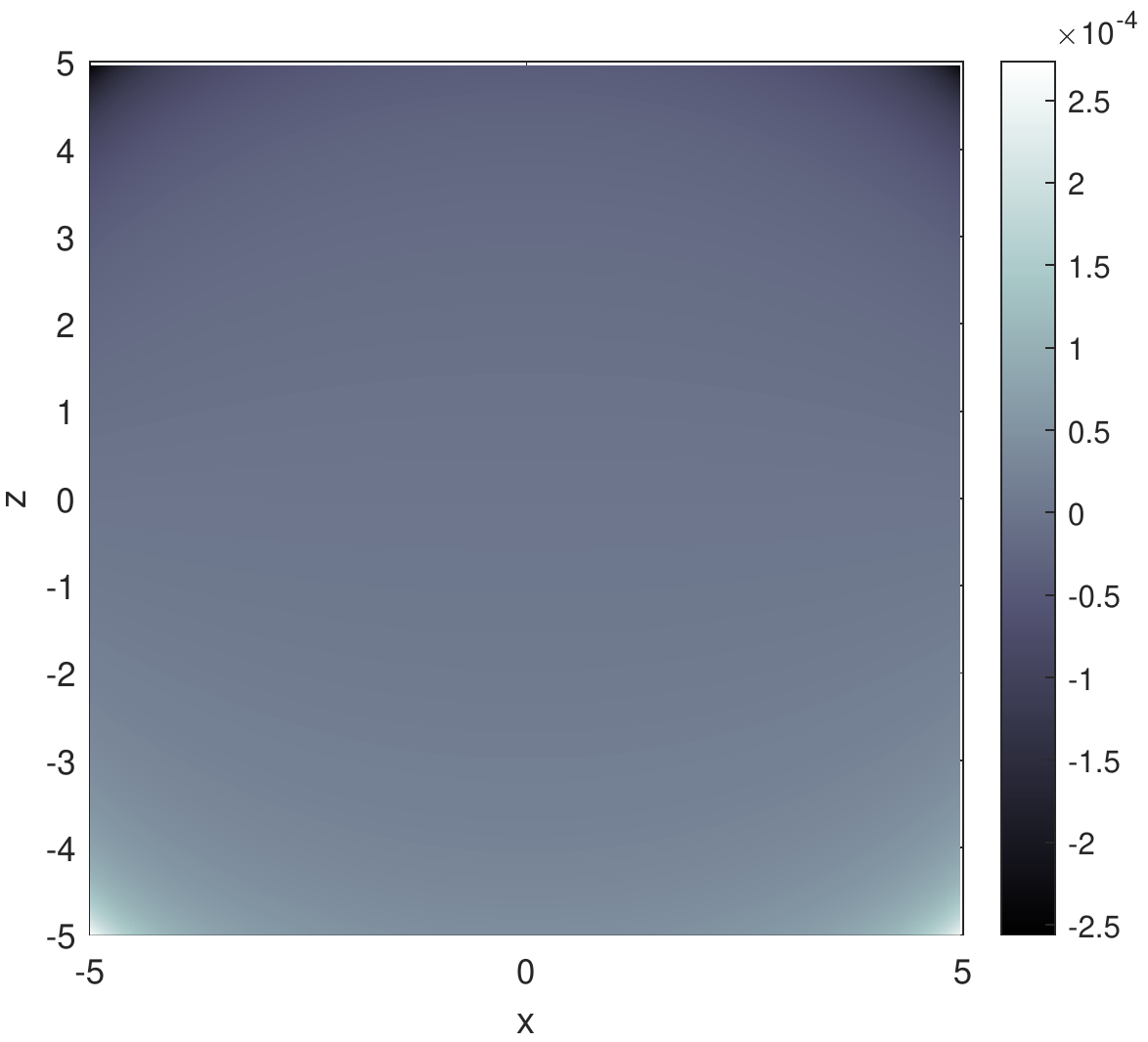}
    \end{minipage}%
    }%
    \subfigure[$t=0.20$s.]{
    \begin{minipage}[t]{0.33\linewidth}
    \centering
    \includegraphics[scale=0.35]{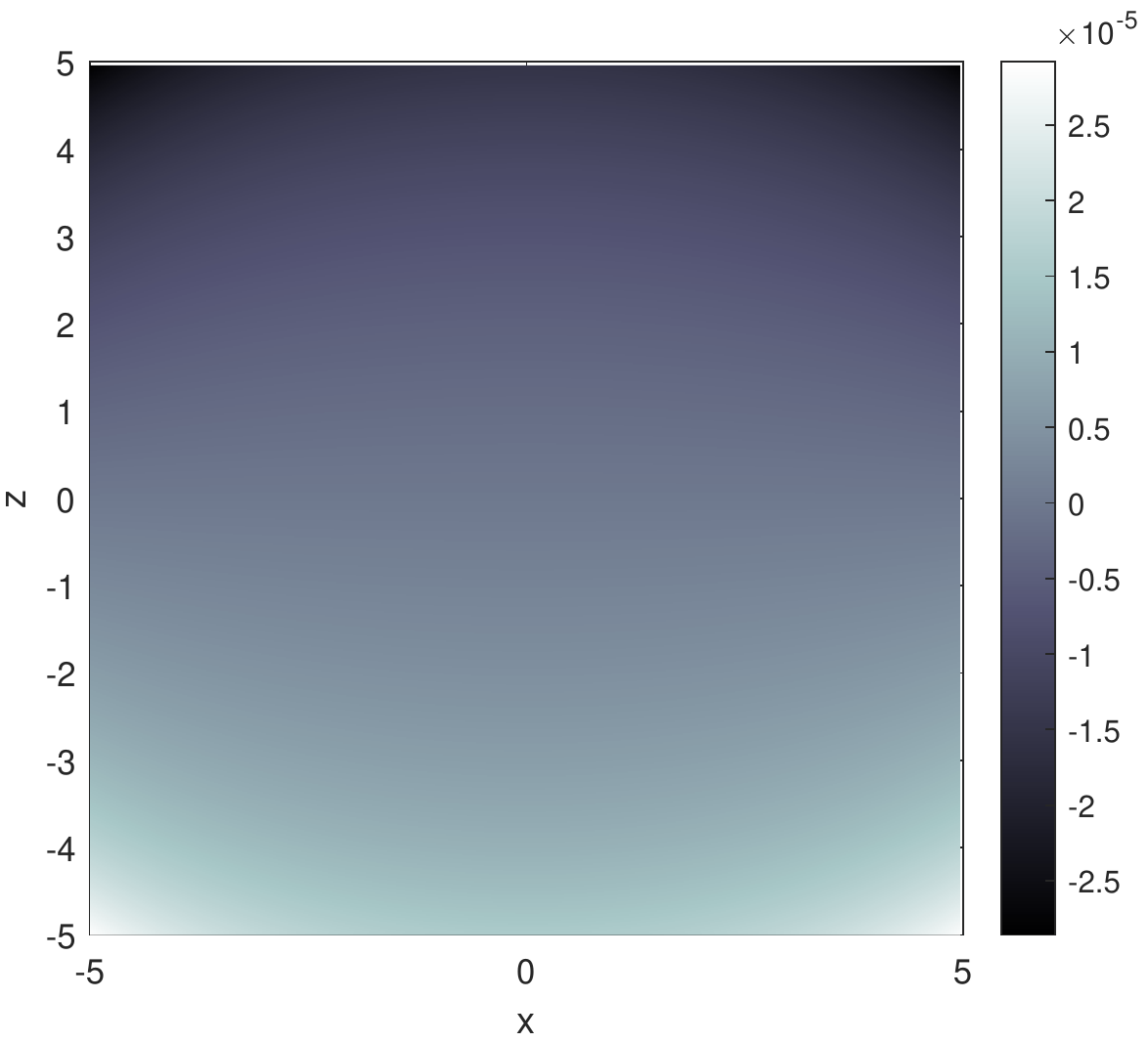}
    \end{minipage}%
    }%
    \centering
    \caption{2-D acoustic wave equation: The snapshot of vibration velocity wavefields $v_3$. The wave begins to be absorbed by PML at $t = 0.1$s. \label{2d_snapshot}}
\end{figure}

\begin{figure}[!h]
    \centering
    \subfigure[$t=0.02$s.]{
    \begin{minipage}[t]{0.33\linewidth}
    \centering
    \includegraphics[scale=0.35]{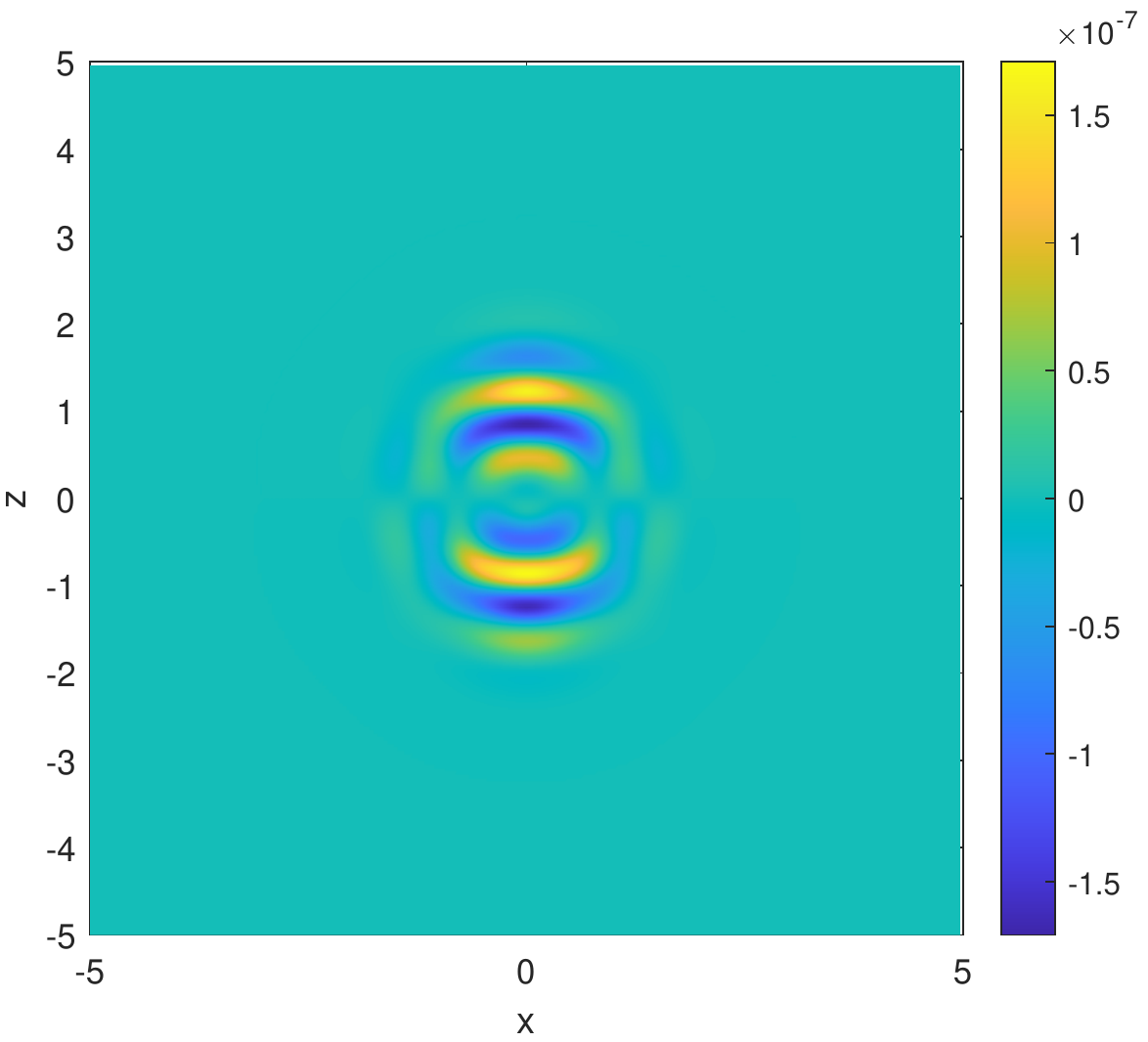}
    \end{minipage}%
    }%
    \subfigure[$t=0.06$s.]{
    \begin{minipage}[t]{0.33\linewidth}
    \centering
    \includegraphics[scale=0.35]{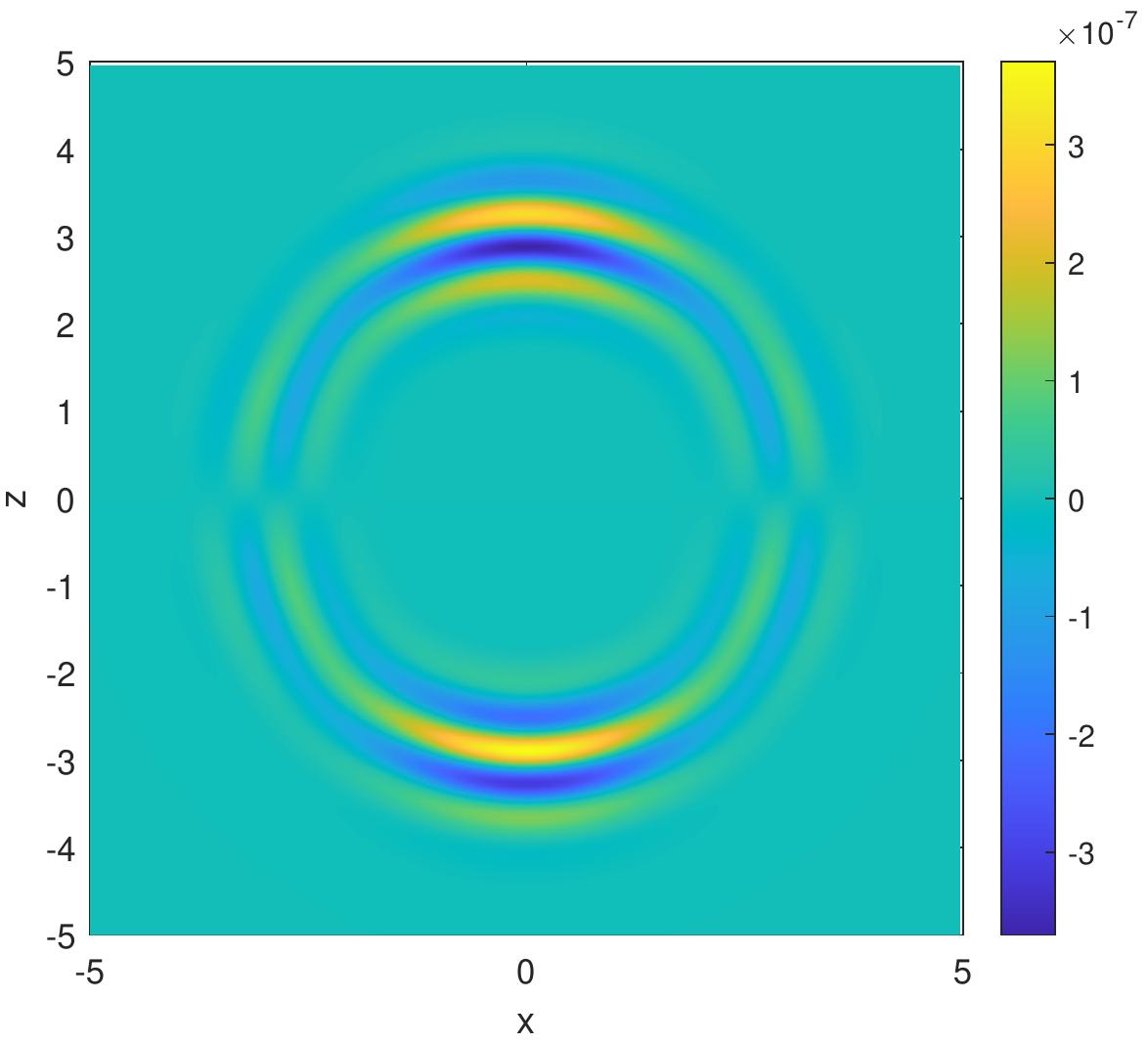}
    \end{minipage}%
    }%
    \subfigure[$t=0.10$s.]{
    \begin{minipage}[t]{0.33\linewidth}
    \centering
    \includegraphics[scale=0.35]{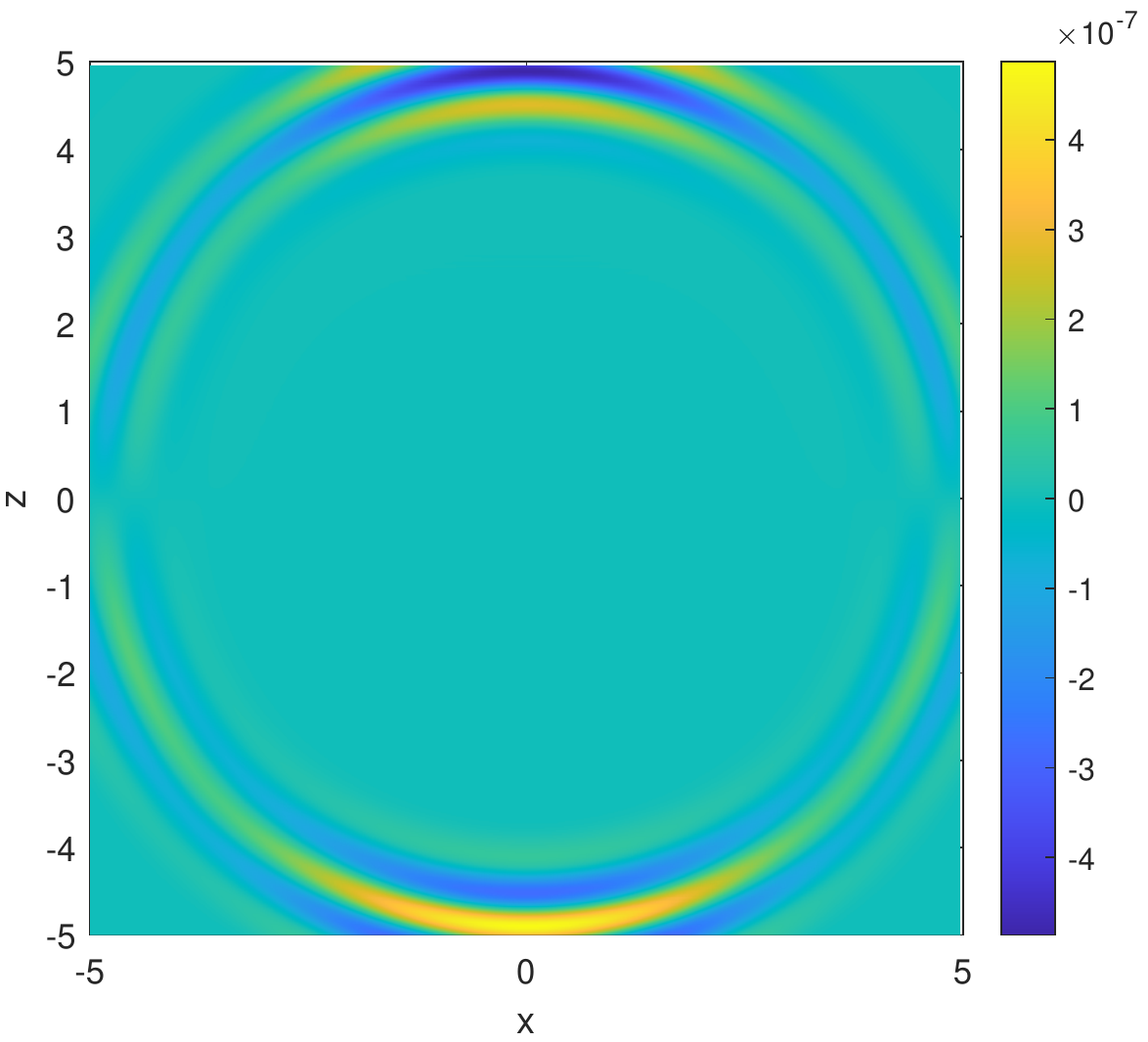}
    \end{minipage}%
    }%

    \subfigure[$t=0.12$s.]{
    \begin{minipage}[t]{0.33\linewidth}
    \centering
    \includegraphics[scale=0.35]{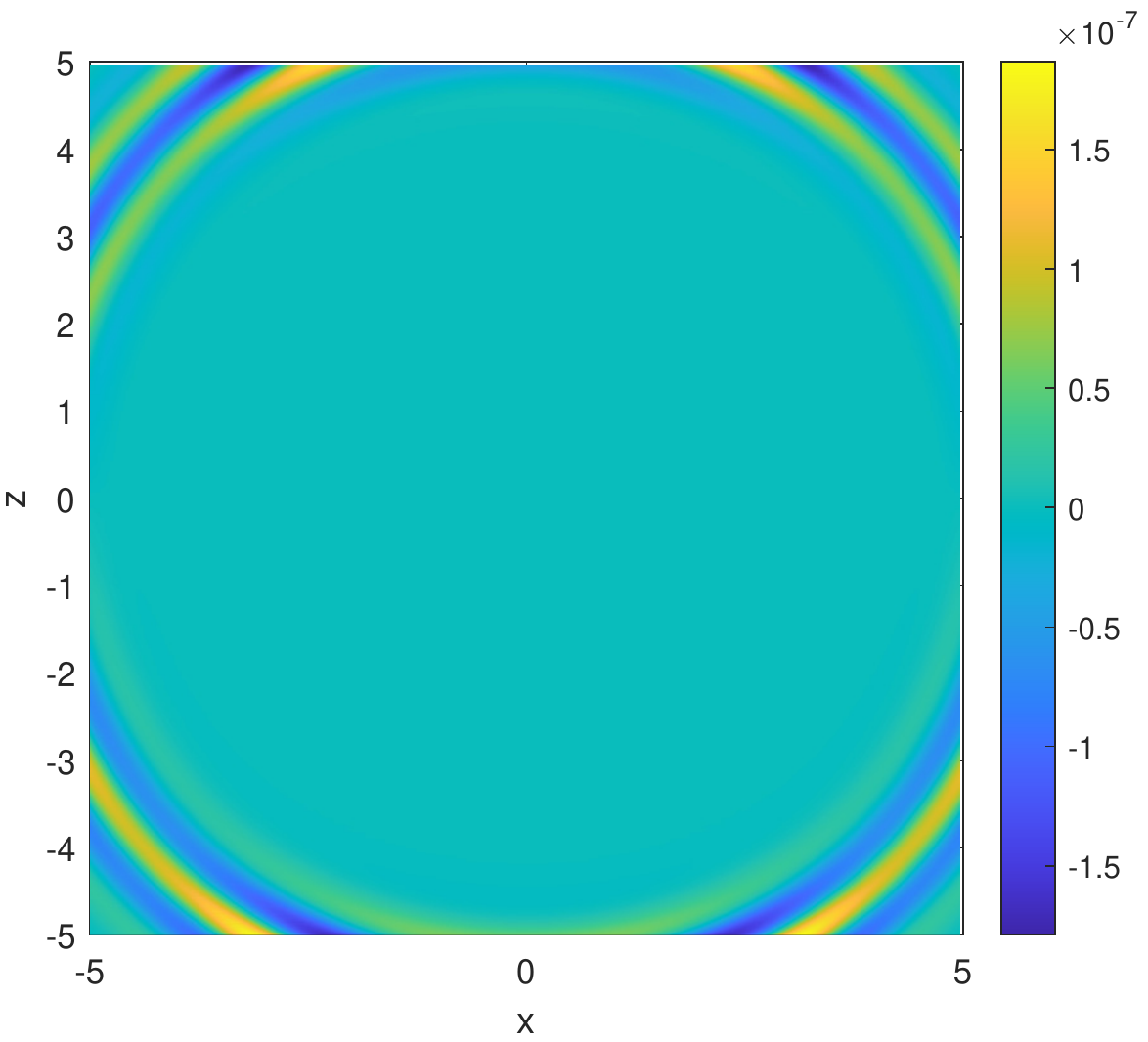}
    \end{minipage}%
    }%
    \subfigure[$t=0.16$s.]{
    \begin{minipage}[t]{0.33\linewidth}
    \centering
    \includegraphics[scale=0.35]{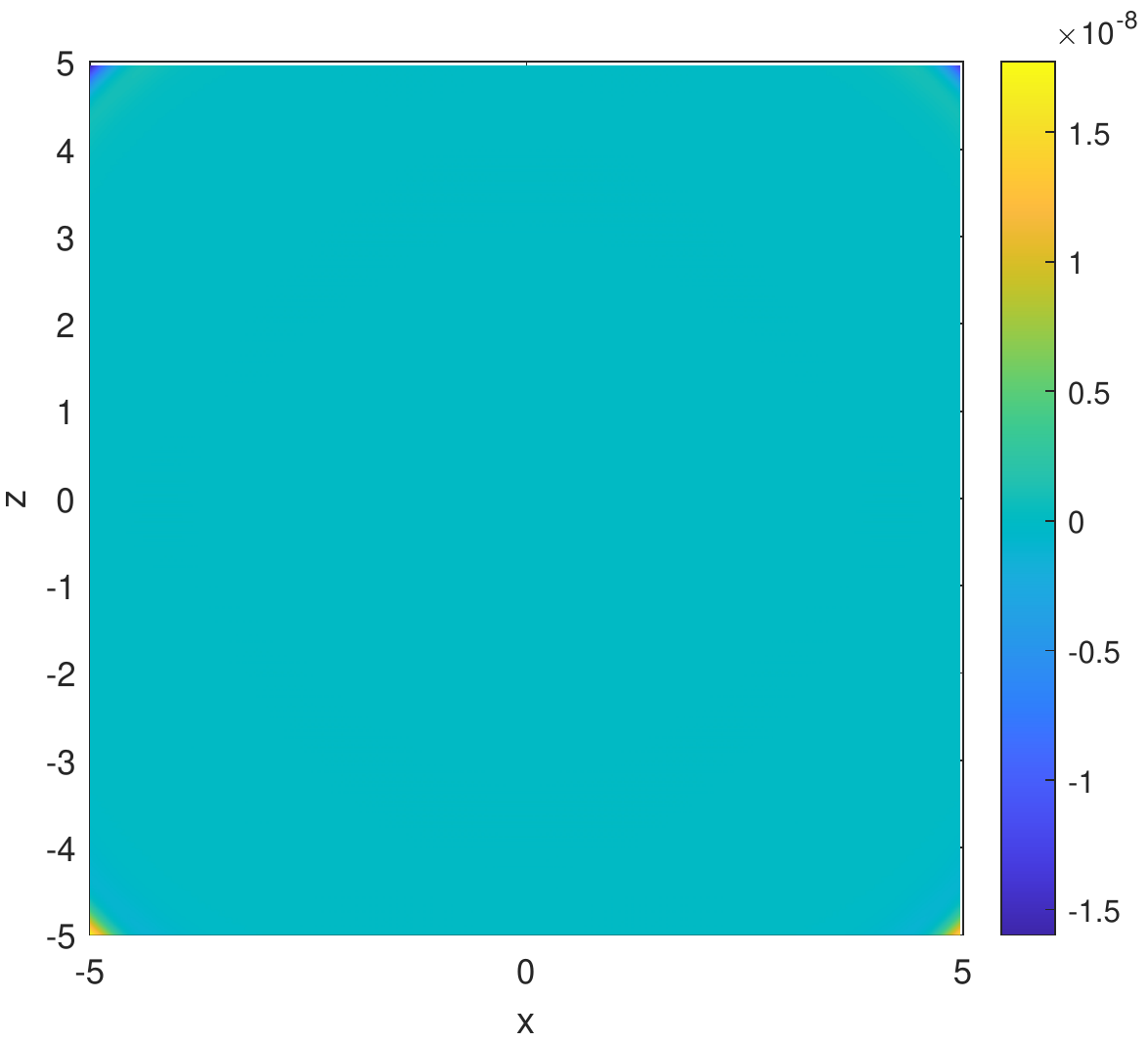}
    \end{minipage}%
    }%
    \subfigure[$t=0.20$s.]{
    \begin{minipage}[t]{0.33\linewidth}
    \centering
    \includegraphics[scale=0.35]{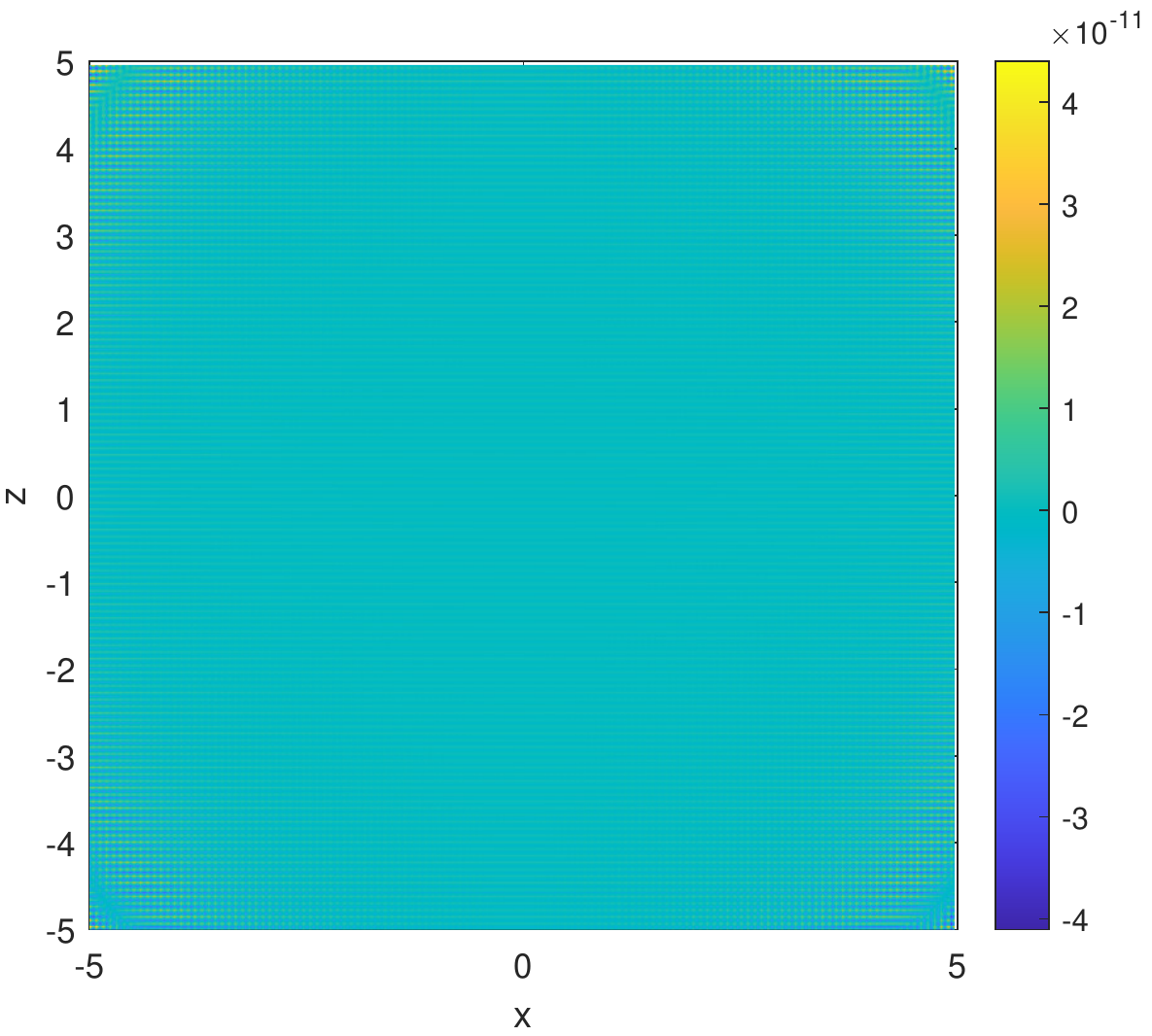}
    \end{minipage}%
    }%
    \centering
    \caption{2-D acoustic wave equation: Visualization of numerical errors in the vibrational velocity wavefields $v_3$. \label{2d_snapshot_error}}
\end{figure}

The snapshots of vibrational velocity wavefield $v_3$ and the distribution of numerical errors are visualized in Figures \ref{2d_snapshot} and \ref{2d_snapshot_error}, respectively. From Figure \ref{2d_snapshot}, the velocity wavefield propagates inside the domain before $t = 0.1$s, and it begins to leave the domain when $t > 0.1$s. Fortunately, the penetrating wavefields are successfully attenuated by ADE-PML and the artificial reflection is almost negligible, as clearly observed in  Figure \ref{2d_snapshot_error}.

\begin{figure}[!h]
    \centering
    \subfigure[When wave propagates in the domain.]{\includegraphics[scale=0.5]{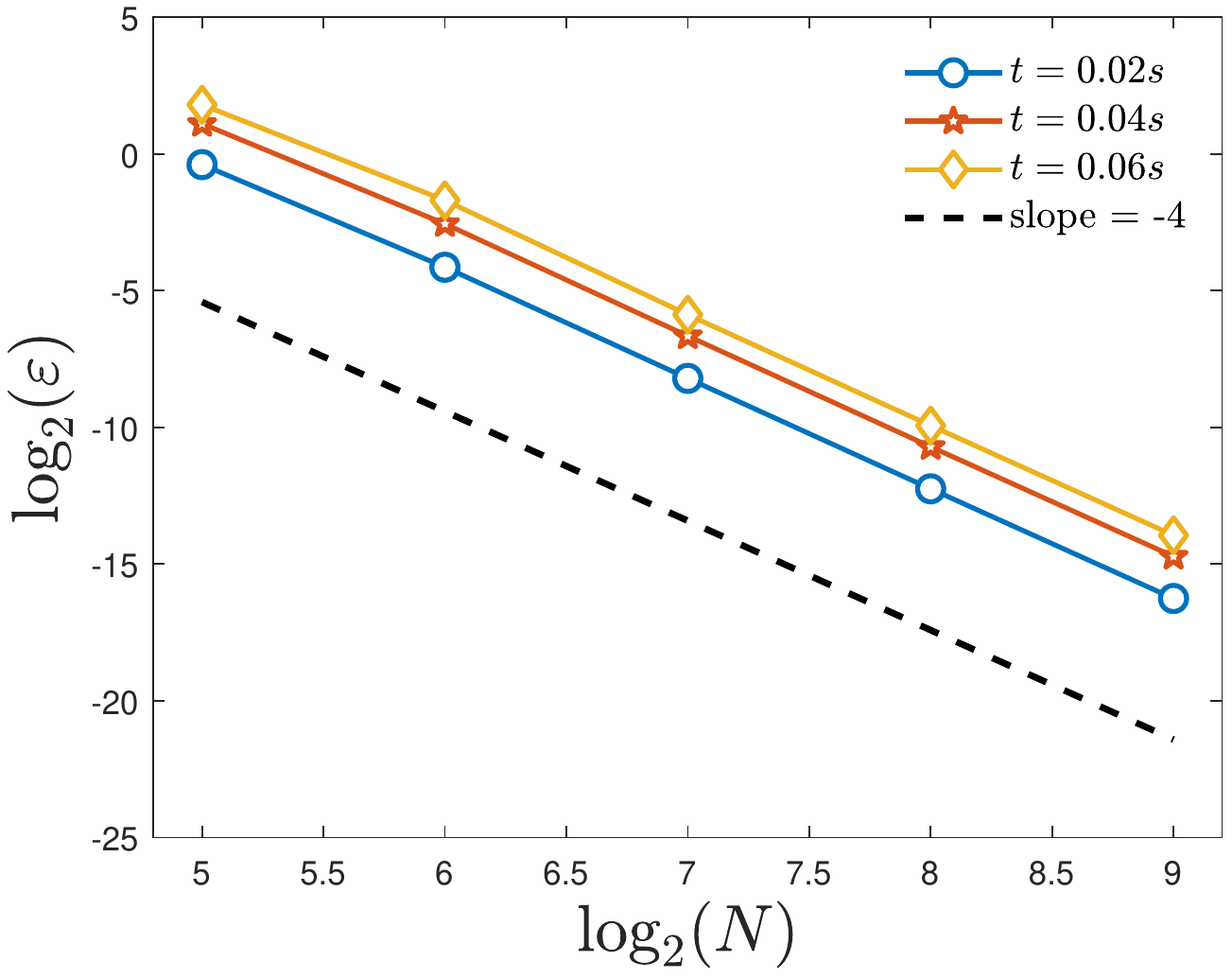}}
    \subfigure[When wave leaves the domain.]{\includegraphics[scale=0.5]{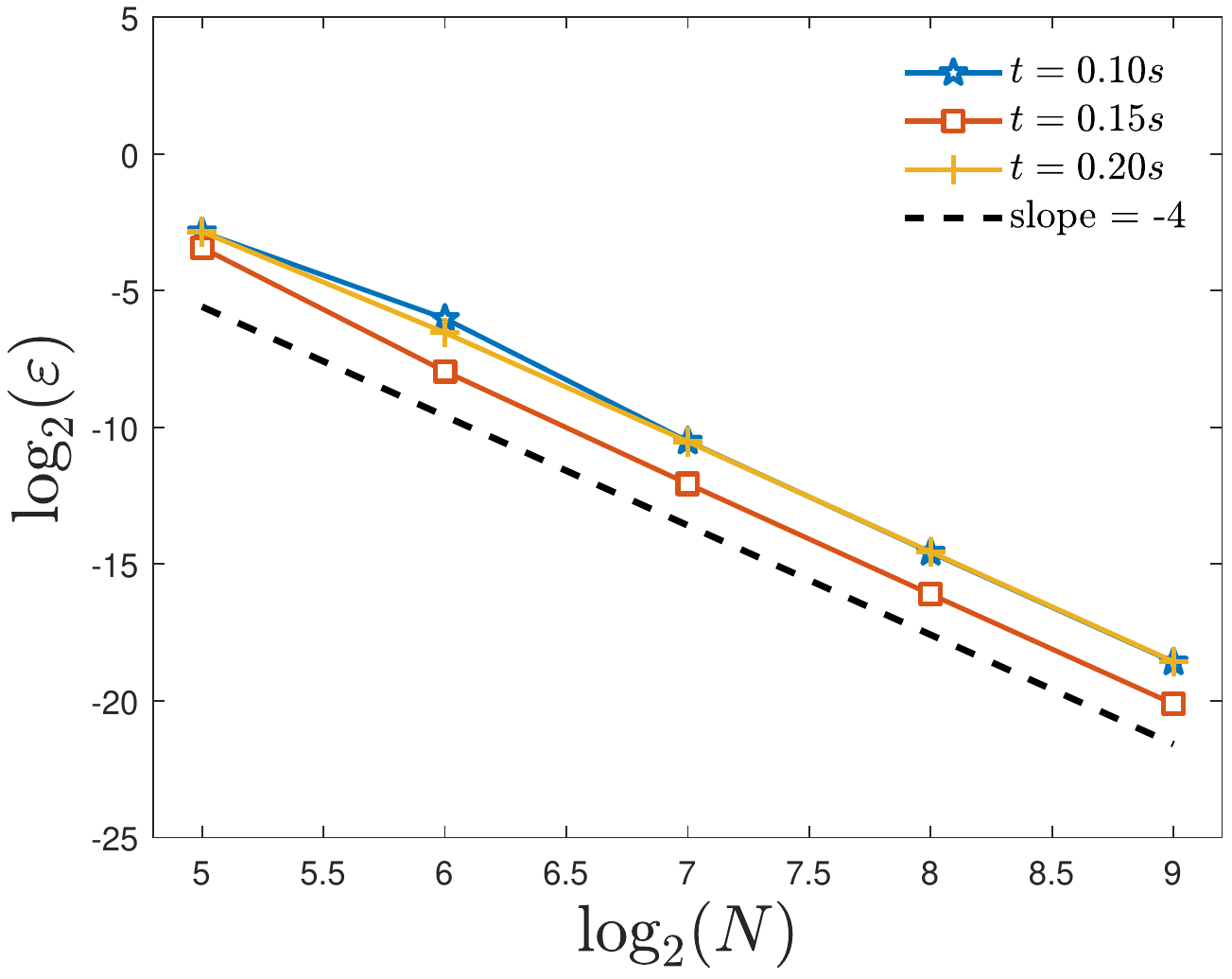}}
    \caption{2-D acoustic wave equation: The fourth-order convergence of the spline collocation method {is verified for ADE-PML}. \label{2d_error_convergence}}
\end{figure}

The time evolution of $l^2$-errors $\varepsilon_2(t)$ are recorded in Table \ref{2d_error} and the convergence rate is given in Figure \ref{2d_error_convergence}. The slope of the dashed line is $-4$, which perfectly matches the theoretical fourth-order convergence. In addition, according to Table \ref{2d_error} and Figure \ref{2d_error_convergence},  the accuracy of ADE-PML can be significantly improved under a finer grid, and the convergence rate is close to $4$. This indicates the performance of ADE-PML depends heavily on the accuracy of the underlying numerical solver.

 \par

\begin{table}[!h]
        \caption{2-D acoustic wave equation: The $l^2$-error $\varepsilon_2(t)$ at different instants. The numerical accuracy of the spline collocation method and ADE-PML can be systematically improved by refining the grid mesh.   \label{2d_error}}
            \centering
    \resizebox{1\columnwidth}{!}{
    \begin{tabular}{|c|c|c|c|c|c|c}
    \hline
    \diagbox{Time}{Grid}&$N=32^2$&$N=64^2$&$N=128^2$&$N=256^2$&$N=512^2$\\ 
    \hline
    \multicolumn{6}{|c|}{When wave propagates inside the domain}\\
    \hline
    $0.02$s		&7.634 $\times 10^{-1}$&5.663 $\times 10^{-2}$&3.382 $\times 10^{-3}$&2.057 $\times 10^{-4}$&1.276 $\times 10^{-5}$
    \\
    \hline
    $0.04$s		&3.634$\times10^{-1}$&2.197$\times10^{-2}$&1.264$\times10^{-3}$&7.640$\times10^{-5}$&4.737$\times10^{-6}$
    \\
    \hline
    $0.06$s		&4.476$\times10^{-1}$&3.638$\times10^{-2}$&2.167$\times10^{-3}$&1.309$\times10^{-4}$&8.118$\times10^{-6}$
    \\
    \hline
    \multicolumn{6}{|c|}{When wave is absorbed by PML outside the domain}\\
    \hline
    $0.10$s		&2.413$\times10^{-1}$&3.099$\times10^{-2}$&1.469$\times10^{-3}$&8.802$\times10^{-5}$&5.466$\times10^{-6}$
    \\
    \hline
    $0.15$s		&2.186$\times10^{-1}$&1.038$\times10^{-2}$&3.837$\times10^{-4}$&2.329$\times10^{-5}$&1.445$\times10^{-6}$
    \\
    \hline
    $0.20$s		&2.441$\times10^{-1}$&1.779$\times10^{-2}$&1.080$\times10^{-3}$&6.611$\times10^{-5}$&4.238$\times10^{-6}$
    \\
    \hline
    \end{tabular}
    }
\end{table}

\subsection{The performance of ADE-PML under spline collocation method}

As mentioned above, the performance of ADE-PML relies on three parameters: the thickness $L$, the reflection coefficient $R$ and the cut-off wavenumber $k_{\max}$ \cite{KomatitschMartin2007,MartinKomatitsch2009,MartinKomatitschGedneyBruthiaux2010}. Therefore, their influence on the absorbing effect, as well as the compatibility with natural boundary conditions for the cubic spline,  deserves a careful investigation.

Among them, the most important parameter is the layer thickness $L$. Intuitively speaking, increasing the layer thickness improves the absorbing effect of PML, at the cost of storing more wavefields and higher computational costs. Therefore, it is expected to use absorbing layers as fewer as possible to maintain the accuracy.

\begin{figure}[!h]
    \centering
    \subfigure[t=0.1s.]{
    \begin{minipage}[t]{0.33\linewidth}
    \centering
    \includegraphics[scale=0.35]{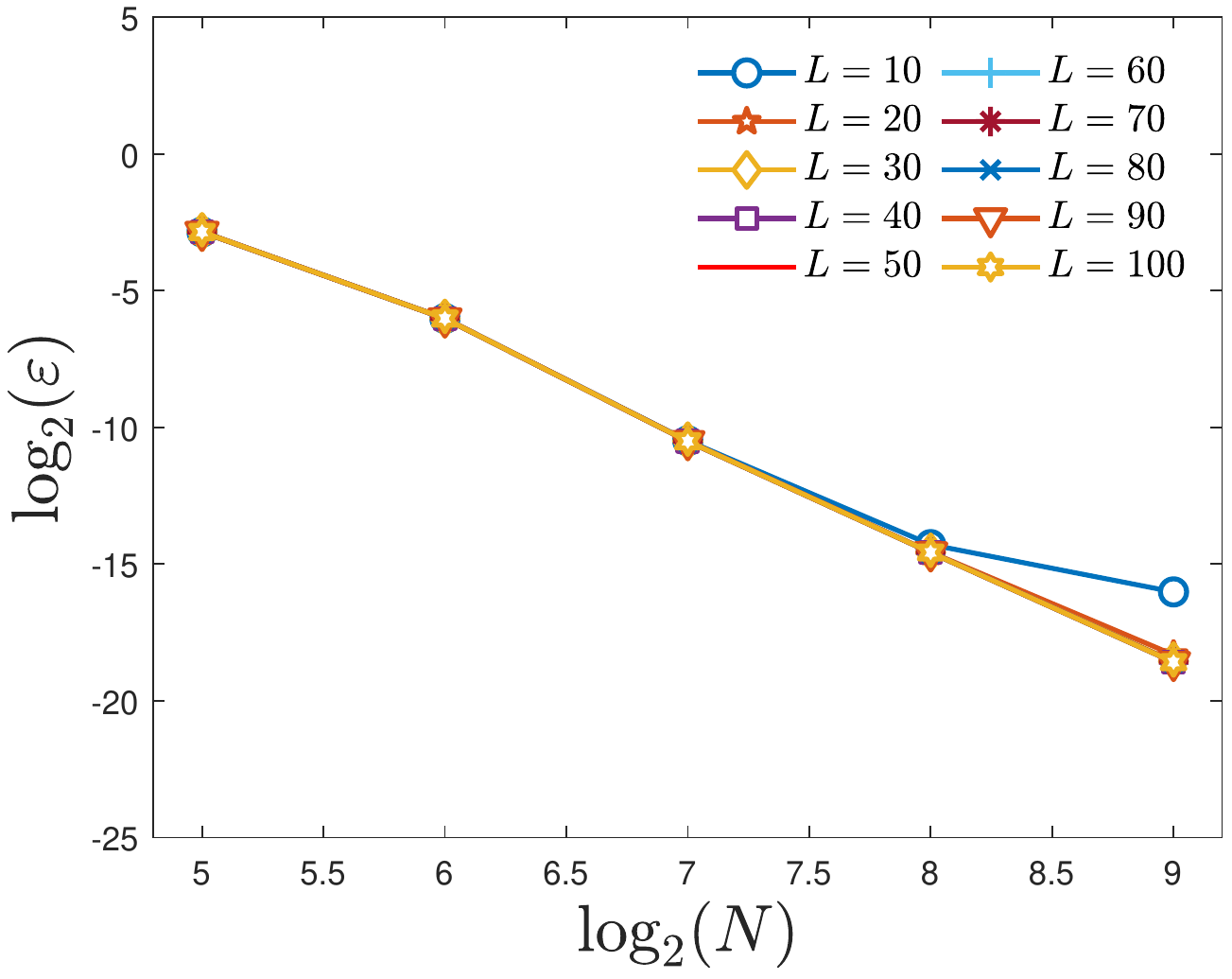}
    \end{minipage}%
    }%
    \subfigure[t=0.15s.]{
    \begin{minipage}[t]{0.33\linewidth}
    \centering
    \includegraphics[scale=0.35]{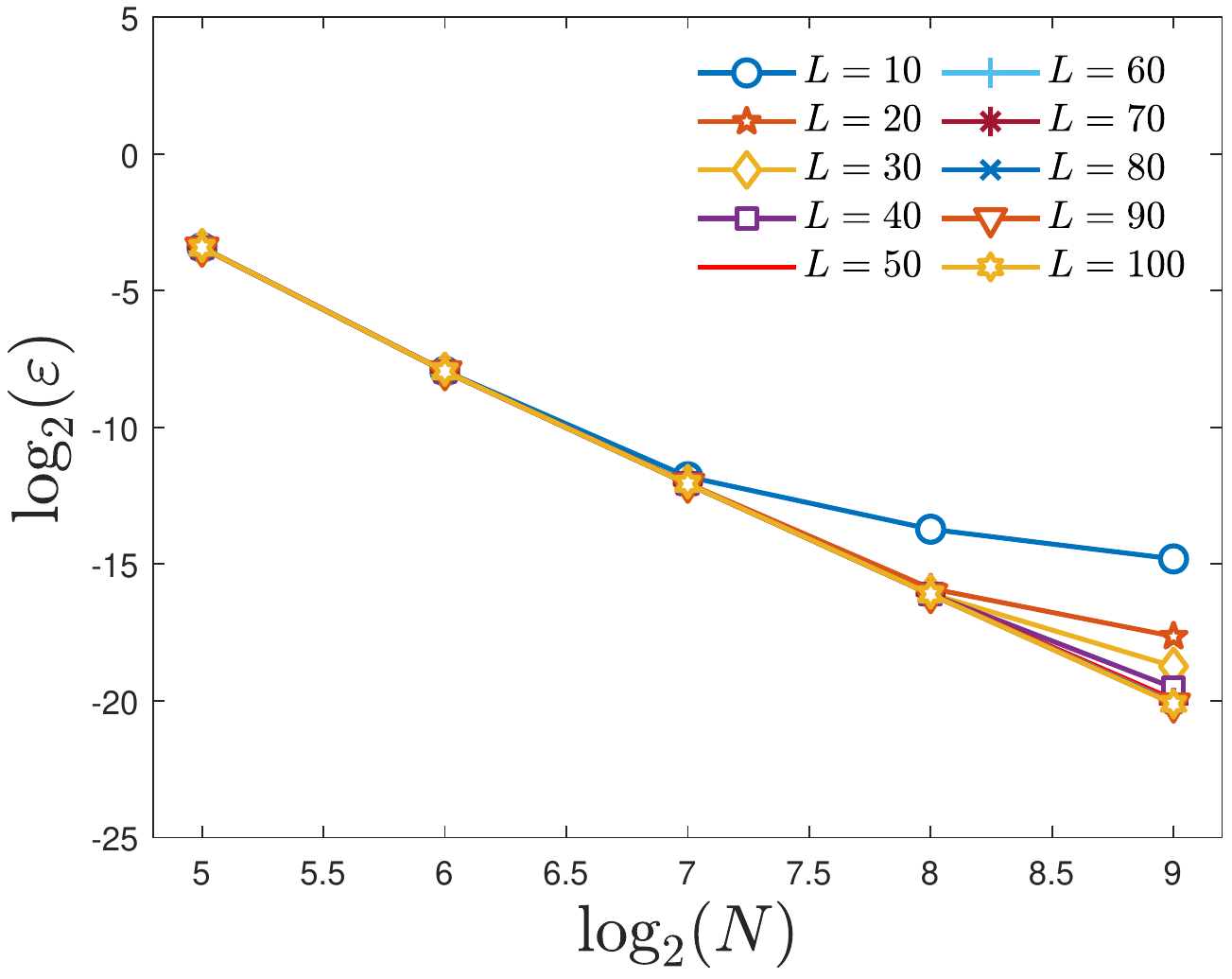}
    \end{minipage}%
    }%
    \subfigure[t=0.20s.]{
    \begin{minipage}[t]{0.33\linewidth}
    \centering
    \includegraphics[scale=0.35]{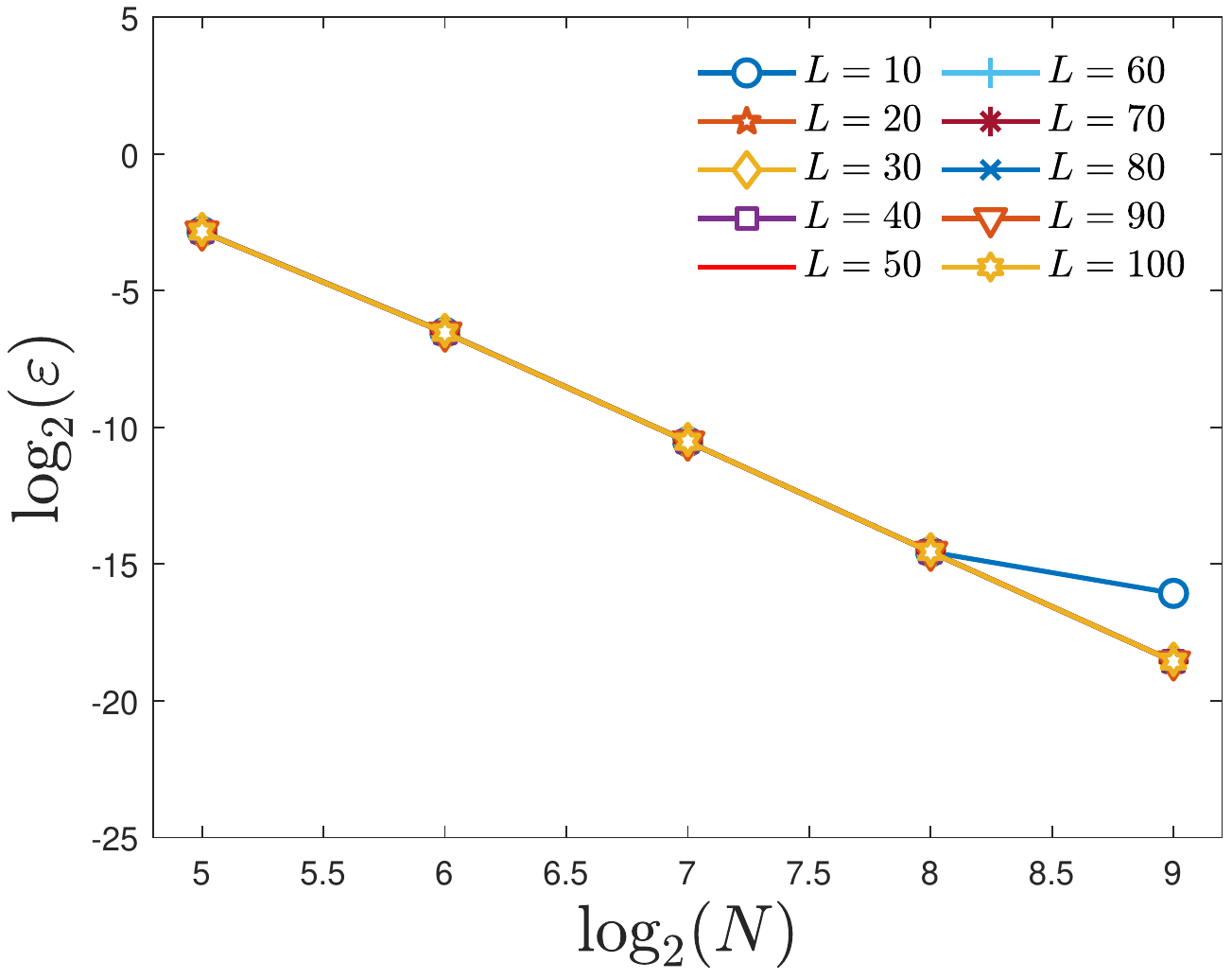}
    \end{minipage}%
    }%
    \centering
    \caption{2-D PML: The convergence of $l^2$-errors  under different boundary layer thicknesses $L$. It is observed that $L=40$-$50$ might be a good choice to strike a balance in accuracy and efficiency. \label{PML_error_L_Nx}}
\end{figure}

To evaluate its effect, we have performed a series of benchmarks under different $L$ and grids. The $l^2$-errors at $t = 0.2$s,  as recorded in Table \ref{PML_error_L_t020}, were adopted to measure the accuracy of ADE-PML. The convergence with respect to the mesh size is plotted in Figure \ref{PML_error_L_Nx}. Clearly, the accuracy of ADE-PML can be systematically improved as $N$ increases, and the convergence rate matches the theoretical order $4$ when $L$ is sufficiently large.
  \par

\begin{table}[!h]
    \caption{2-D PML: The $l^2$-errors at $t = 0.2$s under different boundary layer thicknesses $L$ and grid meshes.  The errors induced by PML are negligible when $L \ge 50$. \label{PML_error_L_t020}}
        \centering
    \resizebox{1.\columnwidth}{!}{
    \begin{tabular}{|c|c|c|c|c|c|}
    \hline
    \diagbox{$L$}{Grid}&$N=32^2$&$N=64^2$&$N=128^2$&$N=256^2$&$N=512^2$\\ 
    \hline
    $10$&1.405$\times10^{-1}$&1.089$\times10^{-2}$&6.785$\times10^{-4}$&4.155$\times10^{-5}$&1.453$\times10^{-5}$

    \\
    \hline
    $20$&1.394$\times10^{-1}$&1.077$\times10^{-2}$&6.785$\times10^{-4}$&4.155$\times10^{-5}$&2.593$\times10^{-6}$

    \\
    \hline
    $30$&1.394$\times10^{-1}$&1.077$\times10^{-2}$&6.787$\times10^{-4}$&4.155$\times10^{-5}$&2.587$\times10^{-6}$

    \\
    \hline
    $40$&1.394$\times10^{-1}$&1.077$\times10^{-2}$&6.786$\times10^{-4}$&4.154$\times10^{-5}$&2.581$\times10^{-6}$

    \\
    \hline
    $50$&1.394$\times10^{-1}$&1.077$\times10^{-2}$&6.786$\times10^{-4}$&4.157$\times10^{-5}$&2.580$\times10^{-6}$

    \\
    \hline
    $60$&1.394$\times10^{-1}$&1.077$\times10^{-2}$&6.786$\times10^{-4}$&4.156$\times10^{-5}$&2.581$\times10^{-6}$

    \\
    \hline
    $70$&1.394$\times10^{-1}$&1.077$\times10^{-2}$&6.786$\times10^{-4}$&4.154$\times10^{-5}$&2.579$\times10^{-6}$

    \\
    \hline
    $80$&1.394$\times10^{-1}$&1.077$\times10^{-2}$&6.786$\times10^{-4}$&4.155$\times10^{-5}$&2.577$\times10^{-6}$

    \\
    \hline
    $90$&1.394$\times10^{-1}$&1.077$\times10^{-2}$&6.786$\times10^{-4}$&4.155$\times10^{-5}$&2.578$\times10^{-6}$

    \\
    \hline
    $100$&1.394$\times10^{-1}$&1.077$\times10^{-2}$&6.786$\times10^{-4}$&4.155$\times10^{-5}$&2.580$\times10^{-6}$

    \\
    \hline
    \end{tabular}
    }
\end{table}
 
 It is observed that too small $L$ may result in a dramatic reduction in accuracy for sufficiently large $N$, while its influence becomes negligible when $N_x = N_z \le 128$ because the discretization errors dominate. Therefore, the choice of $L$ should match the size of grid mesh. When a coarse grid is used, the thickness might be not too large for the sake of efficiency. By contrast, when a fine grid mesh is used, one should use a thicker absorbing layer to ensure the accuracy. Based on the results in Figure \ref{PML_error_L_Nx}, $L=40$ to $50$ might be a good choice to strike a balance in accuracy and efficiency, while too large $L$ might not bring in evident improvements as the accuracy is still limited by the numerical solver under the specified grid mesh.

\begin{figure}[!h]
    \centering
   \subfigure[Different reflection coefficients $R$.\label{err_evo_PML_R}]{ \includegraphics[width=0.47\textwidth,height=0.25\textwidth]{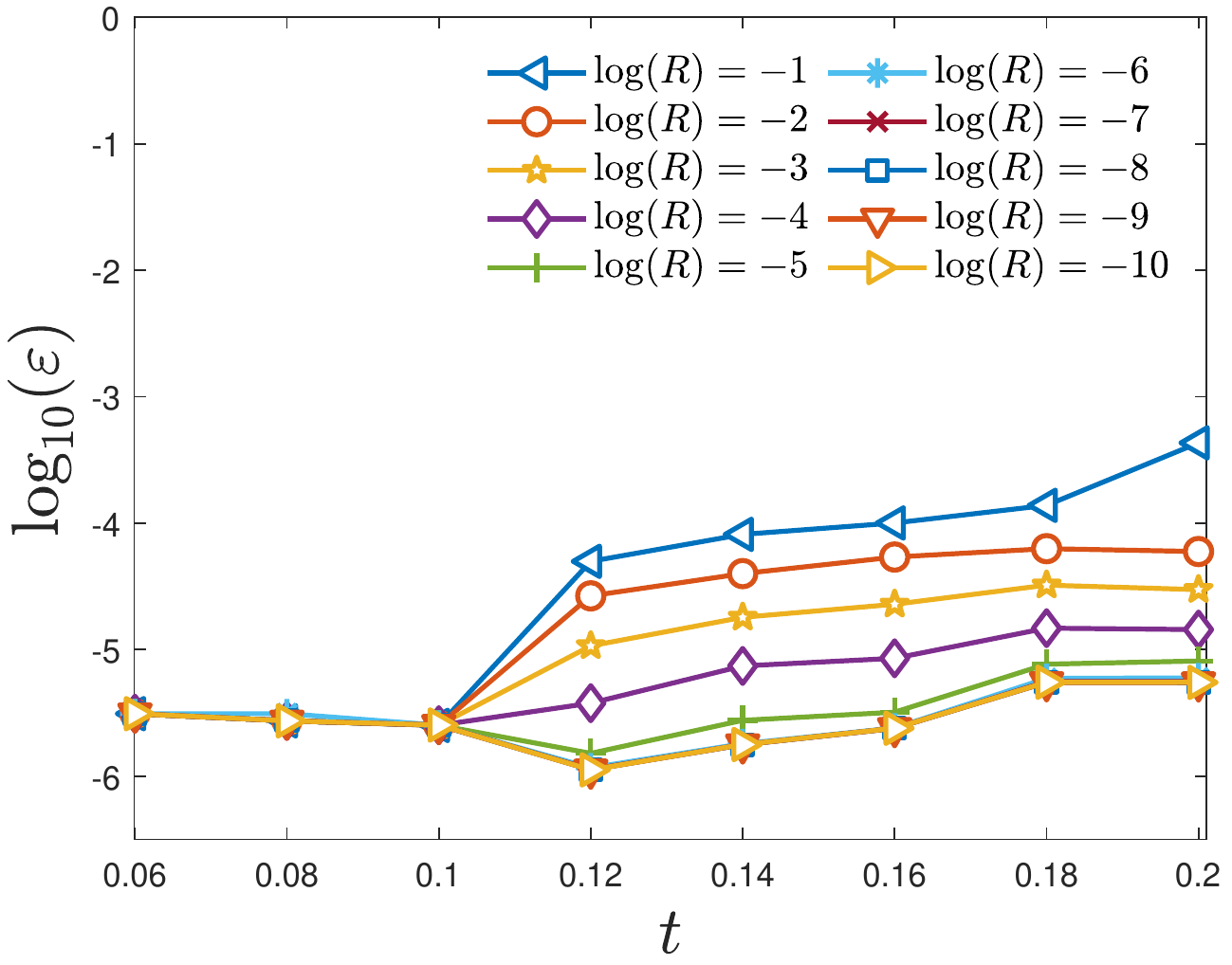}}
   \subfigure[Different cut-off wavenumber $k_{\max}$.\label{err_evo_PML_kmax}]{ \includegraphics[width=0.47\textwidth,height=0.25\textwidth]{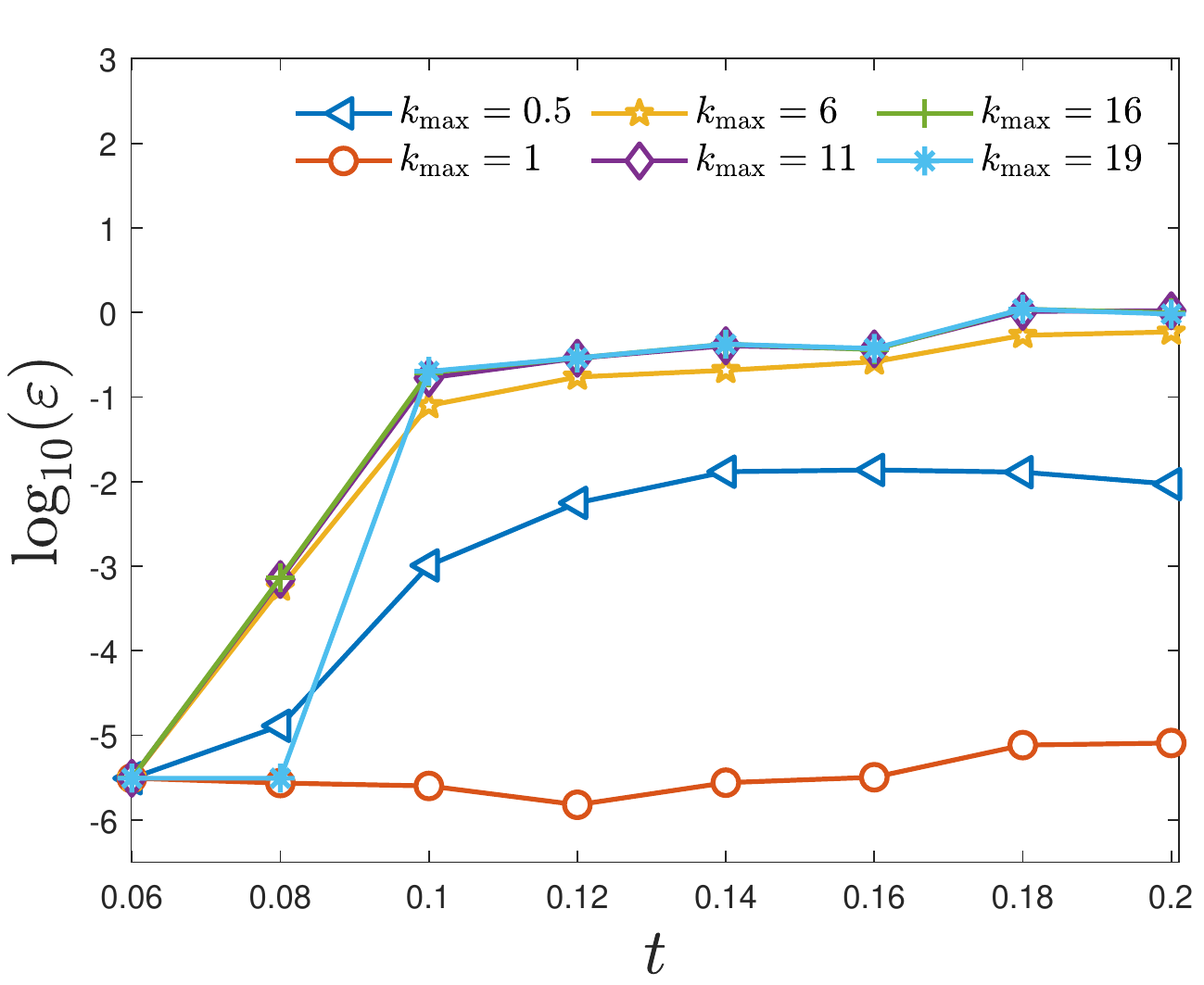}}
    \caption{2-D PML: Time evolution of $l^2$-errors under different reflection coefficients $R$ and cut-off wavenumber $k_{\max}$. The mesh size is fixed to be $N = 512^2$. From numerical perspective, The choice of $R$ should match the thickness $L$ of PML, and cut-off wavenumber $k_{\max}= 1$ achieves the optimal performance.}
\end{figure}

Second, we studied the effect of the reflection coefficient $R$ on the absorbing boundary. The mesh size was fixed to be $N_x \times N_z = 512^2$ with the layer thickness $L = 50$. The time evolution of the $l^2$-errors are plotted in Figure \ref{err_evo_PML_R}, corresponding to the data in Table \ref{err_PML_R}. It is observed that $R = 10^{-6}$ achieves $l^2$-error about $10^{-6}$, whereas it seems to reach the limitation of accuracy as the errors vary slightly when $R$ further decreases. These findings coincide with the observation made in  \cite{CollinoTsogka2001,MartinKomatitschGedneyBruthiaux2010} as the choice of $R$ should match the thickness $L$ of PML.

\begin{table}[!htbp]
    \caption{2-D PML: The $l^2$-errors under different reflection coefficients $R$, with $N_x = N_z = 512$. Here $R = 10^{-6}$ is suggested to ensure the effectiveness of PML.\label{err_PML_R}}
    \centering
    \resizebox{1\columnwidth}{!}{
    \begin{tabular}{|c|c|c|c|c|c|c|c|c|}
    \hline
    \diagbox{$ R$}{Time}&$0.06$s&$0.08$s&$0.10$s&$0.12$s&$0.14$s&$0.16$s&$0.18$s&$0.20$s\\ 
    \hline
    $10^{-1}$&3.11$\times10^{-6}$&2.73$\times10^{-6}$&2.52$\times10^{-6}$&4.99$\times10^{-5}$&8.15$\times10^{-5}$&1.00$\times10^{-4}$&1.38$\times10^{-4}$&4.32$\times10^{-4}$
    \\
    \hline
    $10^{-2}$&3.11$\times10^{-6}$&2.73$\times10^{-6}$&2.52$\times10^{-6}$&2.66$\times10^{-5}$&3.99$\times10^{-5}$&5.40$\times10^{-5}$&6.28$\times10^{-5}$&5.97$\times10^{-5}$
    \\
    \hline
    $10^{-3}$&3.11$\times10^{-6}$&2.73$\times10^{-6}$&2.52$\times10^{-6}$&1.07$\times10^{-5}$&1.80$\times10^{-5}$&2.28$\times10^{-5}$&3.24$\times10^{-5}$&2.98$\times10^{-5}$
    \\
    \hline
    $10^{-4}$&3.11$\times10^{-6}$&2.73$\times10^{-6}$&2.52$\times10^{-6}$&3.77$\times10^{-6}$&7.43$\times10^{-6}$&8.54$\times10^{-6}$&1.48$\times10^{-5}$&1.44$\times10^{-5}$
    \\
    \hline
    $10^{-5}$&3.11$\times10^{-6}$&2.73$\times10^{-6}$&2.52$\times10^{-6}$&1.51$\times10^{-6}$&2.76$\times10^{-6}$&3.20$\times10^{-6}$&7.67$\times10^{-6}$&8.13$\times10^{-6}$
    \\
    \hline
    $10^{-6}$&3.11$\times10^{-6}$&3.11$\times10^{-6}$&2.53$\times10^{-6}$&1.16$\times10^{-6}$&1.80$\times10^{-6}$&2.40$\times10^{-6}$&5.92$\times10^{-6}$&5.98$\times10^{-6}$
    \\
    \hline
    $10^{-7}$&3.11$\times10^{-6}$&2.73$\times10^{-6}$&2.53$\times10^{-6}$&1.13$\times10^{-6}$&1.77$\times10^{-6}$&2.39$\times10^{-6}$&5.58$\times10^{-6}$&5.58$\times10^{-6}$
    \\
    \hline
    $10^{-8}$&3.11$\times10^{-6}$&2.73$\times10^{-6}$&2.53$\times10^{-6}$&1.12$\times10^{-6}$&1.77$\times10^{-6}$&2.39$\times10^{-6}$&5.51$\times10^{-6}$&5.48$\times10^{-6}$
    \\
    \hline
    $10^{-9}$&3.11$\times10^{-6}$&2.73$\times10^{-6}$&2.53$\times10^{-6}$&1.12$\times10^{-6}$&1.77$\times10^{-6}$&2.39$\times10^{-6}$&5.51$\times10^{-6}$&5.48$\times10^{-6}$
    \\
    \hline
    $10^{-10}$&3.11$\times10^{-6}$&2.73$\times10^{-6}$&2.53$\times10^{-6}$&1.12$\times10^{-6}$&1.78$\times10^{-6}$&2.39$\times10^{-6}$&5.53$\times10^{-6}$&5.51$\times10^{-6}$
    \\
    \hline
    \end{tabular}
    }
\end{table}

Finally, for the cut-off wavenumber $k_{\max}$, we also  plot the time evolution of the $l^2$-errors under different $k_{\max}$ in Figure \ref{err_evo_PML_kmax}, with data collected in Table \ref{err_PML_kmax}. The mesh size was fixed to be $N_x \times N_z = 512^2$ and the layer thickness was again set as $L = 50$. The results show that ADE-PML can work only when $k_{\max} = 1$. In fact, just as pointed out in \cite{MartinKomatitsch2009}, the sharp variations of the profile of the $k_x$ and $k_z$ functions might augment the reflection coefficient of waves impinging on the boundary, so that too large $k_{\max}$ is not recommended.

\begin{table}[!h]
    \caption{2-D PML: The $l^2$-errors under different cut-off wavenumber $k_{\max}$, with $N_x = N_z = 512$. The ADE-PML may only work when  $k_{\max} = 1$. \label{err_PML_kmax}}
    \centering
    \resizebox{1.\columnwidth}{!}{
    \begin{tabular}{|c|c|c|c|c|c|c|c|c|}
    \hline
    \diagbox{$k_{\max}$}{Time}&$0.06$s&$0.08$s&$0.10$s&$0.12$s&$0.14$s&$0.16$s&$0.18$s&$0.20$s\\ 
    \hline
    $0.5$&3.11$\times10^{-6}$&1.30$\times10^{-5}$&1.02$\times10^{-3}$&5.61$\times10^{-3}$&1.31$\times10^{-2}$&1.37$\times10^{-2}$&1.29$\times10^{-2}$&9.34$\times10^{-3}$
    \\
    \hline
    $1$&3.11$\times10^{-6}$&2.73$\times10^{-6}$&2.52$\times10^{-6}$&1.51$\times10^{-6}$&2.76$\times10^{-6}$&3.20$\times10^{-6}$&7.67$\times10^{-6}$&8.13$\times10^{-6}$
    \\
    \hline
    $6$&3.11$\times10^{-6}$&5.52$\times10^{-4}$&7.95$\times10^{-2}$&1.73$\times10^{-1}$&2.07$\times10^{-1}$&2.61$\times10^{-1}$&5.37$\times10^{-1}$&5.89$\times10^{-1}$
    \\
    \hline
    $11$&3.11$\times10^{-6}$&6.96$\times10^{-4}$&1.68$\times10^{-1}$&2.87$\times10^{-1}$&4.04$\times10^{-1}$&3.71$\times10^{-1}$&1.03&1.05
    \\
    \hline
    $16$&3.11$\times10^{-6}$&7.27$\times10^{-4}$&1.94$\times10^{-1}$&2.90$\times10^{-1}$&4.21$\times10^{-1}$&3.64$\times10^{-1}$&1.09&9.91$\times10^{-1}$
    \\
    \hline
    $19$&3.11$\times10^{-6}$&3.11$\times10^{-6}$&2.01$\times10^{-1}$&2.90$\times10^{-1}$&4.21$\times10^{-1}$&3.76$\times10^{-1}$&1.09&9.59$\times10^{-1}$
    \\
    \hline
    \end{tabular}
    }
\end{table}

\subsection{Elastic wave propagation in 3-D homogenous media}
\label{sec.homo_3d}

To study the wave equation in 3-D homogenous media, we set the initial velocity and stress tensors at equilibrium state and simulated the wave propagation activated by source functions with a Ricker-type wavelet history,
\begin{equation}
f_1(\bx, t) = f_2(\bx, t) = f_3(\bx, t) = A(\bx) f_r(t),
\end{equation}
where the amplitude $A(\bx)$ was a Gaussian profile
\begin{equation}
A(\bx) = \me^{-(x-x_0)^2 - (y - y_0)^2 - (z - z_0)^2} 
\end{equation}
and the Ricker wavelet was given by
\begin{equation}
f_r(t) = (1 - 2(\pi f_P (t- d_r)^2)) \me^{-(\pi f_P(t- d_r))^2},
\end{equation}
where $f_P$ was the peak frequency and $d_r$ was the temporal delay. Here we chose $f_P = 100$Hz, $d_r = 0$ and the centre position $x_0 = 0$, $y_0 = 0$, $z_0 = 10$. The model parameters in a homogenous media were set as follows: The group velocities for $P$-wave and $S$-wave were  $c_P = 2.614$km/s and $c_S = 0.802$km/s, respectively, and the mass density was $\rho = 2.2 \textup{kg}/\textup{m}^3$.

The computational domain was $[-40, 40]^3$ ($80\textup{km}\times 80 \textup{km}\times 80 \textup{km}$).
Four groups of simulations have been performed under $N = N_x \times N_y \times N_z = 129^3, 257^3, 385^3, 513^3$. 
The Strang splitting was adopted with time step $\Delta t = 0.005$s and the final instant was $T = 10$s (2000 steps in total). 
 The reference solution was produced by FSM with a fine grid mesh $N = 512^3$ for testing the convergence of LOSS. 

As a pretreatment in LOSS, the domain was evenly decomposed into $4\times 4\times 4$ patches, thereby achieving a perfect balance in overload.  Here PMBCs were assembled by  truncating the neighborhoods at $n_{nb} = 20$ \cite{XiongZhangShao2022_arXiv}. A visualization of the domain decomposition in 3-D space is presented in Figure \ref{domain_decomposition_3d}, when the whole computational domain is decomposed into $2\times 2\times 2$ cells. For the calculations of first-order spatial derivatives, it only requires the local spline expansion in one direction, where PMBCs are assembled at the shared faces.  For instance, when one needs the first-order derivatives in $x$-direction, it only requires the communications in $1 \leftrightarrow 2$,  $3 \leftrightarrow 4$, $5 \leftrightarrow 6$ and $7 \leftrightarrow 8$, thereby greatly reducing the communication cost.

\begin{figure}[!h]
    \centering
    \subfigure[Convergence with respect to $N_z$. (left: maximal error $\varepsilon_{\infty}(10)$, right: $l^2$-error $\varepsilon_{2}(10)$)\label{3d_homo_convergnce}]{
    {\includegraphics[width=0.49\textwidth,height=0.24\textwidth]{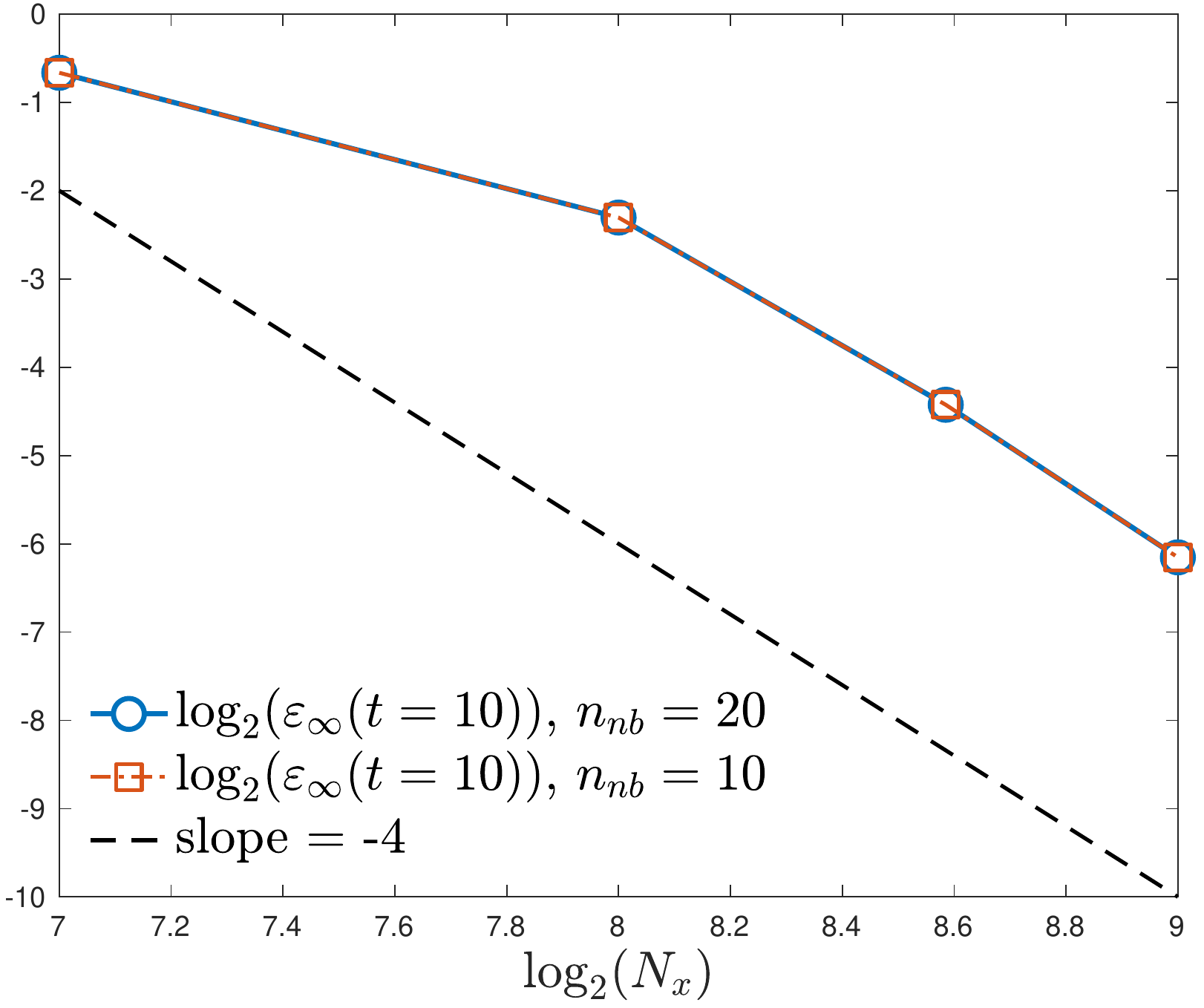}}
     {\includegraphics[width=0.49\textwidth,height=0.24\textwidth]{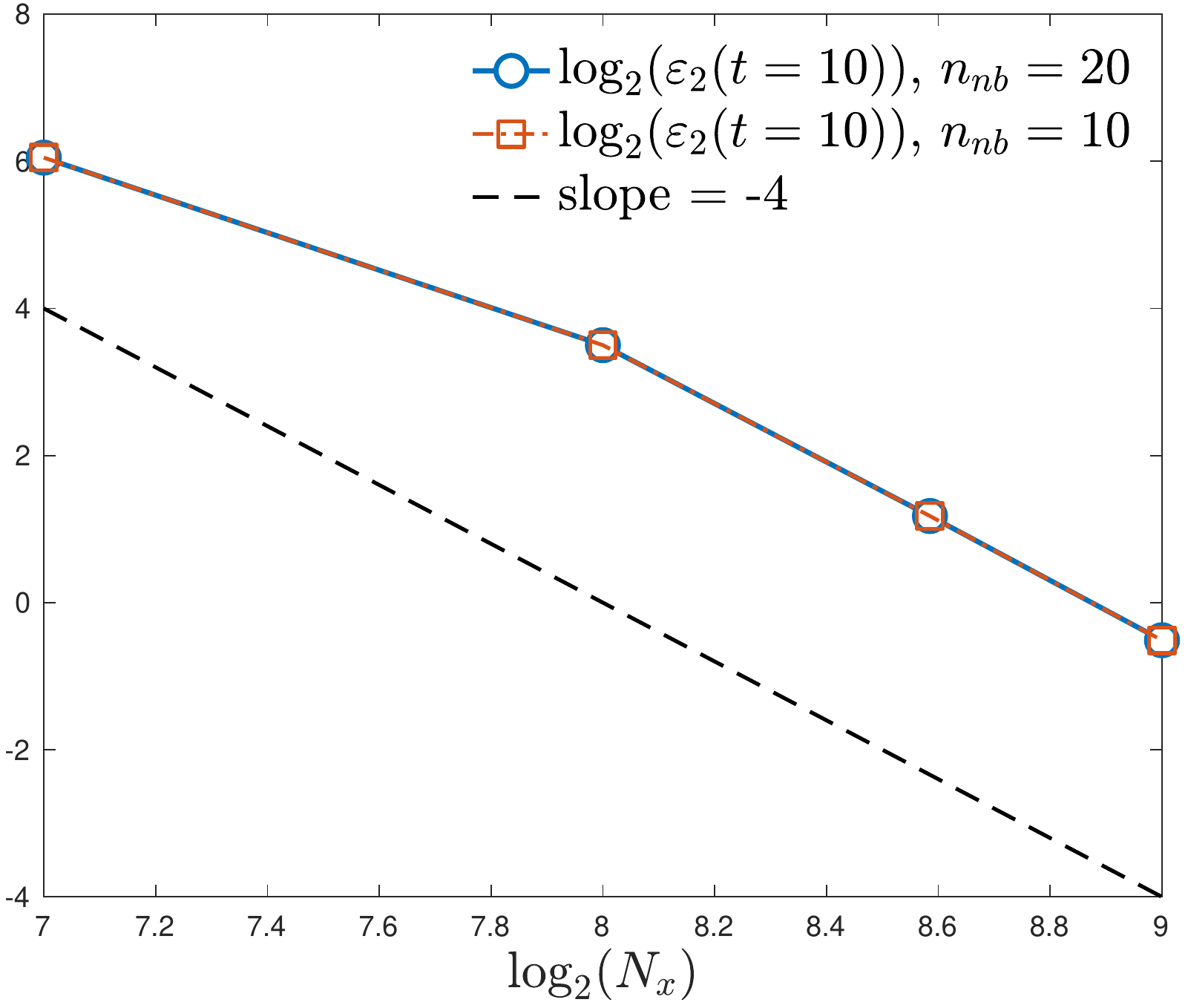}}}
     \\
     \centering
    \subfigure[$v_3(0, 0, z)$  in the homogenuous media (left) and the relative errors (right) at $t=10$s. \label{v3_homo_t10}] 
    {\includegraphics[width=0.49\textwidth,height=0.24\textwidth]{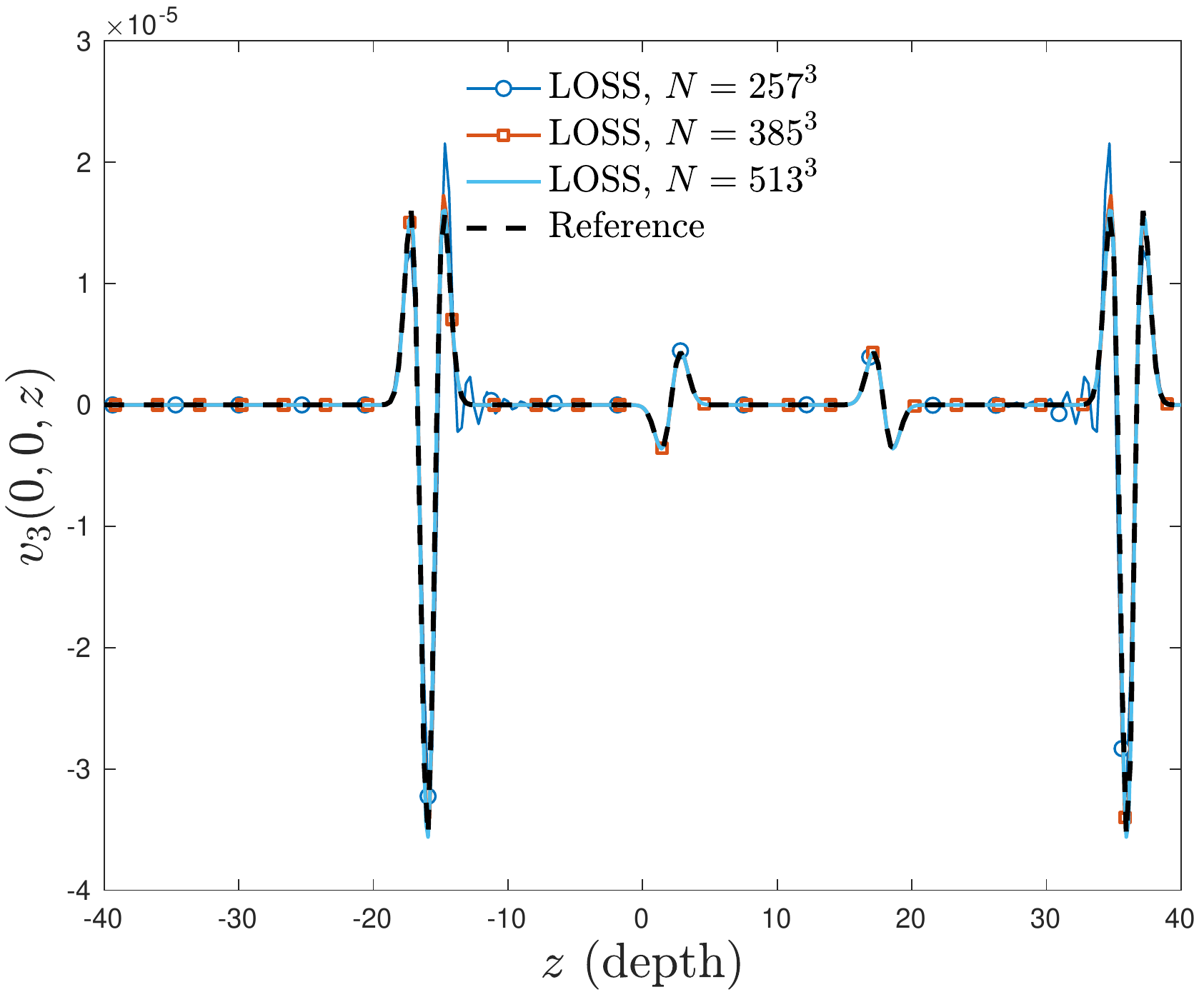}
     {\includegraphics[width=0.49\textwidth,height=0.24\textwidth]{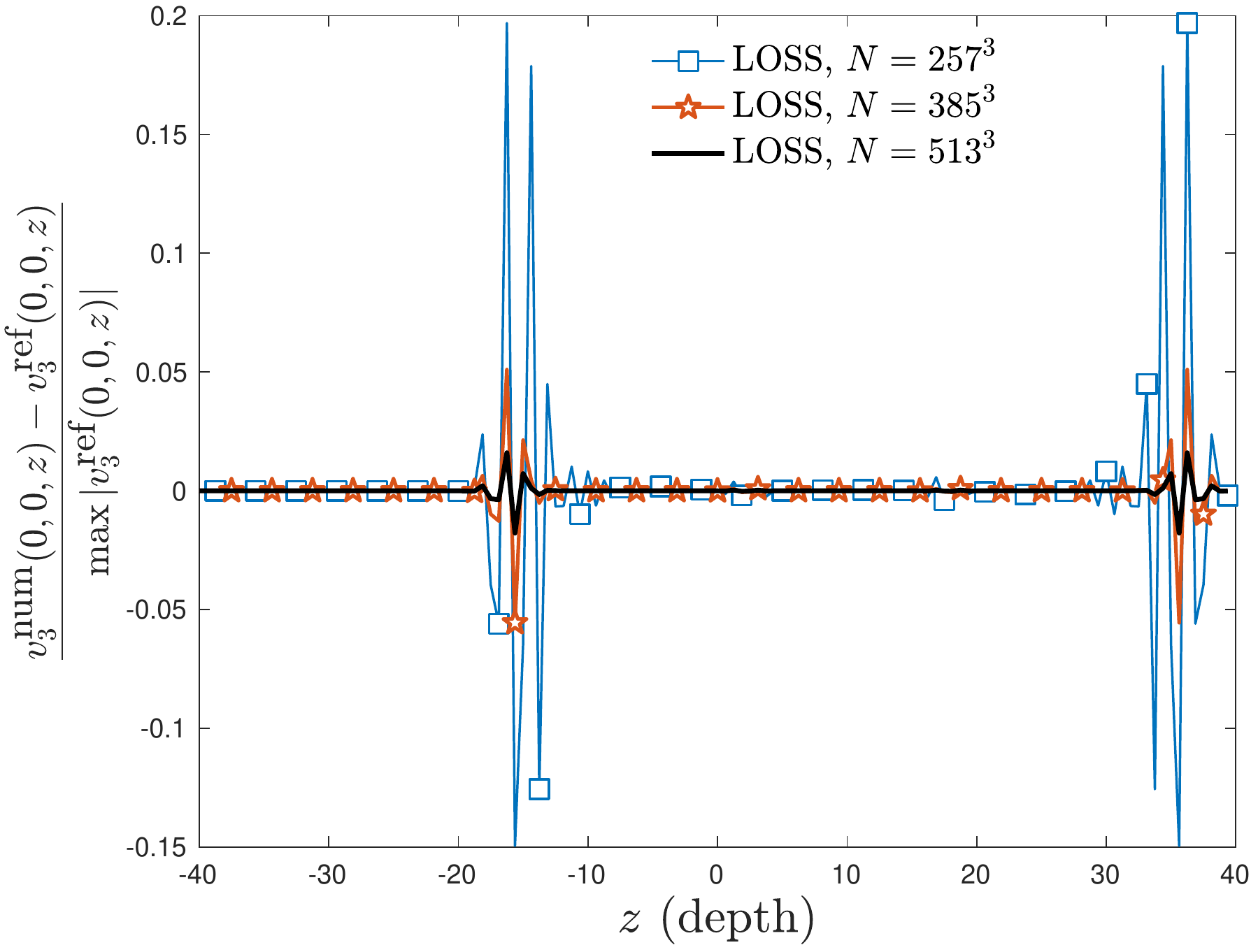}}}
    \caption{\small Wave propagation in 3-D homogenous media: A comparison of the vibrational wavefield $v_3(0, 0, z)$ at $t =10$s. The convergence of LOSS is verified for smooth data and coefficients.}
\end{figure} 

\begin{figure}[!h]
    \centering
    \subfigure[$v_3(x, y, z)$ at $t = 4$s. (left: FSM, $N=512^3$, right: LOSS, $N=513^3$)]
    {\includegraphics[width=0.49\textwidth,height=0.24\textwidth]{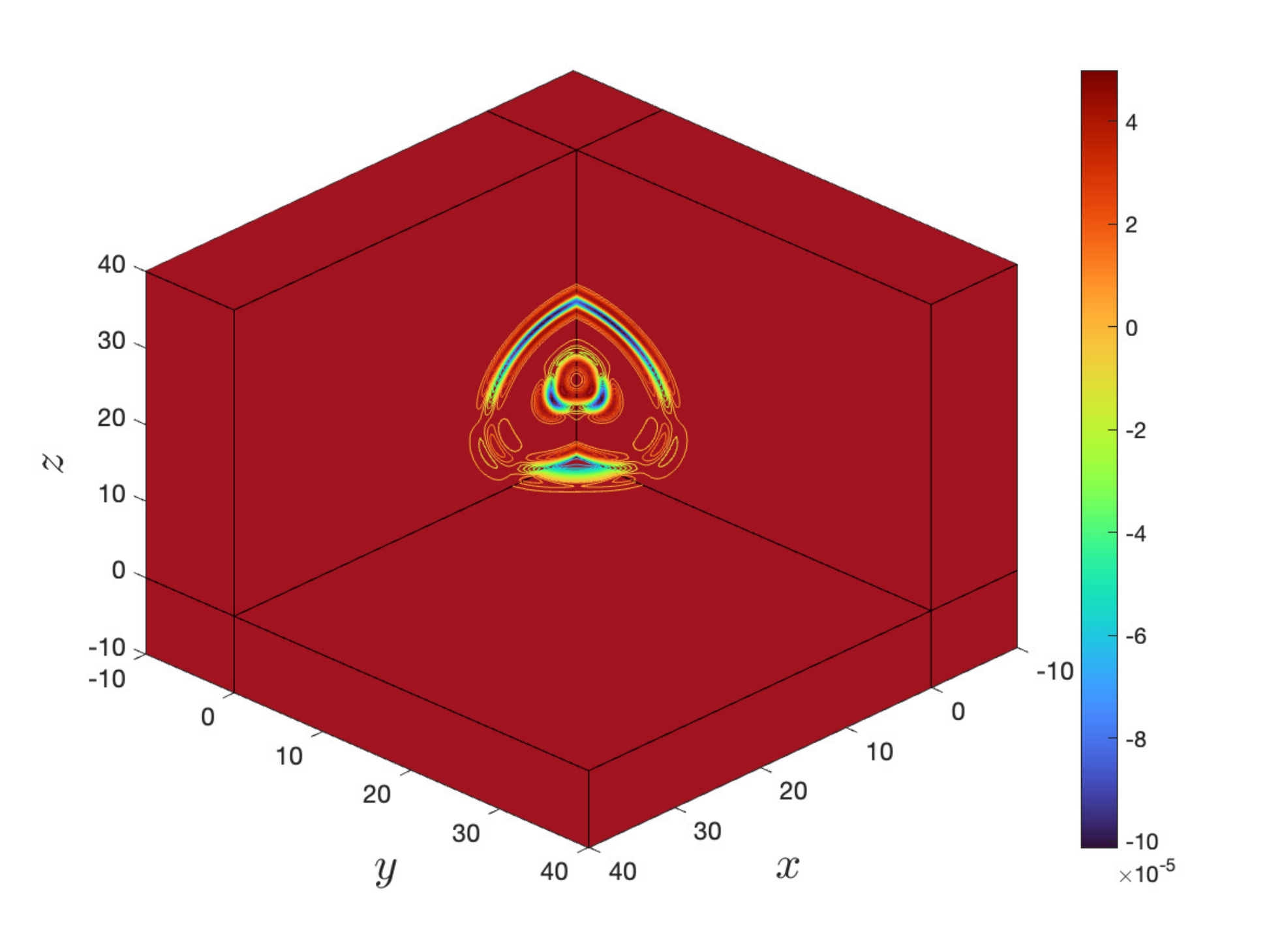}
     \includegraphics[width=0.49\textwidth,height=0.24\textwidth]{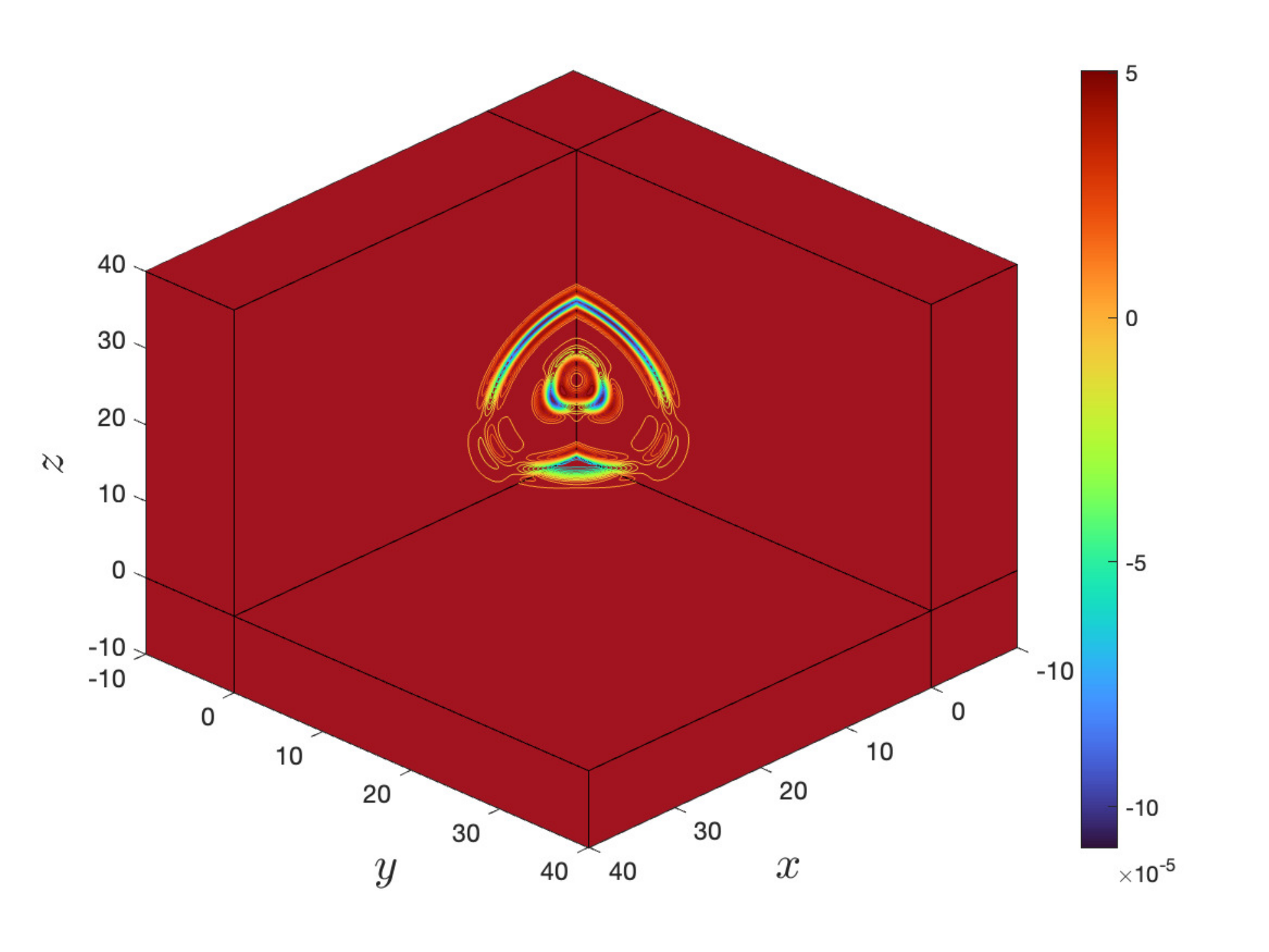}}
     \\
     \centering
    \subfigure[$v_3(x, y, z)$ at $t = 6$s. (left: FSM, $N=512^3$, right: LOSS, $N=513^3$)]
    {\includegraphics[width=0.49\textwidth,height=0.24\textwidth]{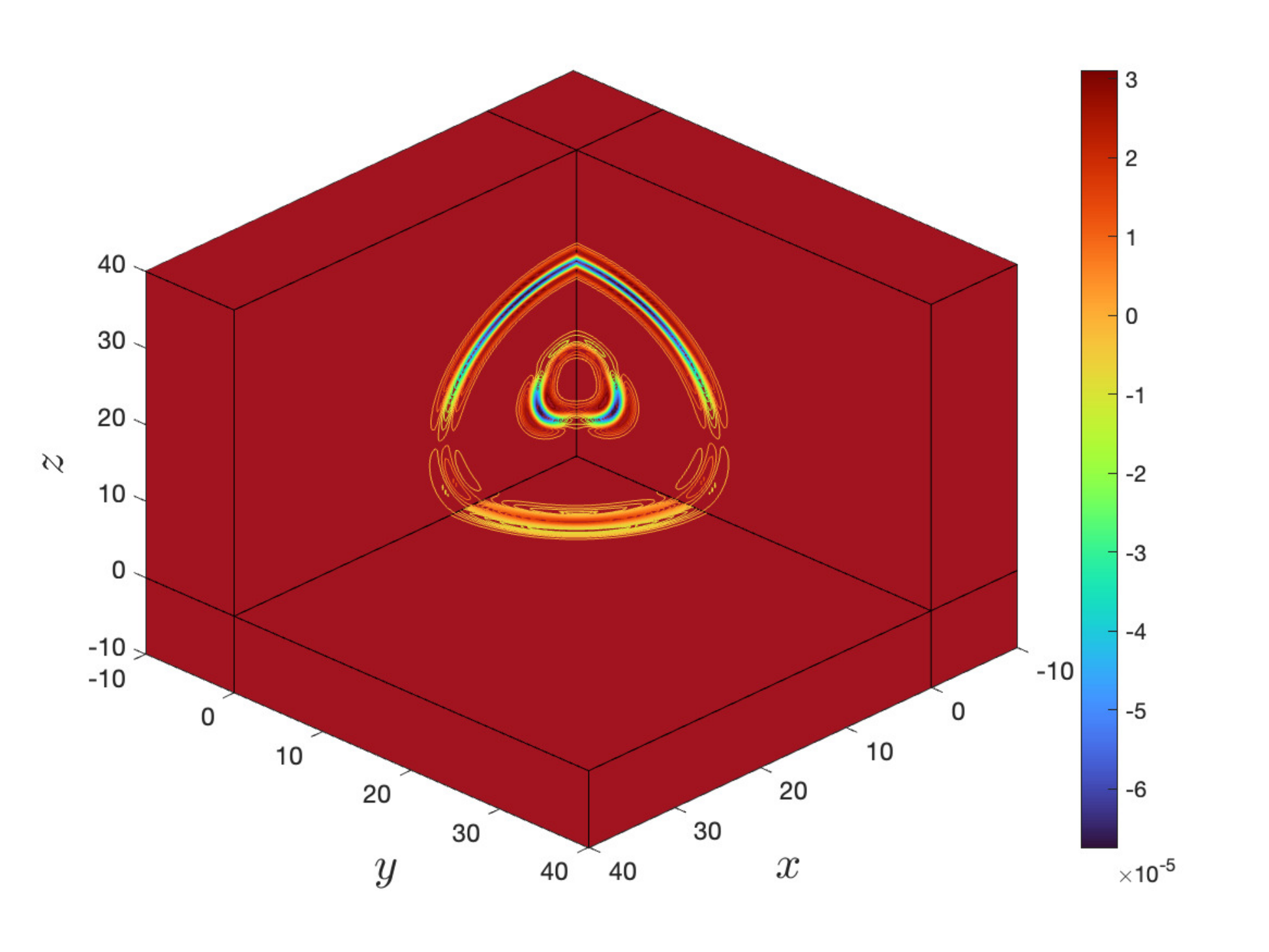}
     \includegraphics[width=0.49\textwidth,height=0.24\textwidth]{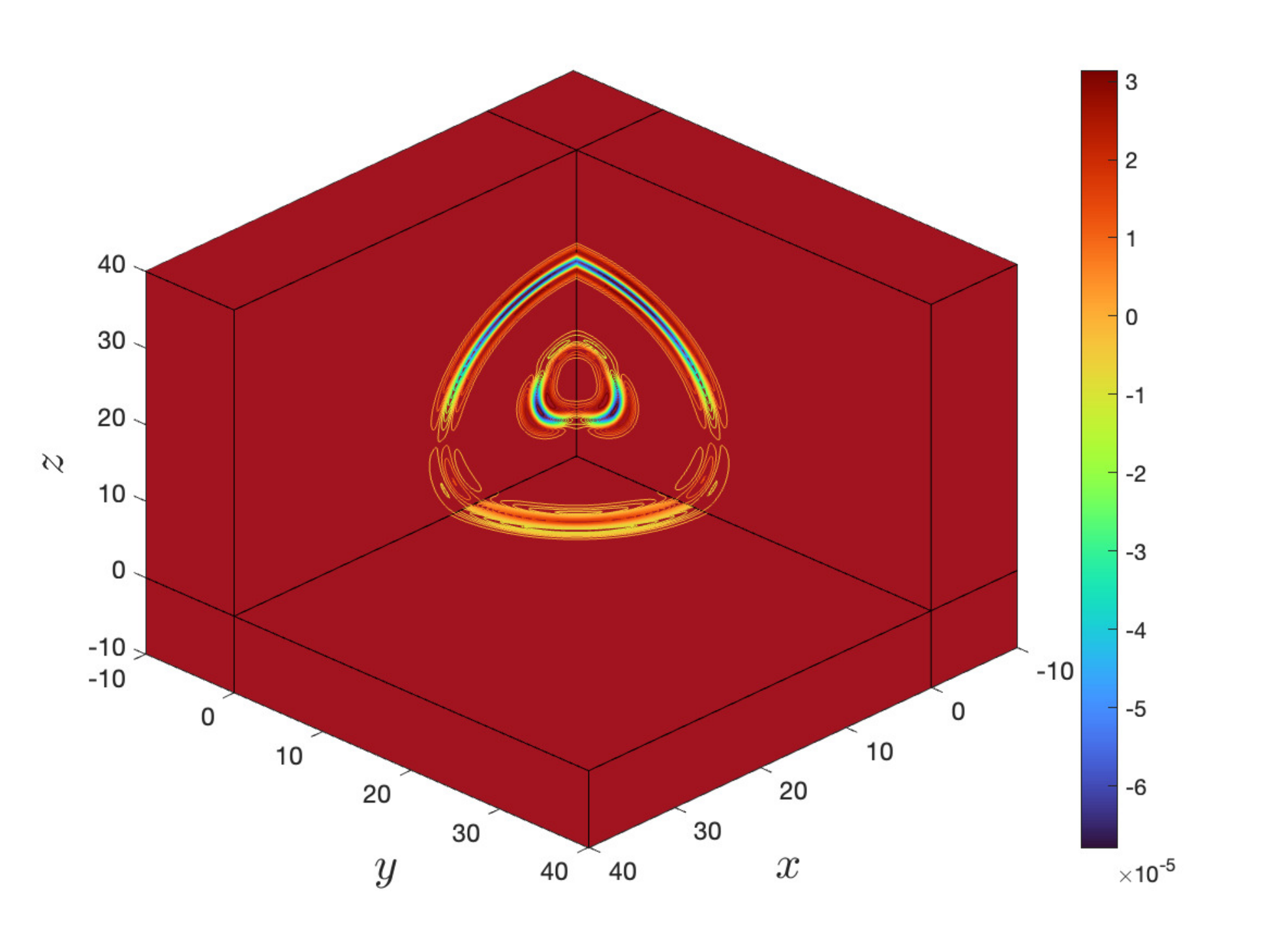}}
     \\
    \centering
    \subfigure[$v_3(x, y, z)$ at $t = 8$s. (left: FSM, $N=512^3$, right: LOSS, $N=513^3$)]
    {\includegraphics[width=0.49\textwidth,height=0.24\textwidth]{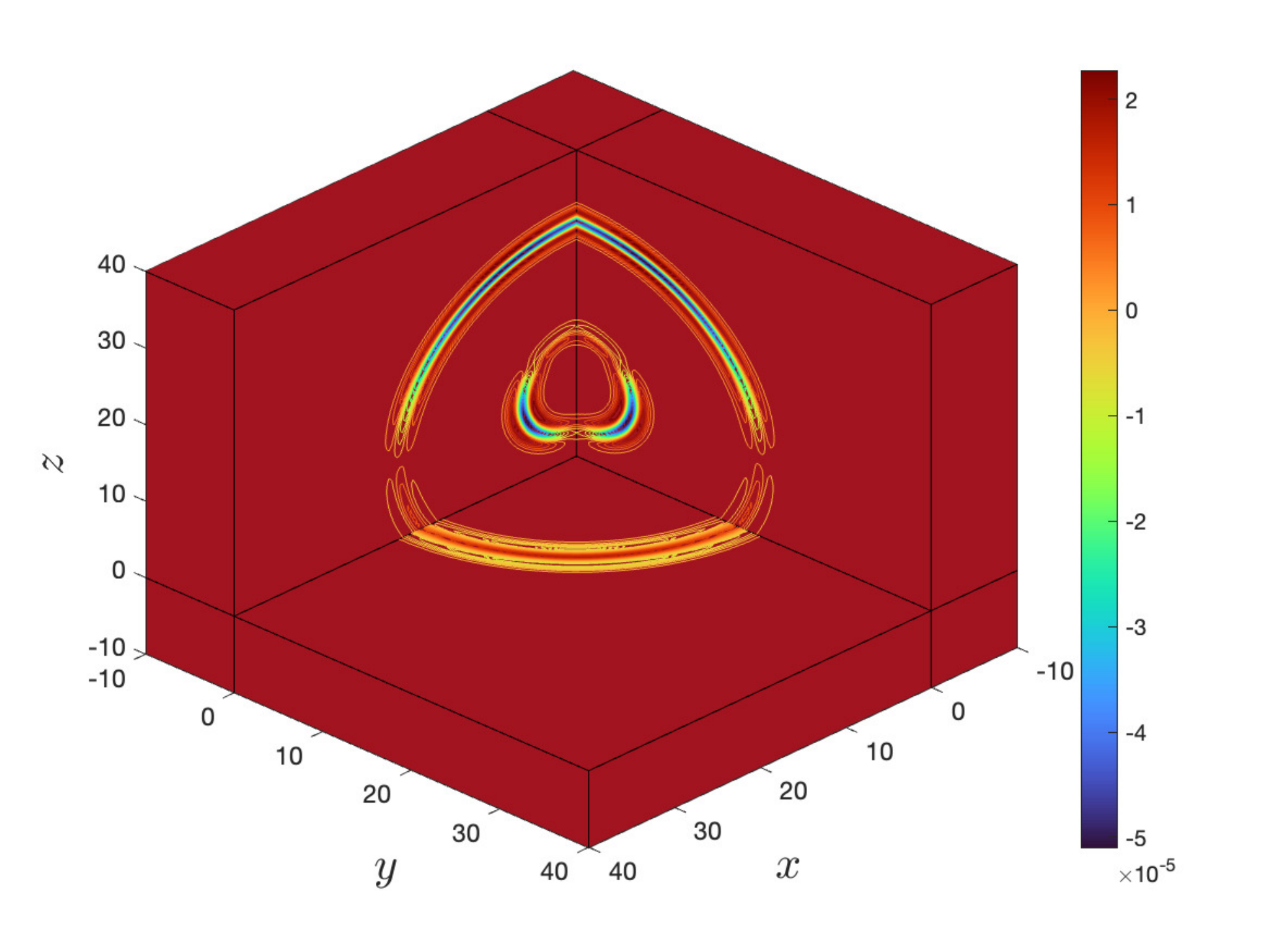}
     \includegraphics[width=0.49\textwidth,height=0.24\textwidth]{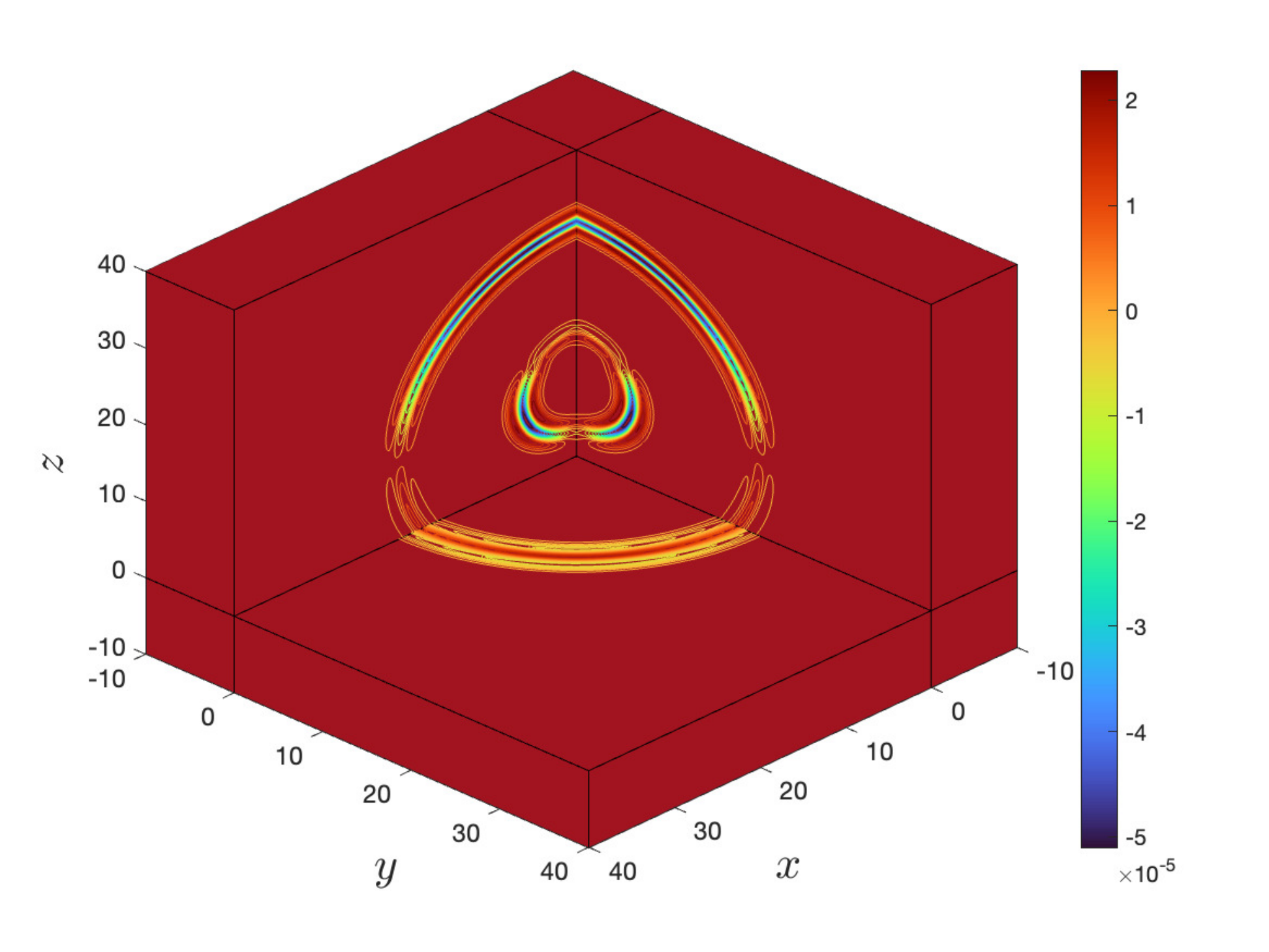}}
     \\
    \centering
    \subfigure[$v_3(x, y, z)$ at $t = 10$s. (left: FSM, $N=512^3$, right: LOSS, $N=513^3$)]
    {\includegraphics[width=0.49\textwidth,height=0.24\textwidth]{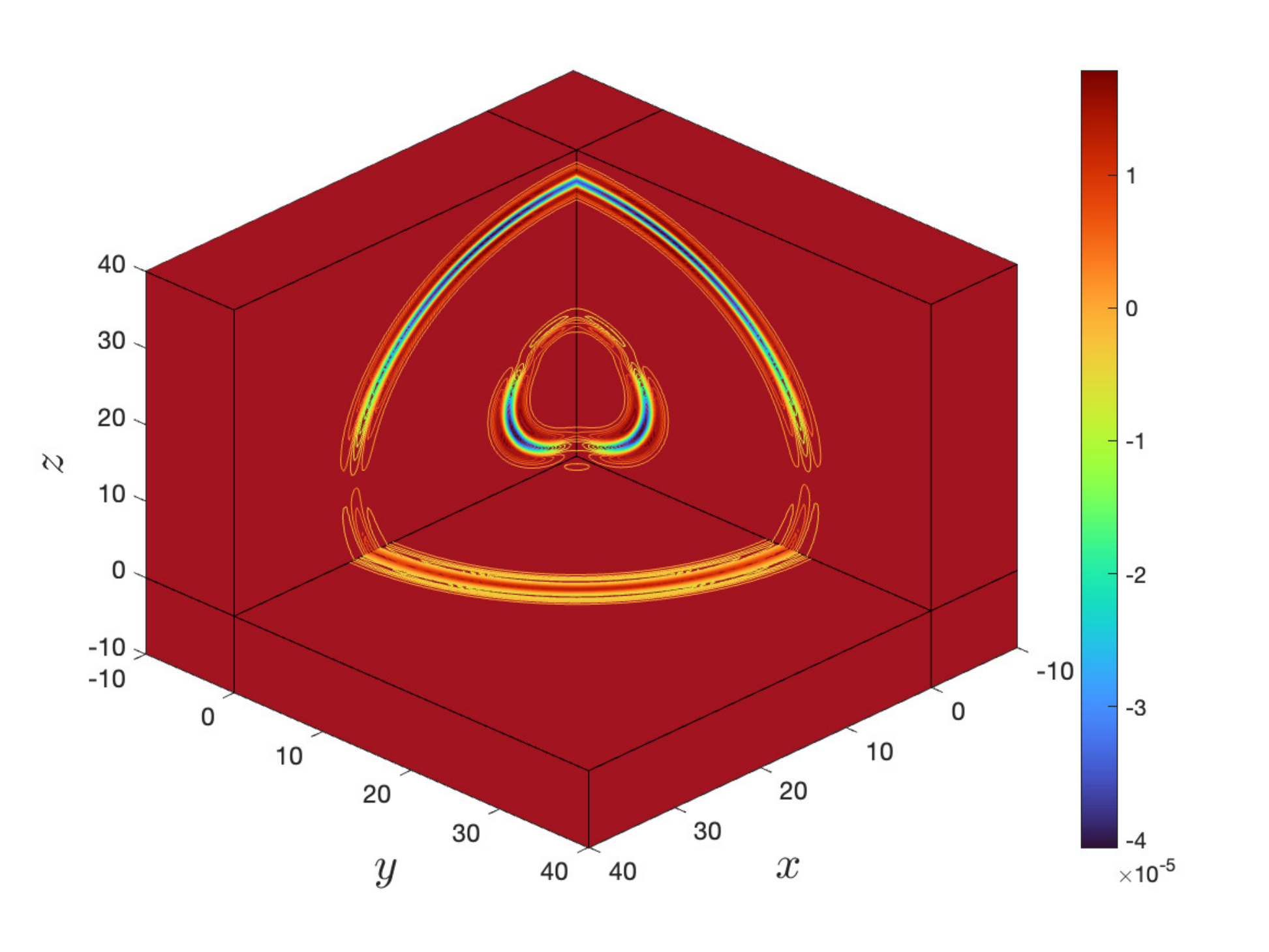}
     \includegraphics[width=0.49\textwidth,height=0.24\textwidth]{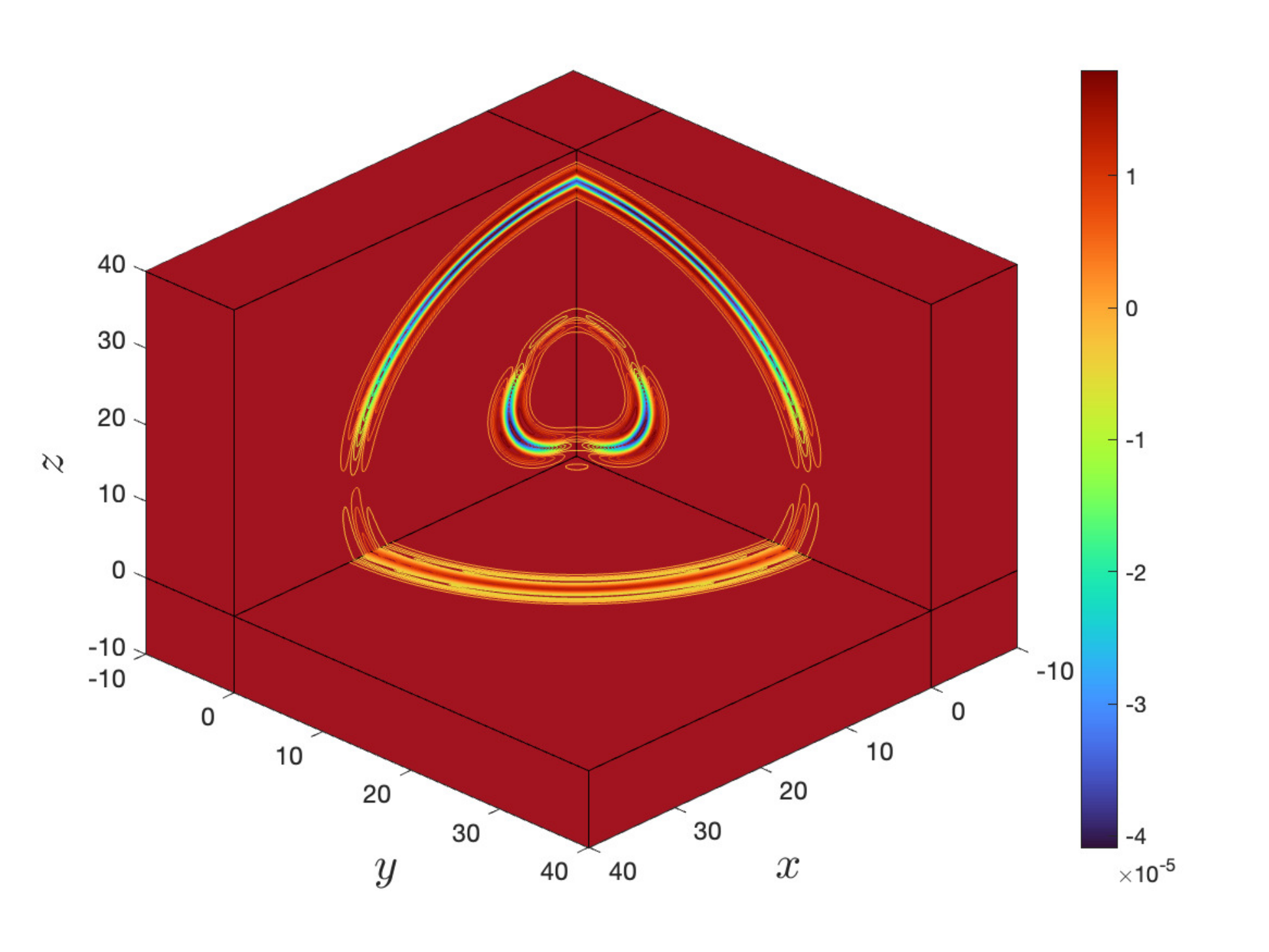}}
    \caption{\small Wave propagation in 3-D homogenous media: A comparison of snapshots of wavefield $v_3(x, y, z)$.\label{homo_PV}}
\end{figure} 

The relative maximal error $\varepsilon_{\infty}(t)$ and $l^2$-error $\varepsilon_2(t)$ were again adopted to measure the numerical accuracy and the spatial convergence {\cf is} plotted in Figure \ref{3d_homo_convergnce}. Clearly, the convergence rate coincides with the theoretical fourth order. We also investigate the influence of the parameter $n_{nb}$ in assembling PMBCs.  As shown in Figure \ref{3d_homo_convergnce}, the errors induced by truncation in PMBC seem negligible even when $n_{nb}=10$, which accords with our observations made in \cite{XiongZhangShao2022_arXiv}.

The snapshots of the vibrational wave-field $v_3$ of the $P$- and $S$- wave are given in Figure \ref{homo_PV}, where three sectional drawings are provided to visualize the wave propagation in the homogenous media. To further demonstrate the errors of LOSS, we plot the synthetic data $v_3(0, 0, z)$ at $t = 10$s in Figure \ref{v3_homo_t10}. Small oscillations are observed around the reference solution when  too coarse grid mesh is used ($N = 257^3$). Nonetheless, they could be suppressed to a large extent when the finer grid mesh was adopted.

\subsection{Elastic wave propagation in a 3-D double-layer media}

\begin{figure}[!h]
\centering
\subfigure[Domain decomposition in 3-D space. \label{domain_decomposition_3d}]{\includegraphics[width=0.49\textwidth,height=0.24\textwidth]{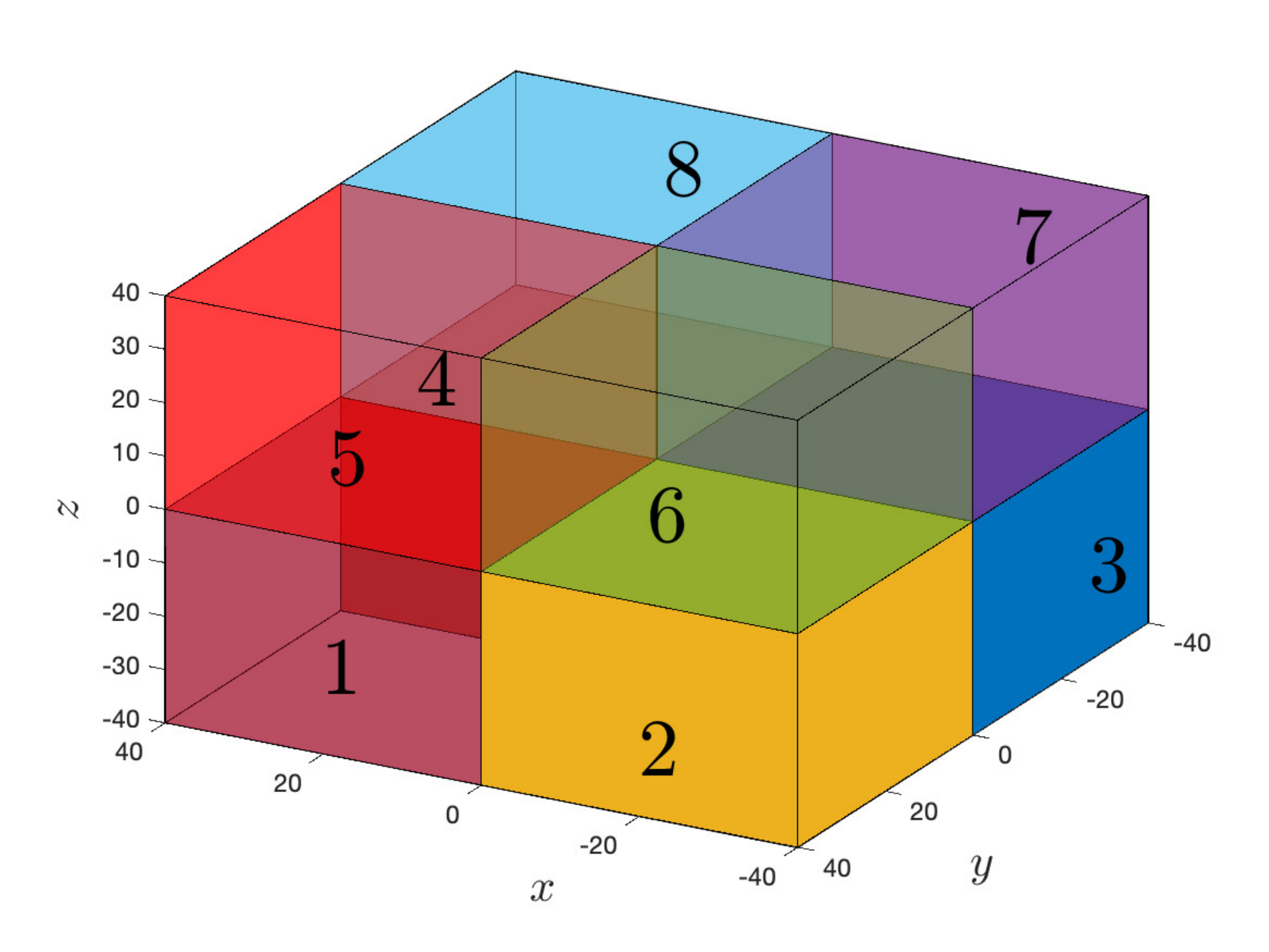}}
\subfigure[A double-layer heterogeneous structure. \label{double_layer}]
{\includegraphics[width=0.49\textwidth,height=0.24\textwidth]{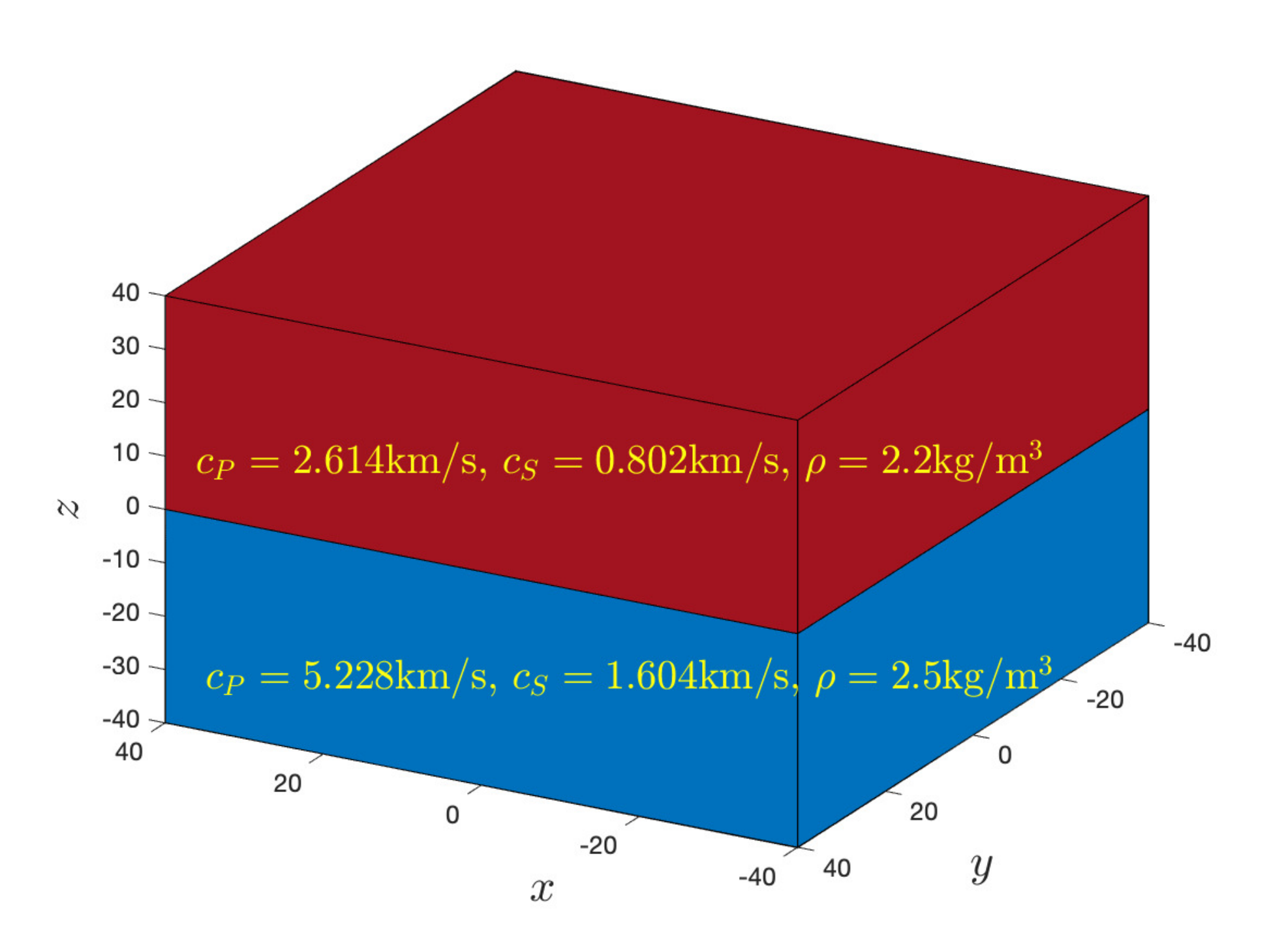}}
\caption{An illustration of domain decomposition of 3-D space and a double-layer heterogeneous structure.}
\end{figure}

As a more challenging example, we also studied the wave propagation in a double-layer structure, with the same computational domain $[-40, 40]^3$ and the same source impulse as adopted in Section \ref{sec.homo_3d}. The model parameters in a heterogenous media, iven in Figure \ref{double_layer}, were set to be: The group velocities for the upper layer $40^2 \times [0, 40]$ were $c_P = 2.614$km/s and $c_S = 0.802$km/s and the mass density was $\rho = 2.2 \textup{kg}/\textup{m}^3$, while for the deeper layer $40^2 \times [-40, 0]$, $c_P = 5.228$km/s, $c_S = 1.604$km/s and $\rho = 2.5\textup{kg}/\textup{m}^3$.

Seven groups of simulations have been performed under $N = N_x \times N_y \times N_z = 129^3, 257^3, 385^3, 513^3, 641^3, 769^3$ and $513^2\times 1025$. 
The Strang splitting was adopted with time step $\Delta t = 0.005$s and the final instant was $T = 10$s (2000 steps in total). 
 The reference solution was still produced by FSM with a fine grid mesh $N = 640^3$.

\begin{figure}[!h]
    \centering
    \subfigure[Convergence with respect to $N_z$  (left: maximal error $\varepsilon_{\infty}(10)$, right: $l^2$-error $\varepsilon_{2}(10)$). \label{3d_hetero_convergnce}]
    {\includegraphics[width=0.49\textwidth,height=0.24\textwidth]{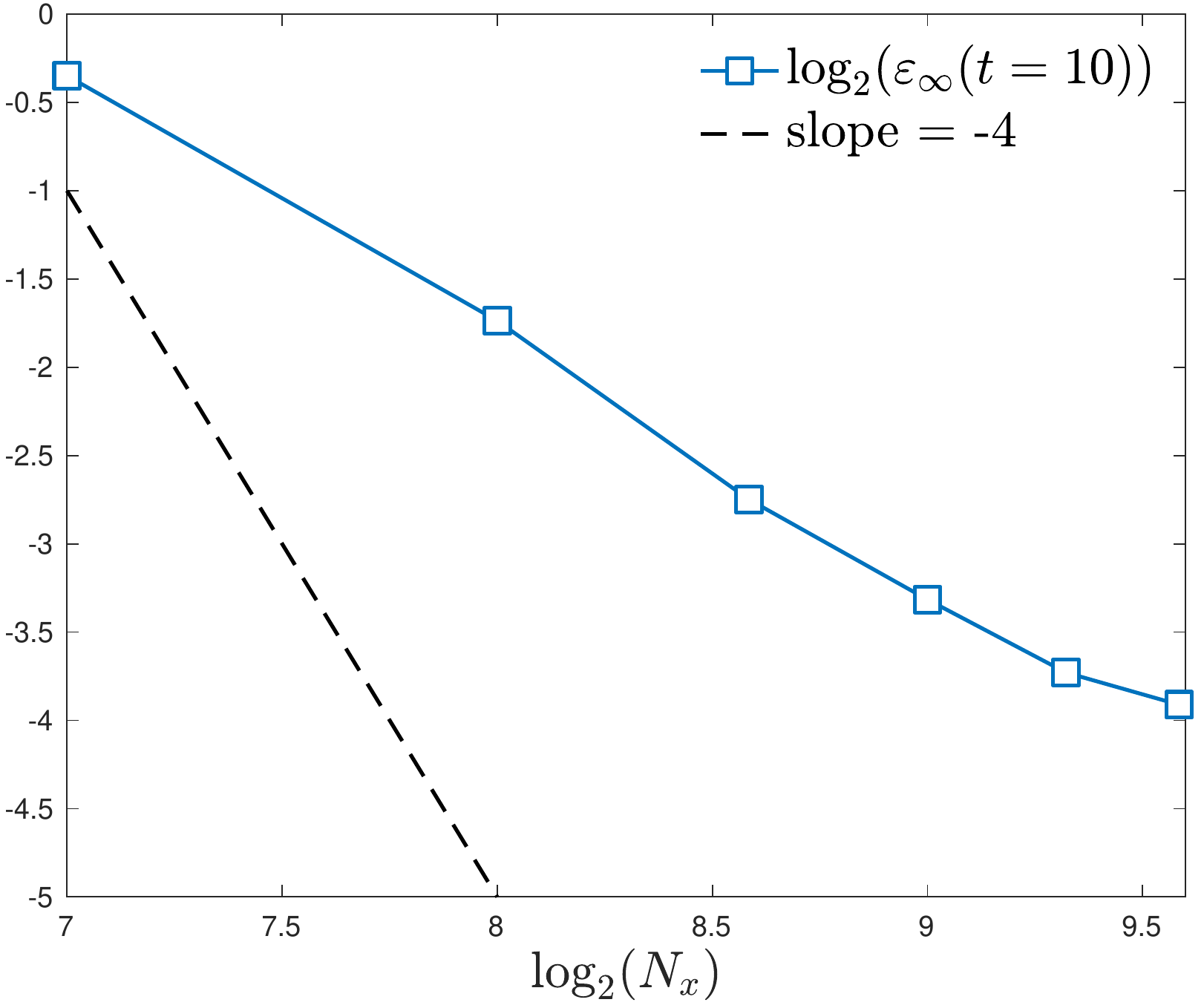}
    {\includegraphics[width=0.49\textwidth,height=0.24\textwidth]{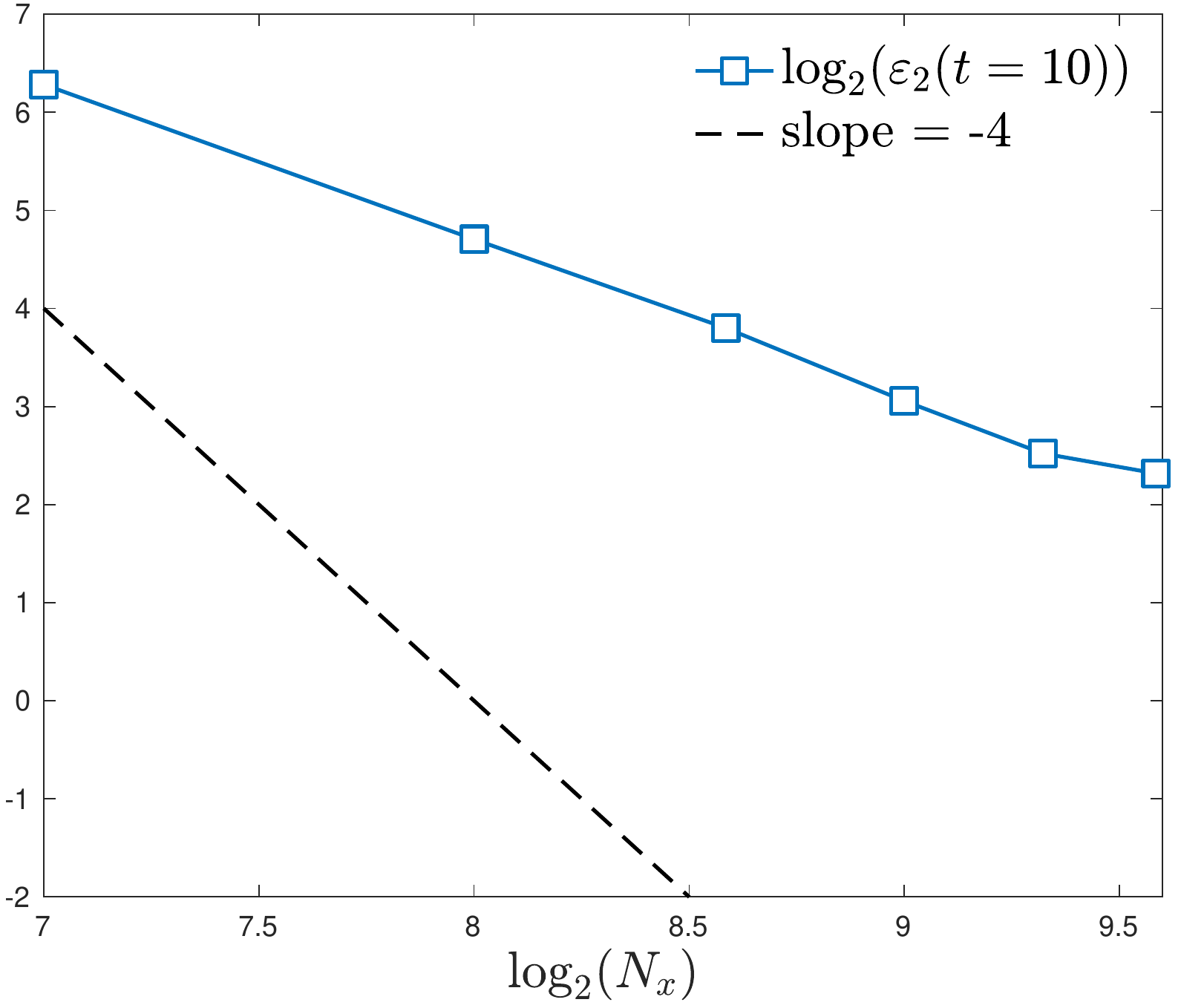}}}
     \\
     \centering
    \subfigure[$v_3(0, 0, z)$  in the heterogenuous media (left) and the relative errors (right) at $t=6$s. \label{v3_hetero_t6}]
    {\includegraphics[width=0.49\textwidth,height=0.24\textwidth]{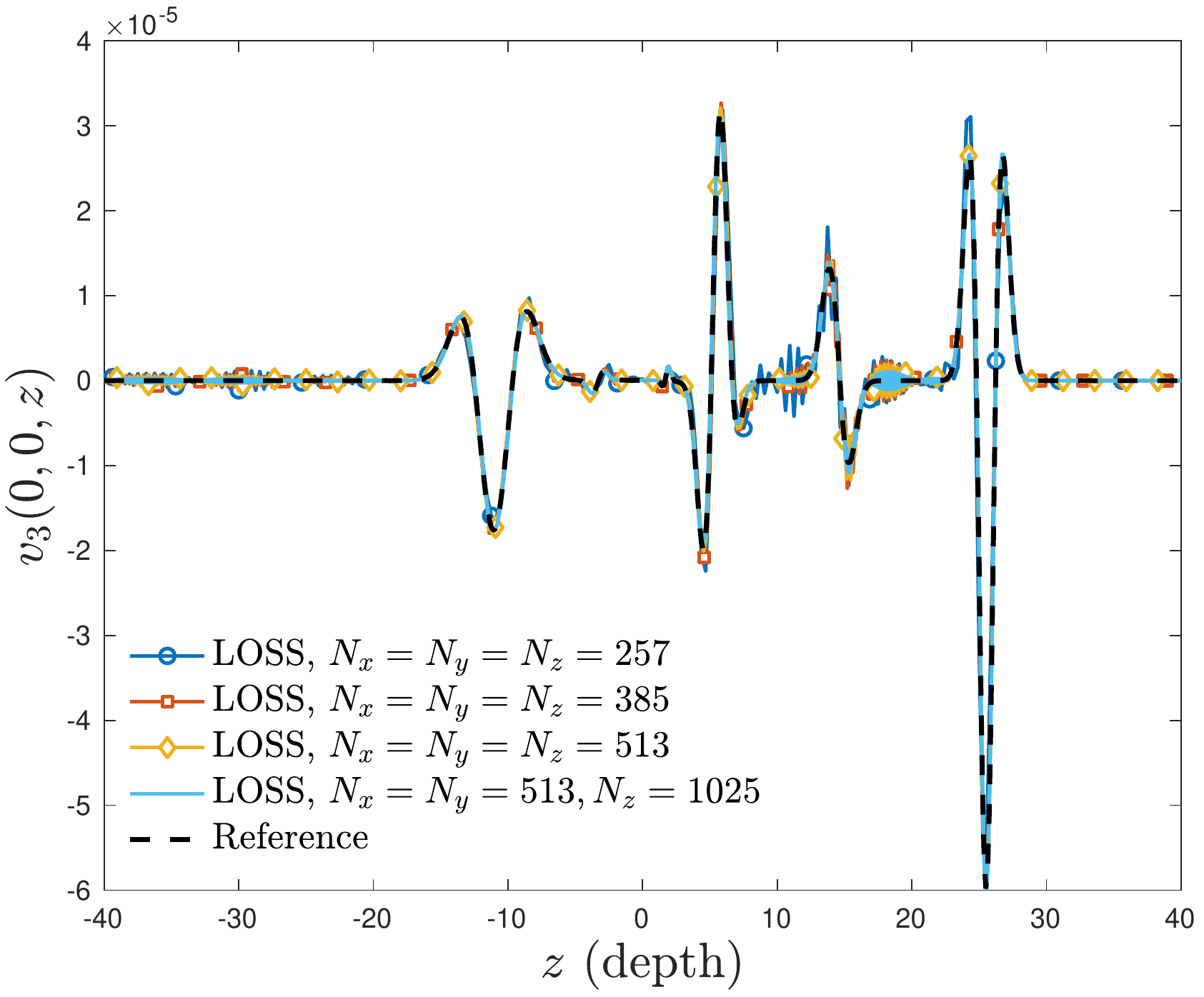}
     \includegraphics[width=0.49\textwidth,height=0.24\textwidth]{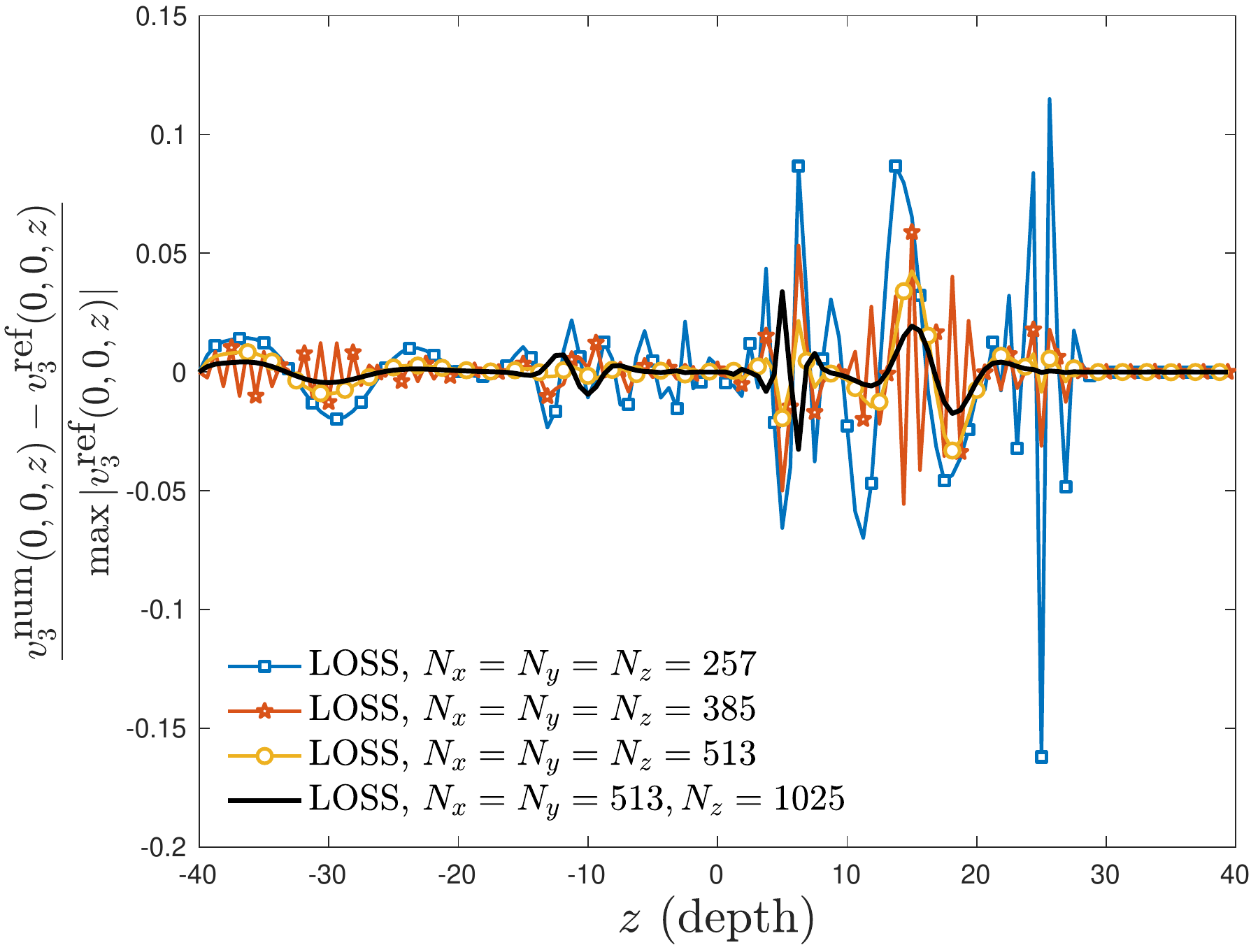}}
     \\
    \centering
    \subfigure[$v_3(0, 0, z)$  in the heterogenuous media (left) and the relative errors (right) at $t=10$s. \label{v3_hetero_t10}]
    {\includegraphics[width=0.49\textwidth,height=0.24\textwidth]{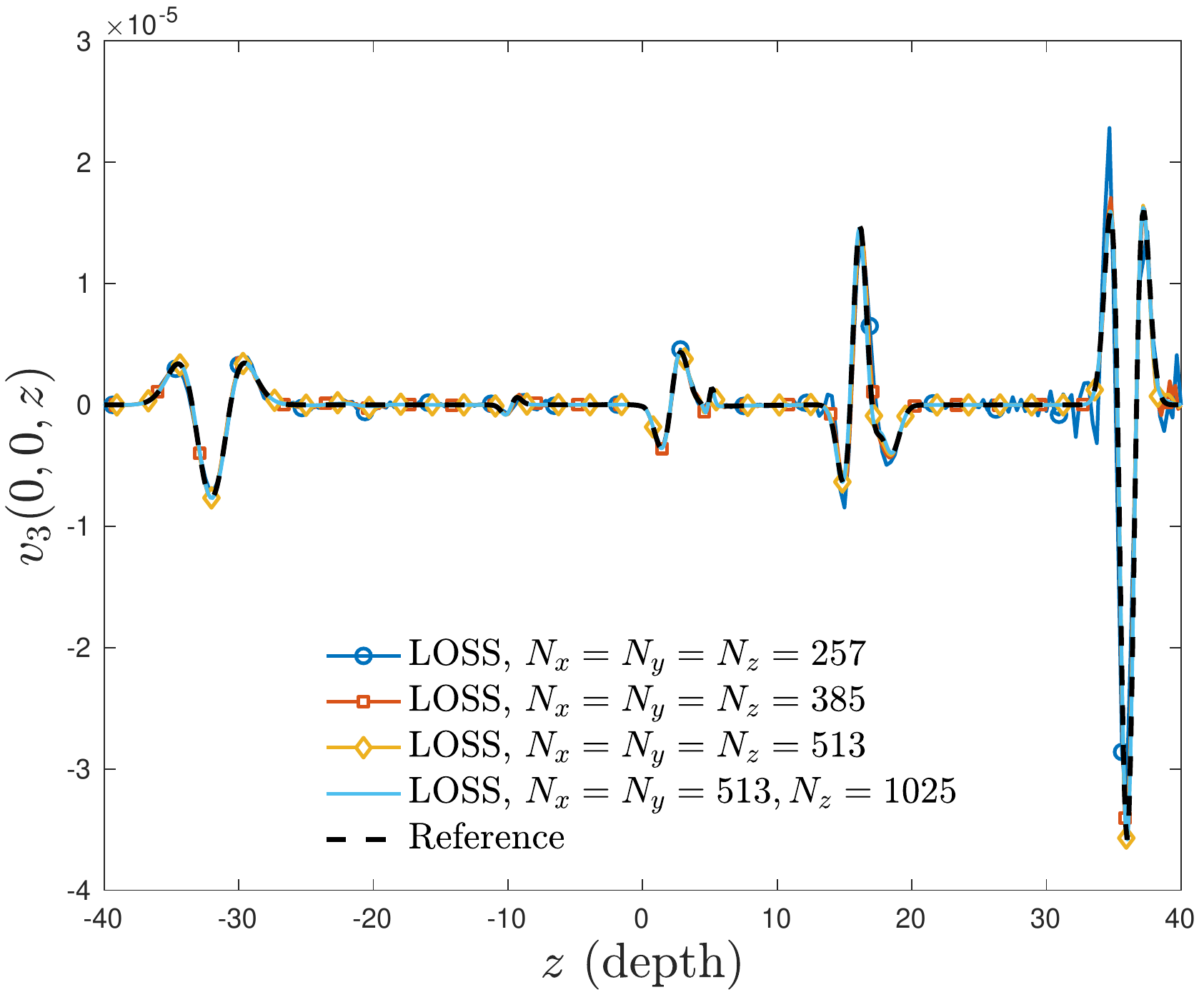}
     \includegraphics[width=0.49\textwidth,height=0.24\textwidth]{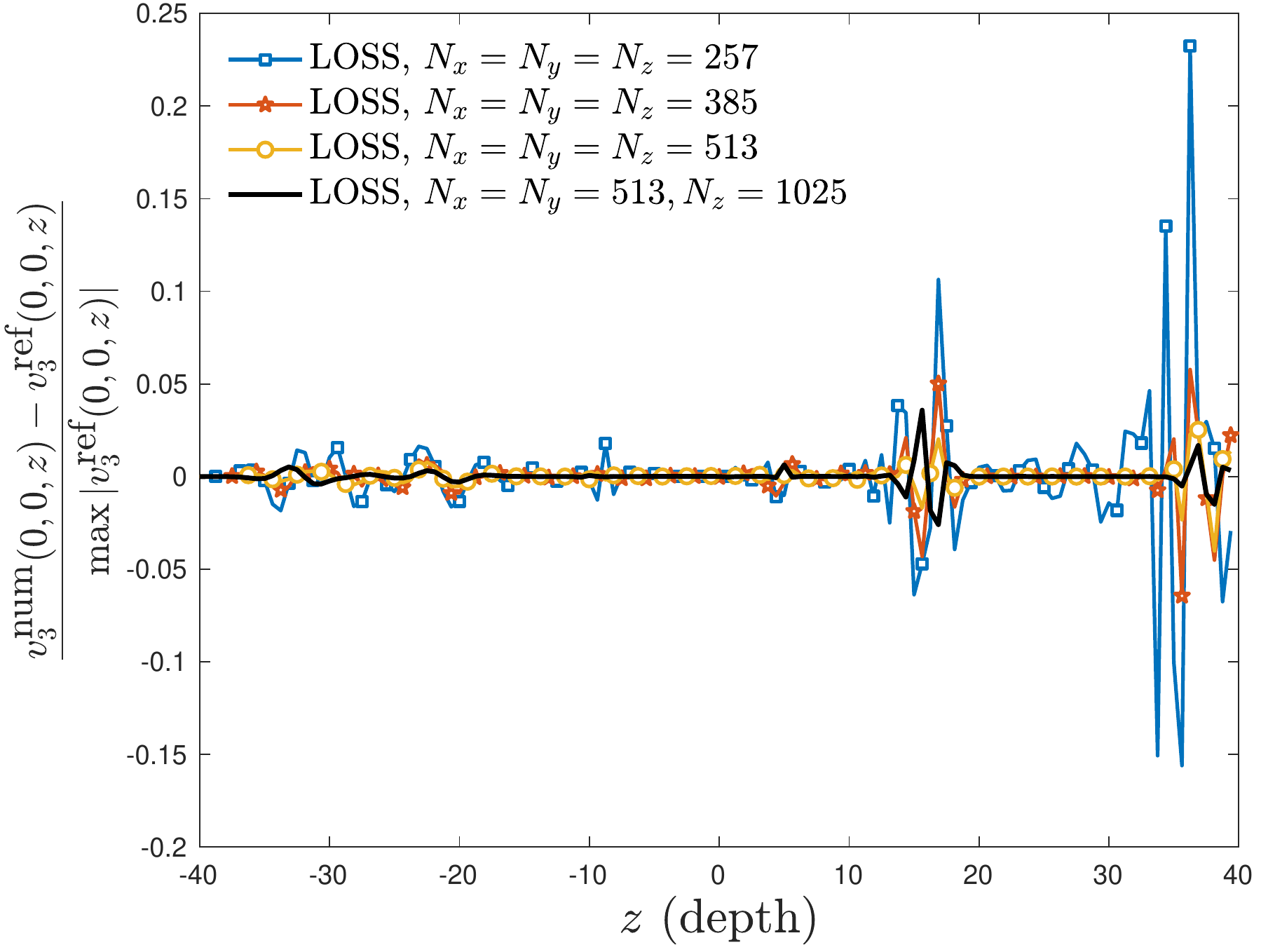}}
    \caption{\small Wave propagation in 3-D heterogeneous media: A comparison of the vibrational data  $v_3(0, 0, z)$ at $t =6$s and $t =10$s. The reduction in the convergence order is observed when model parameters have sharp variations, and small oscillations are produced by the spline construction. Nevertheless, the errors can be diminished as a finer grid mesh is adopted.}
\end{figure} 

\begin{figure}[!h]
    \centering
    \subfigure[$v_3(x, y, z)$ at $t = 4$s. (left: FSM, $N=512^3$, right: LOSS, $N=513^2\times 1025$)]
    {\includegraphics[width=0.49\textwidth,height=0.24\textwidth]{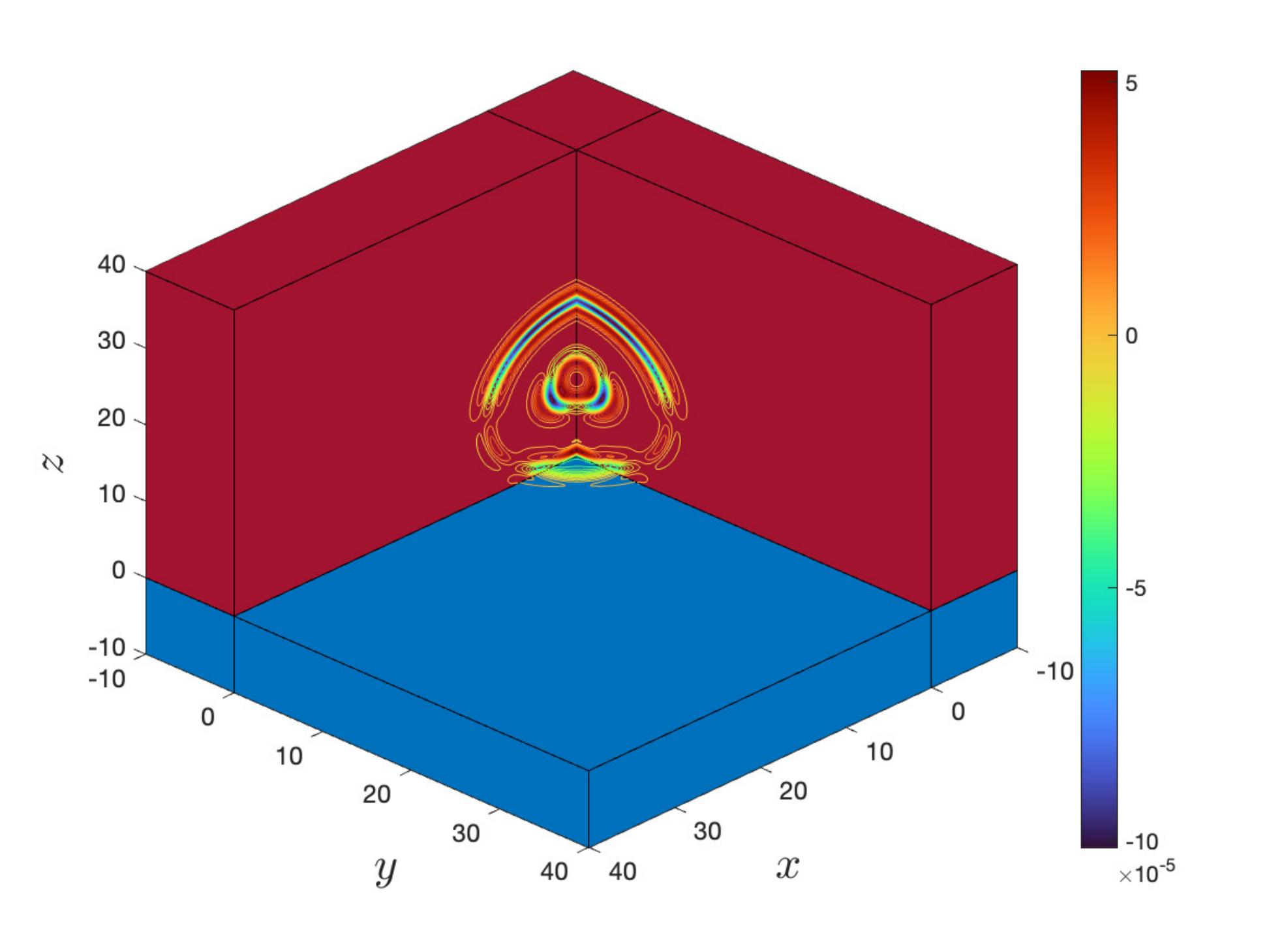}
     \includegraphics[width=0.49\textwidth,height=0.24\textwidth]{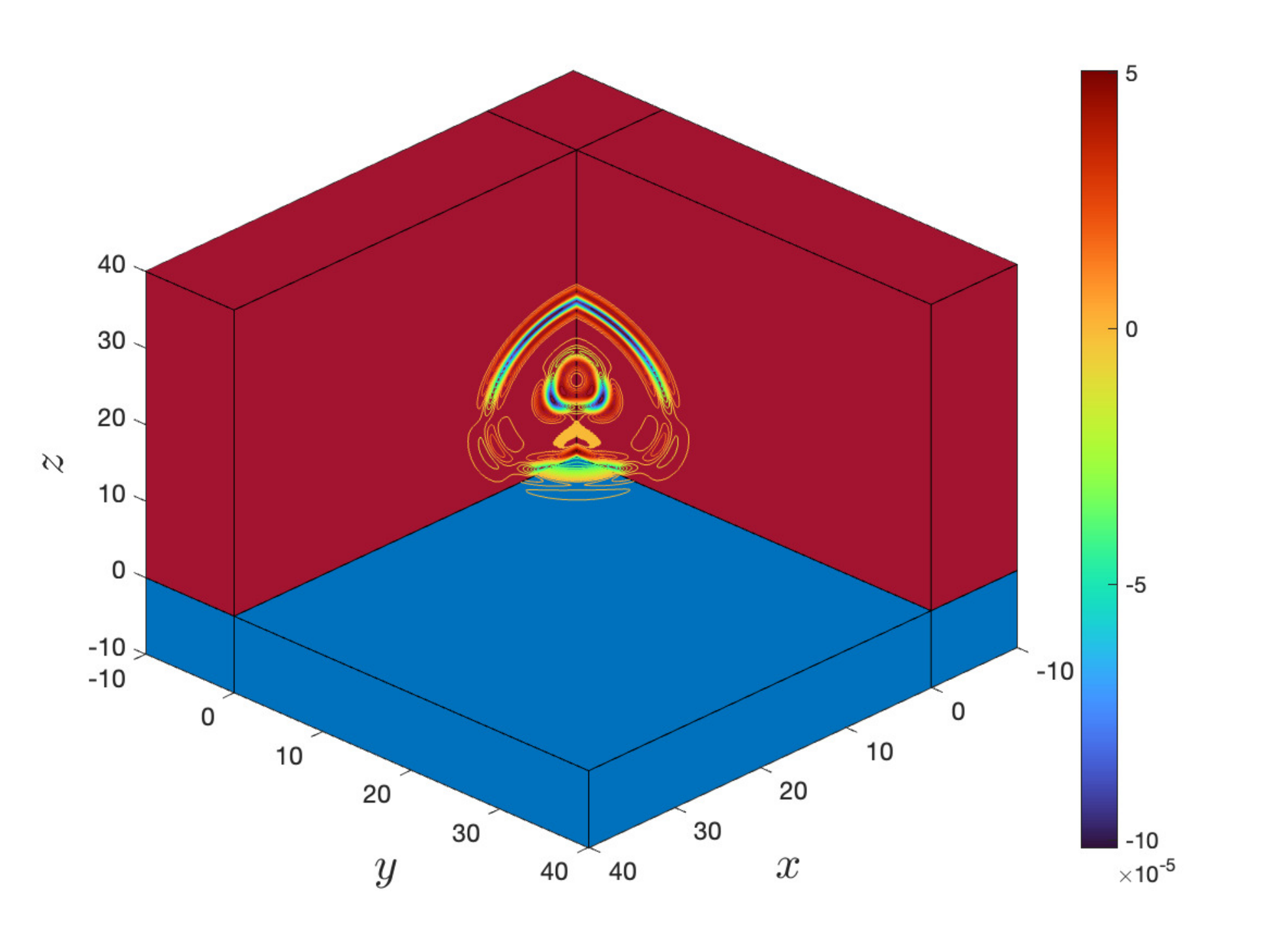}}
     \\
     \centering
    \subfigure[$v_3(x, y, z)$ at $t = 6$s. (left: FSM, $N=512^3$, right: LOSS, $N=513^2\times 1025$)]
    {\includegraphics[width=0.49\textwidth,height=0.24\textwidth]{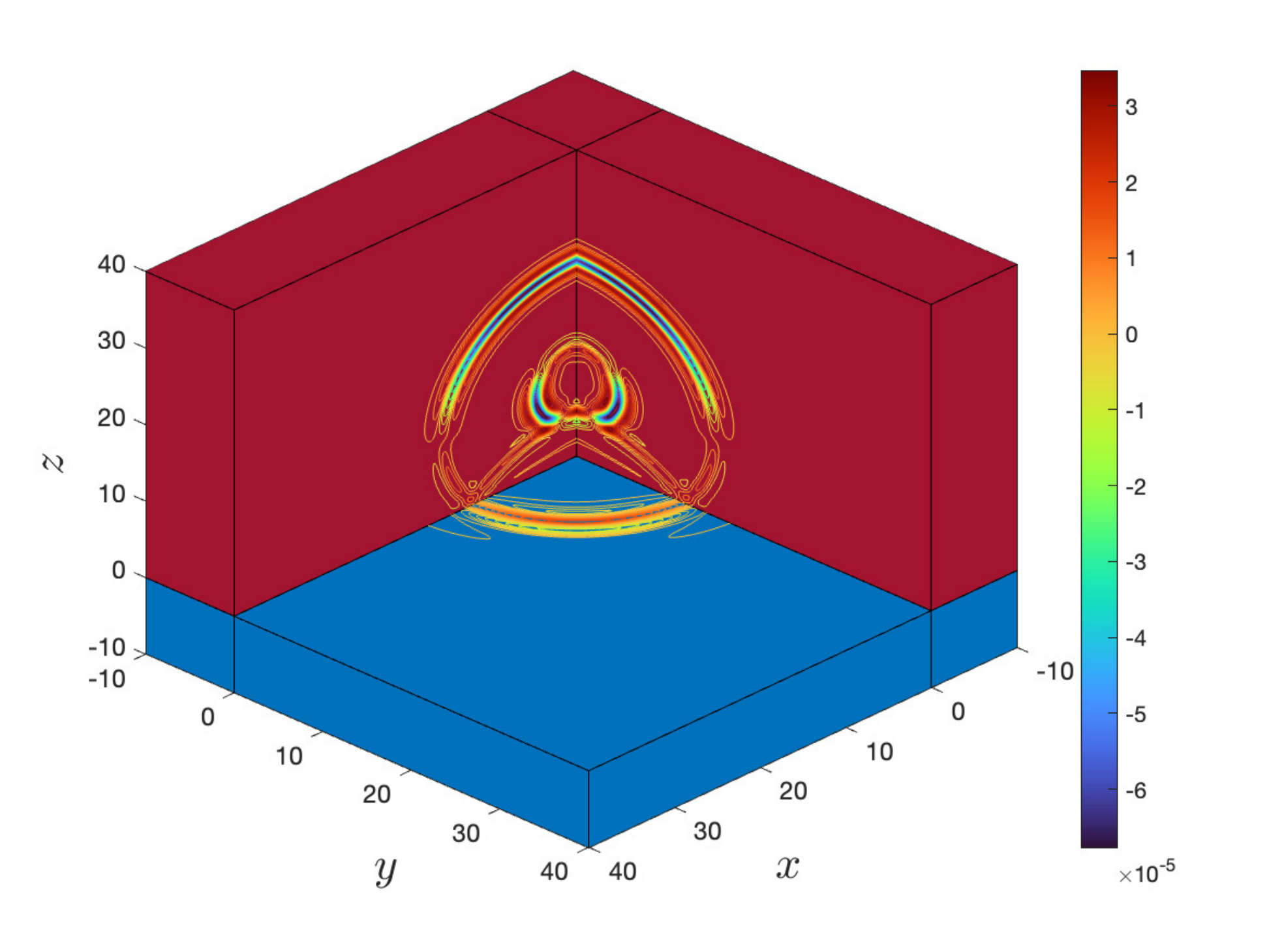}
     \includegraphics[width=0.49\textwidth,height=0.24\textwidth]{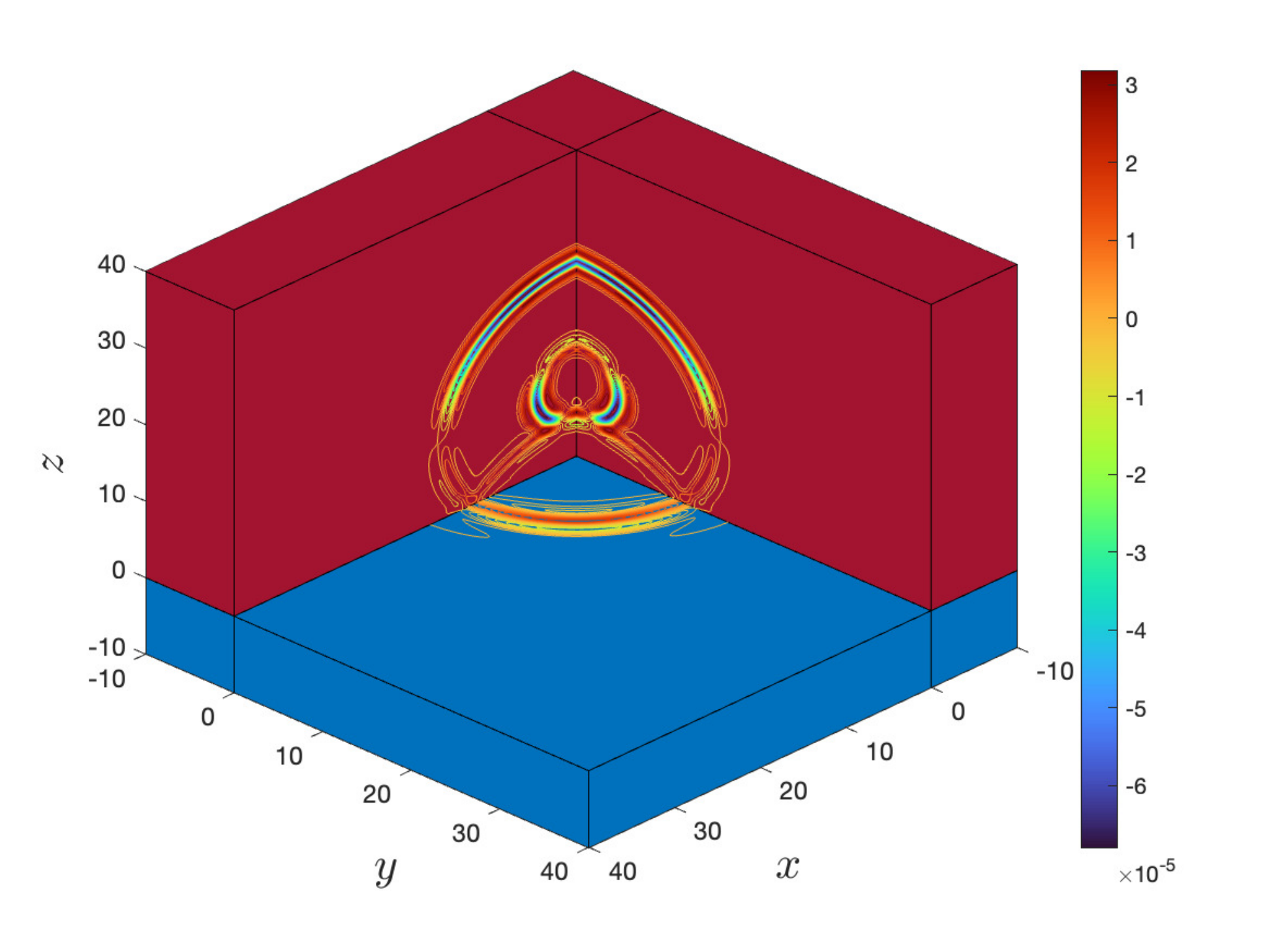}}
     \\
    \centering
    \subfigure[$v_3(x, y, z)$ at $t = 8$s. (left: FSM, $N=512^3$, right: LOSS, $N=513^2\times 1025$)]
    {\includegraphics[width=0.49\textwidth,height=0.24\textwidth]{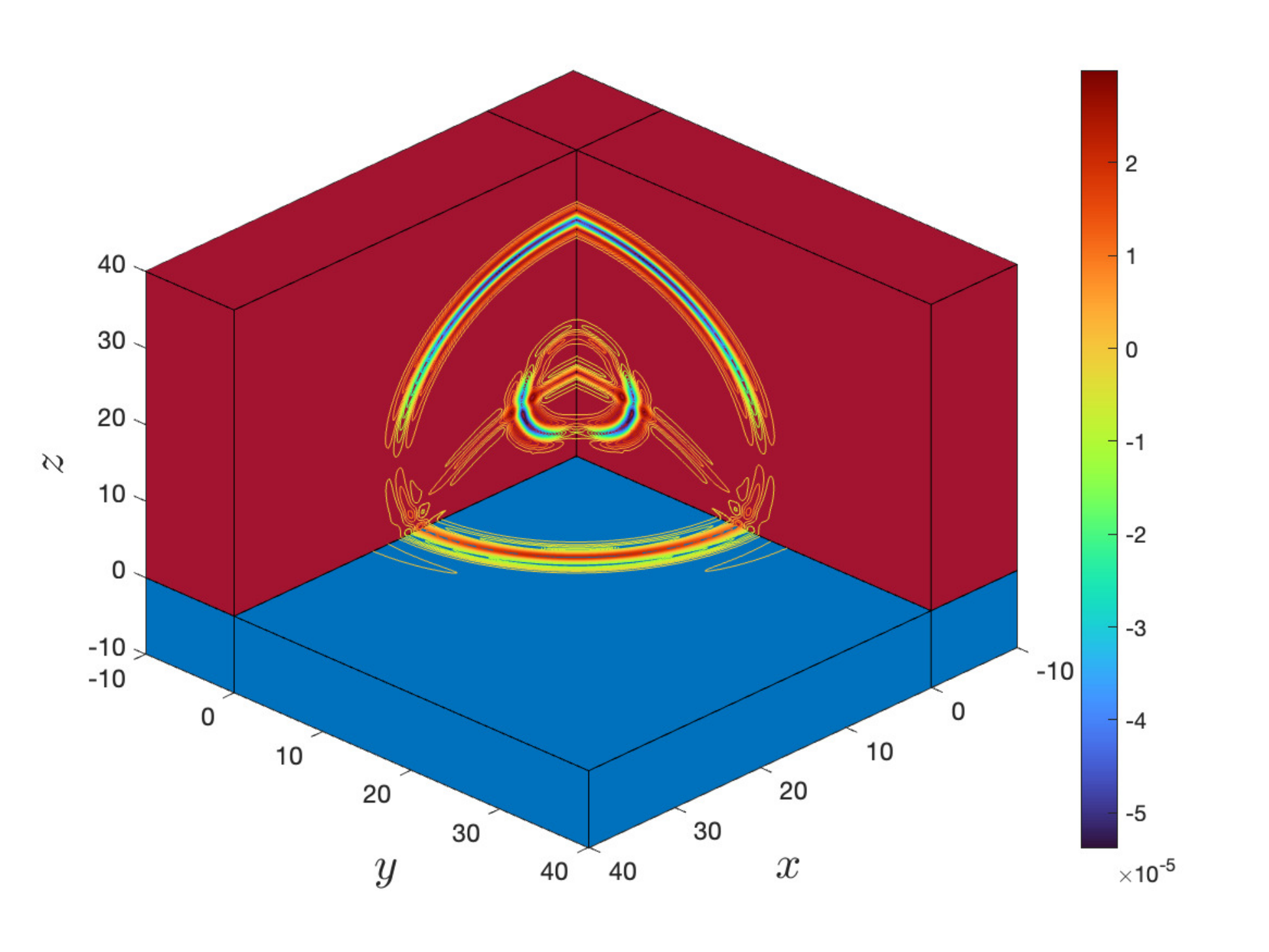}
     \includegraphics[width=0.49\textwidth,height=0.24\textwidth]{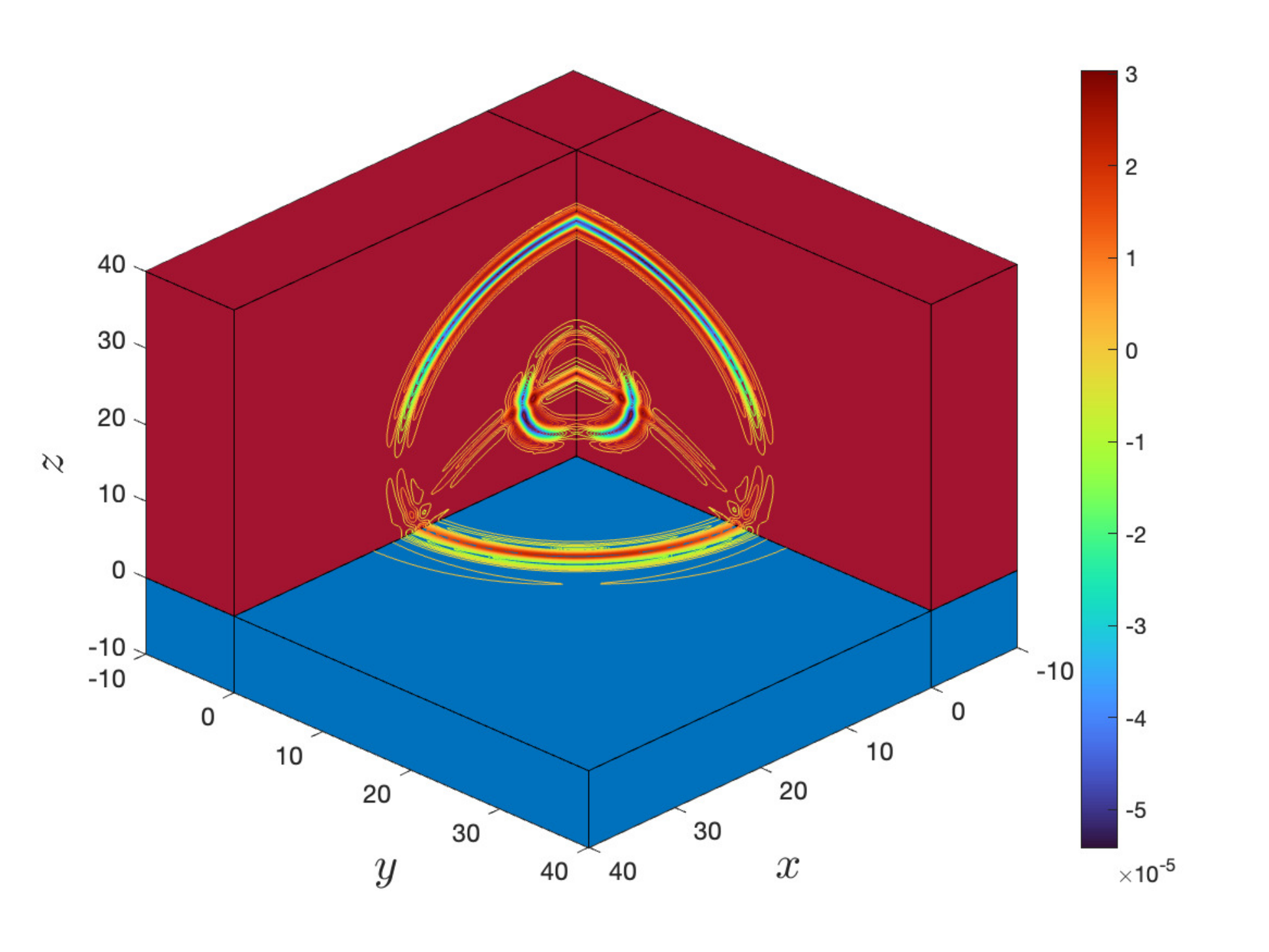}}
     \\
    \centering
    \subfigure[$v_3(x, y, z)$ at $t = 10$s. (left: FSM, $N=512^3$, right: LOSS, $N=513^2\times 1025$)]{\includegraphics[width=0.49\textwidth,height=0.24\textwidth]{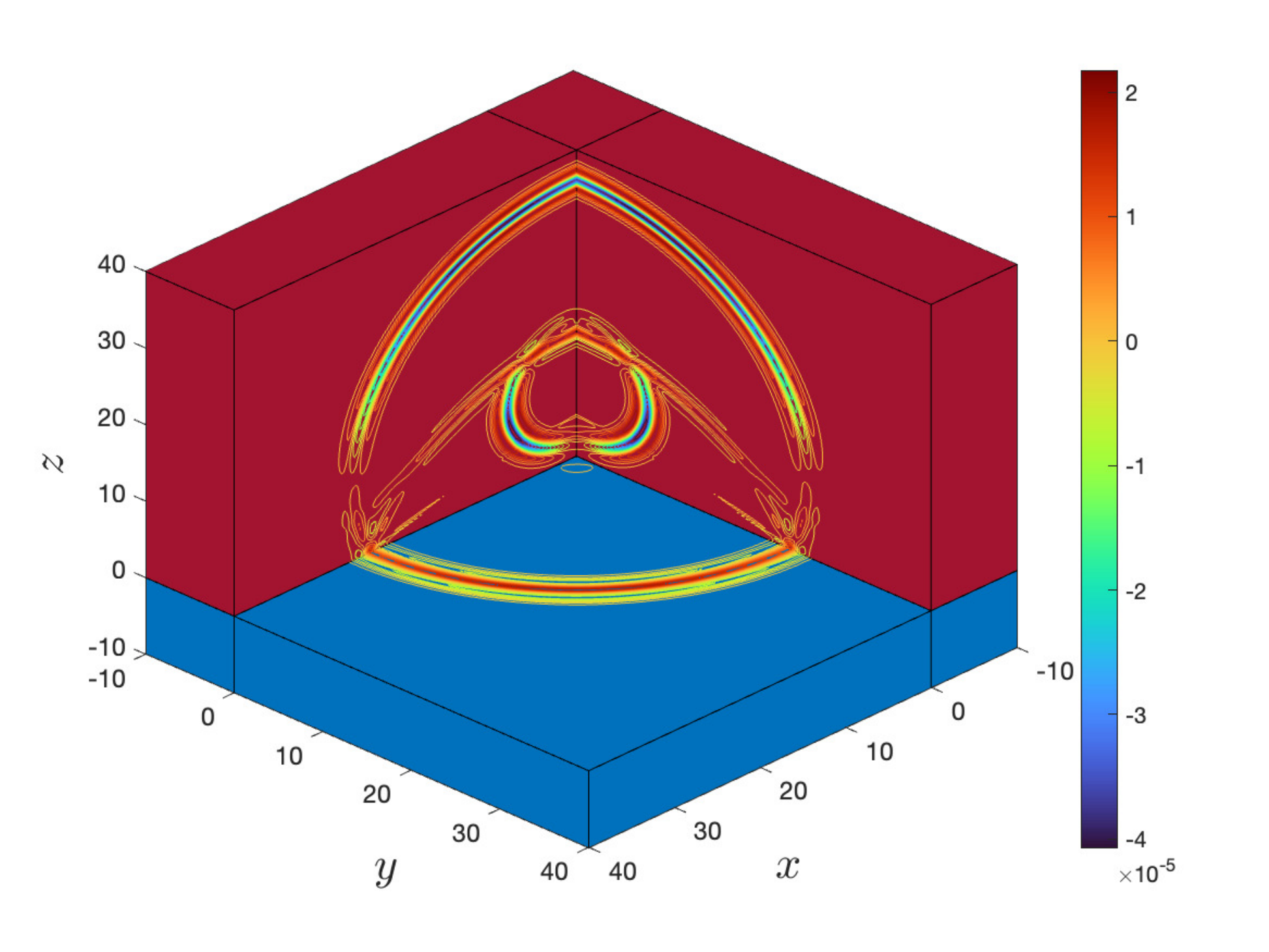}
     \includegraphics[width=0.49\textwidth,height=0.24\textwidth]{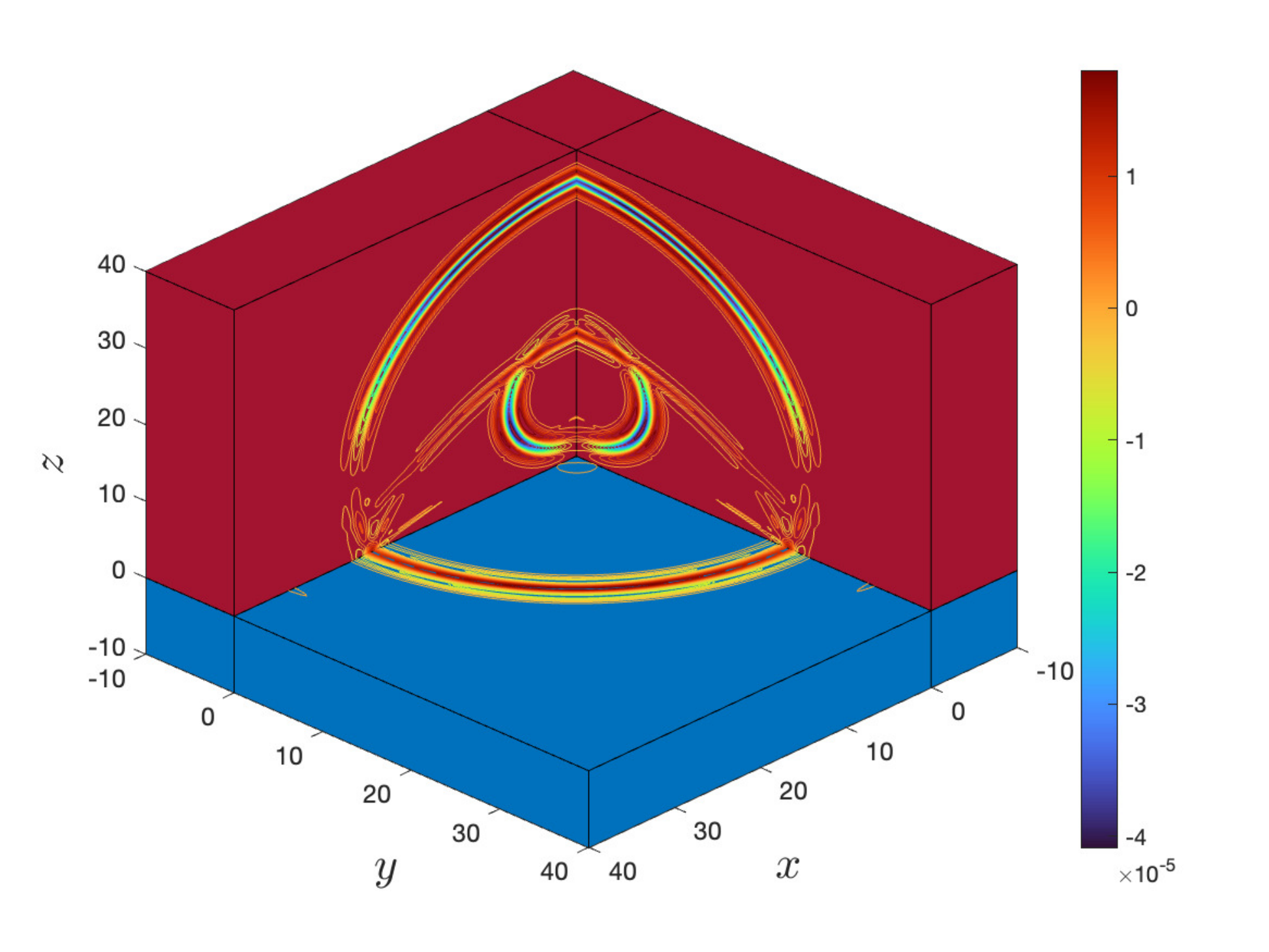}}
    \caption{\small Wave propagation in 3-D double-layer media: A comparison of snapshots of wavefield $v_3(x, y, z)$. \label{hetero_PV}}
\end{figure}

The snapshots of the vibrational wavefield $v_3$ of the $P$- and $S$-wave are given in Figure \ref{hetero_PV}. Unlike the propagation in the homogenous media, the scattering of wave occurs when it enters a different media, resulting in the reflection of wavepackets. From the synthetic data in Figures \ref{v3_hetero_t6} and \ref{v3_hetero_t10}, the  superposition of the impulse and reflected wave is clearly observed.

The convergence of LOSS is still verified in Figure \ref{3d_hetero_convergnce}, albeit with an evident reduction in the order. This is caused by the sharp variation of the model parameters. As visualized in Figures \ref{v3_hetero_t6} and \ref{v3_hetero_t10}, small oscillations in the reflected waves are produced by the spline construction under a coarse grid mesh ($N_x \le 385$). Fortunately, such phenomenon  can be alleviated when a finer grid mesh was adopted.  In particular, we tried to refine the grid mesh in $z$-direction, where there was a sharp variation in density and group velocities, and found that it succeeded in further reducing the errors (see Table \ref{cpu_time}). This manifests that LOSS is still capable to deal with the multi-layer geological model, albeit the grid mesh must be refined near the discontinuity of model parameters.
\begin{table}[!h]
  \caption{\small The memory requirement for storing nine wavefields (in both serial setting and $4\times 4\times 4$ decomposition) and the numerical errors up to $T=10$s (with $\Delta t = 0.005$s, 2000 steps).   \label{cpu_time}}
    \centering
    \resizebox{1\columnwidth}{!}{
    \begin{tabular}{|c|c|c|c|c|c|c|c}
    \hline
    Grid&$129^3$&$257^3$&$513^3$&$641^3$		&$513^2\times1025$		&$769^3$	\\ 
    \hline
    Memory (Serial)	&	$0.14$GB	&	$1.14$GB	&	$9.05$GB	&	$17.66$GB	& 	$18.09$GB &	$30.49$GB	\\
    \hline
    Memory (MPI)	&	$0.17$GB	&	$1.23$GB	&	$9.43$GB	&	$18.25$GB	& 	$18.64$GB&	$31.33$GB	\\
    \hline
    $\varepsilon_2(t=10)$	&	$77.845$	&$26.078$	&$8.320$	&	$5.742$	&	$5.334$	&	$4.985$	\\  
    \hline
    $\varepsilon_\infty(t=10)$	&	$0.785$	&$0.300$		&$0.100$	&	$0.075$	&	$0.083$	&	$0.067$	\\
    \hline
    \end{tabular}
    }
\end{table}

\begin{figure}[!h]
    \centering
   \subfigure[Computational complexity.\label{complexity}]{ \includegraphics[width=0.49\textwidth,height=0.24\textwidth]{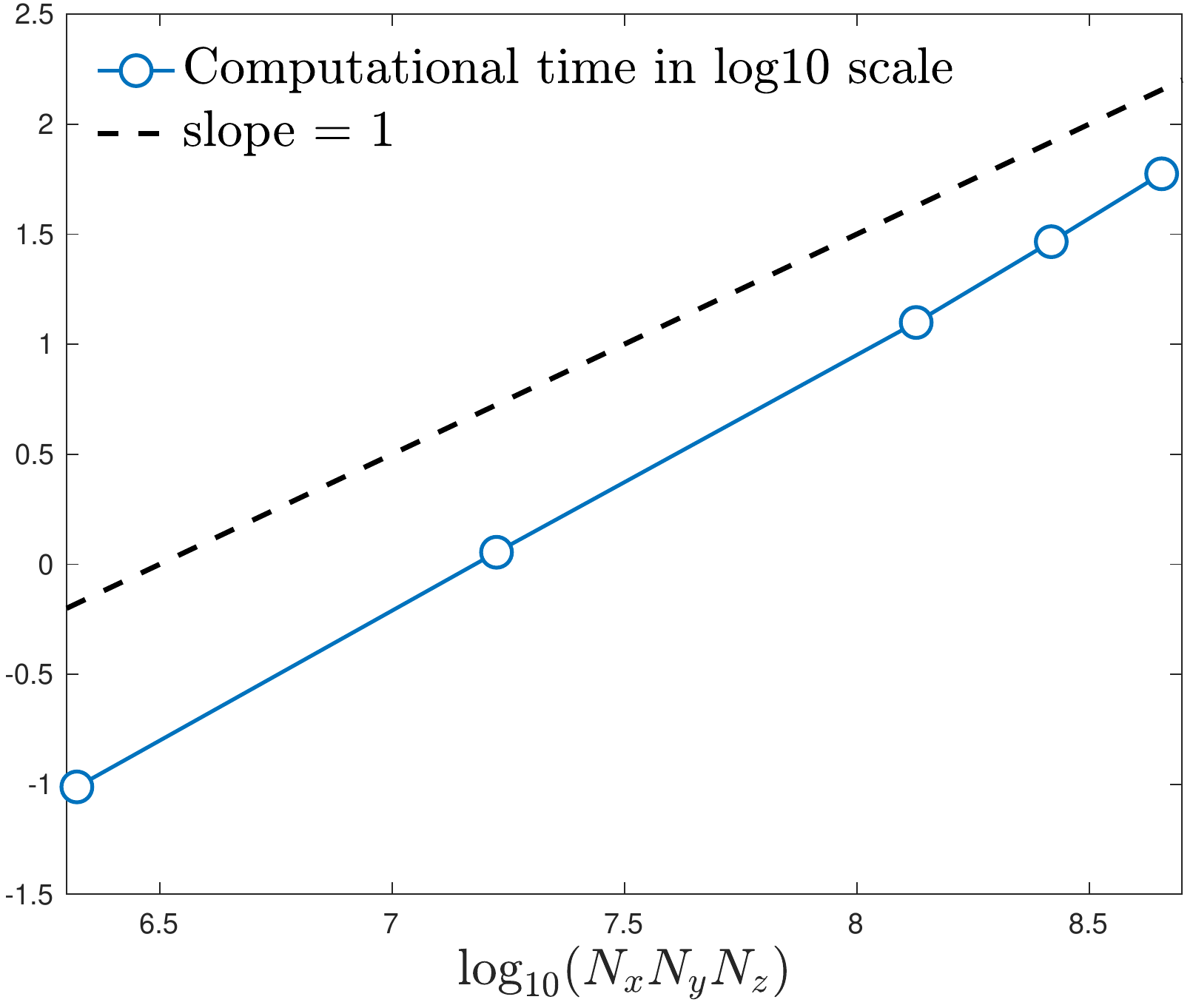}}
   \subfigure[Speedup ratio.\label{Speedup}]{\includegraphics[width=0.49\textwidth,height=0.24\textwidth]{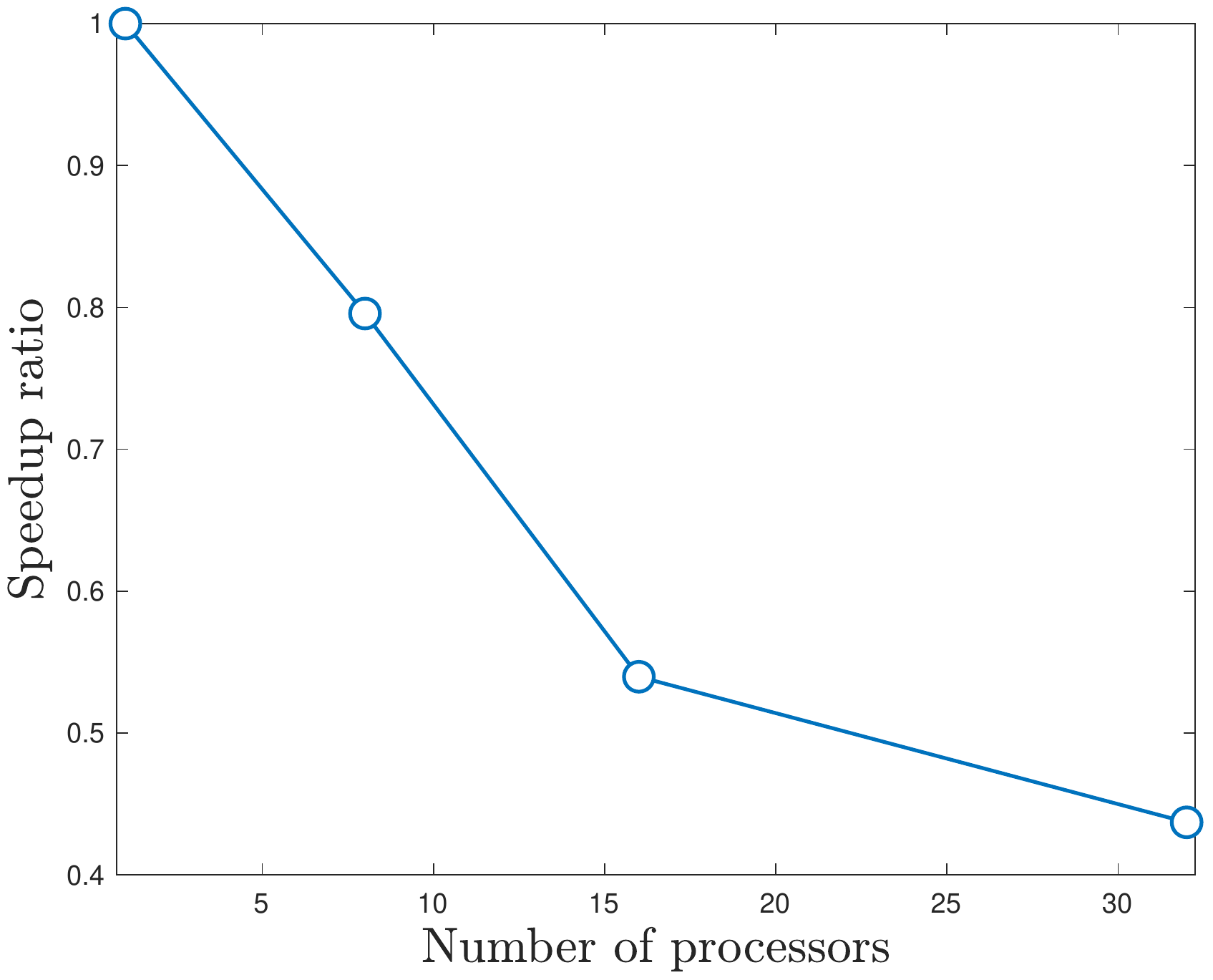}}
    \caption{The complexity of LOSS is almost proportional to $N_x N_y N_z$. In addition, it achieves a speedup ratio about $43.71\%$ using $32$ cores. The platform is AMD Ryzen 7950X (4.50GHz, 64MB Cache, 16 Cores, 32 Threads) with 64GB Memory (4800Mhz). }
\end{figure}

Finally, we tested the complexity of LOSS by recording the computational time under $N = N_x N_y N_z = 129^3, 257^3, 513^3, 641^3, 769^3$. As shown in Figure \ref{complexity}, the complexity of LOSS is almost proportional to the mesh size $N$. Besides, we also calculated the speedup ratio of LOSS in Figure \ref{Speedup} by performing the simulations under the fixed grid mesh $N = 241^3$ and time step $\Delta t = 0.005$s. For serial realization, it required $11.85$ hours to reach the final time $T = 10$s (2000 steps in total). By contrast, when the domain was decomposed into $4\times 4\times 4$ patches and 
$32$ cores ($64$ tasks) were used, it only took $0.85$ hour and the speedup ratio was about $43.71\%$. In fact, the storage of shared grids is relatively small compared with the total memory requirement (see Table \ref{cpu_time}). Thus it is expected that LOSS can achieve even higher speedup ratio when the latency caused by the Hyper-Threading technique is precluded.

%

\section{Conclusions and discussions}
\label{sec.con}

This paper proposes the distributed local spline simulator (LOSS) for solving the elastic wave propagation, where the wavefields are expanded by patched local cubic B-splines and the first-order spatial derivatives can be calculated accurately with low complexity. A perfectly matched boundary condition (PMBC) is introduced by exploiting the exponential decay property of the wavelet basis in its dual space. In this manner, the local spline is able to recover the global spline as accurately as possible with only local communication costs, thereby greatly facilitating  the distributed parallelization. Several typical 2-D and 3-D examples are provided to validate the accuracy, efficiency and parallel scalability of LOSS.

For brevity, we only discuss the implementation under a uniform grid mesh, but  a nonuniform grid is much more desirable when the model parameters have discontinuities or large variations. Actually, the settings of PMBCs can be readily generalized to the non-uniform grid owing to the scaling and flexibility of wavelet basis, and  the implementation of LOSS on a structure-driven grid mesh may be a topic of our future work.

\section*{Acknowledgement}

This research was supported  by the National Natural Science Foundation of China (Nos.~12271303). X. Guo was partially supported by the Natural Science Foundation of Shandong Province for Excellent Youth Scholars (Nos.~ZR2020YQ02) and the Taishan Scholars Program of Shandong Province of China (Nos.~tsqn201909044). 
Y. Xiong was partially supported by the National Natural Science Foundation of China (No.~1210010642) and the Fundamental Research Funds for the Central Universities (Nos.~310421125).






\end{document}